\documentclass[11pt]{article}

\usepackage{latexsym}
\usepackage{setspace}
\usepackage{fancyhdr}
\usepackage{epsfig}
\usepackage{graphicx}
\usepackage{subfigure}
\usepackage{pgfplots}
\usepackage{fancybox}
\usepackage{tikz}
\usepackage{enumitem}
\usepackage{authblk}

\usepackage[utf8]{inputenc}
\usepackage{pgfplots}
\DeclareUnicodeCharacter{2212}{−}
\usepgfplotslibrary{groupplots,dateplot}
\usetikzlibrary{patterns,shapes.arrows}
\pgfplotsset{compat=newest}

\usepackage{amsmath,amsfonts,amsthm}
\usepackage{multirow, tabularx}
\usepackage{color}
\usepackage{rotating,amssymb}
\usepackage{graphics}
\usepackage{wrapfig}

\usepackage{floatflt}
\usepackage{pgfgantt}
\usepackage{amsmath}
\usepackage{amsbsy} 
\usepackage{url}
\usepackage[numbers,square,sort&compress]{natbib}
\usepackage{psfrag}
\usepackage[labelfont=bf,small]{caption}
\usepackage{matlab-prettifier}
\usepackage{algorithm}
\usepackage{algpseudocode}
\usepackage{footnote}
\usepackage{upgreek}

\definecolor{mygrey}{gray}{0.8}
\definecolor{myyellow}{rgb}{.690,.768,.870}

\newcommand{\bU}{\mathbf{U}}
\newcommand{\bY}{\mathbf{Y}}

\newcommand{\bm}[1]{\ensuremath{\mathbf{#1}}}
\newcommand{\bs}[1]{\ensuremath{\boldsymbol{#1}}}
\newcommand{\ol}[1]{\ensuremath{\overline{#1}}}

\def\urltilde{\kern -.15em\lower .7ex\hbox{\~{}}\kern .04em}
\newcommand{\Vhat}{\hat{\bm V}}

\newcommand{\T}{\mathrm{T}}
\newcommand{\D}{\mathrm d}
\newcommand{\rt}{\tilde{r}}
\newcommand{\cb}[1]{{\color{black} #1}}

\newcommand{\Langle}{\mathrel{\tikz[baseline=-0.3ex]{\draw (0,0) -- (0,1ex); \draw (0,0) -- (1ex,0);}}}
\newcommand{\Oangle}{\mathrel{\tikz[baseline=-0.3ex]{\draw (0,0) -- (.75ex,.75ex); \draw (0,0) -- (.75ex,0);}}}

\newcommand{\eat}[1]{}

{\begin{myfont}\begin{center}\begin{Sbox}\begin{minipage}{0.95\textwidth}}%
{\end{minipage}\end{Sbox}%
\colorbox{myyellow}{\TheSbox}\end{center}\end{myfont}}
\newenvironment{short-item}%
{\begin{list}{$\bullet$}{
\setlength{\leftmargin}{\labelwidth}%
}}%
{\end{list}}

\newcommand{\squishlist}{%
   \begin{list}{$\bullet$}
    { \setlength{\itemsep}{0pt}      \setlength{\parsep}{3pt}
      \setlength{\topsep}{3pt}       \setlength{\partopsep}{0pt}
      \setlength{\leftmargin}{1.5em} \setlength{\labelwidth}{1em}
      \setlength{\labelsep}{0.5em} } }

\newcommand{\squishlisttwo}{%
   \begin{list}{$\bullet$}
    { \setlength{\itemsep}{0pt}    \setlength{\parsep}{0pt}
      \setlength{\topsep}{0pt}     \setlength{\partopsep}{0pt}
      \setlength{\leftmargin}{2em} \setlength{\labelwidth}{1.5em}
      \setlength{\labelsep}{0.5em} } }

\newcommand{\squishend}{%
    \end{list}  }

\newcommand{\be}{\begin{equation}}
\newcommand{\ee}{\end{equation}}

\newtheorem{theorem}{Theorem}

\newtheorem{definition}{Definition}

\newtheorem{remark}{Remark}

\usepackage{listings}

\usepackage[font=scriptsize]{caption}

\begin{document}

\title{Oblique projection for  scalable rank-adaptive reduced-order modeling of nonlinear stochastic PDEs  with time-dependent bases}

\author[1]{M. Donello}
\author[1]{G. Palkar}
\author[1]{M.H. Naderi}
\author[2]{D. C. Del Rey Fern\'andez}
\author[1]{Hessam Babaee\thanks{Corresponding author. Email:h.babaee@pitt.edu.}}
\affil[1]{\footnotesize Department of Mechanical Engineering and Materials Science, University of Pittsburgh}
\affil[2]{\footnotesize Department of Applied Mathematics, University of Waterloo}

\date{}
\maketitle

\begin{abstract}
Time-dependent basis reduced order models (TDB ROMs) have successfully been used for approximating the solution to nonlinear stochastic partial differential equations (PDEs). For many practical problems of interest, discretizing these PDEs results in massive matrix differential equations (MDEs) that are too expensive to solve using conventional methods. While TDB ROMs have the potential to significantly reduce this computational burden, they still suffer from the following challenges: (i) inefficient for general nonlinearities, (ii) intrusive implementation, (iii) ill-conditioned in the presence of small singular values, and (iv) error accumulation due to fixed rank. To this end, we present a scalable method \cb{based on oblique projections} for solving TDB ROMs that is computationally efficient, minimally intrusive, robust in the presence of small singular values, rank-adaptive, and highly parallelizable. These favorable properties are achieved via low-rank approximation of the time discrete MDE. Using the discrete empirical interpolation method (DEIM), a low-rank decomposition is computed at each iteration of the time stepping scheme, enabling a near-optimal approximation at a fraction of the cost. \cb{We coin the new approach TDB-CUR since it is equivalent to a CUR decomposition based on sparse row and column samples of the MDE.}  We also propose a rank-adaptive procedure to control the error on-the-fly. Numerical results demonstrate the accuracy, efficiency, and robustness of the new method for a diverse set of problems.

\end{abstract}

\section{Introduction}\label{sec:Intro}
 Discretizations of many time-dependent partial differential equations (PDEs) result in  matrix differential equations (MDEs) in the form of $\D\bm{V}/\D t = \mathcal{F}(\bm{V})$, where $\bm{V} \in \mathbb{R}^{n\times s}$ is the solution matrix and $\mathcal{F}(\bm{V}) \in \mathbb{R}^{n\times s}$ is obtained by discretizing the PDE in all dimensions except time. One such example is  the uncertainty propagation of random parameters into the PDEs, which requires solving the PDEs for a large number of random realizations \cite{SL09,KL07}. Discretization of this problem can be formulated as an MDE, where the rows of the matrix are obtained by discretizing the PDE in the physical domain and the columns of the matrix are samples of the discretized equation for a particular choice of random parameters. For high-dimensional PDEs subject to high-dimensional random parameters, the resulting MDEs can be massive. For example, uncertainty quantification of a 3D time-dependent fluid flow typically requires solving an MDE with $n \sim \mathcal{O}(10^6)-\mathcal{O}(10^9)$ grid points (rows) and $s \sim \mathcal{O}(10^4)-\mathcal{O}(10^7)$ random samples (columns).  Therefore, the solution to these massive MDEs is cost prohibitive due to the floating point operations (flops), memory, and storage requirements.  The discretization of many other PDEs can also be cast as MDEs, for example, \cb{kinetics equations \cite{EL18,HW22,KS23}}, linear sensitivity analyses \cite{DCB22} and species transport equations in turbulent combustion \cite{RNB21}.

For many practical applications, $\bm{V}(t)$ is instantaneously low-rank. Therefore, low-rank approximations  using time-dependent bases (TDBs) have the  potential to significantly reduce the computational cost of solving massive MDEs. For these systems, a TDB-based low-rank approximation  extracts low-rank structures via TDBs for the column and row spaces of $\bm{V}$. A reduced-order model (ROM) is then constructed by projecting the full-order model (FOM) onto the column and row TDBs. Low-rank approximation based on TDB was first introduced in the quantum chemistry field to solve the Schr\"{o}dinger equation  \cite{Beck:2000aa}, and it is commonly known as the multiconfiguration time-dependent Hartree (MCTDH) method. The MCTDH methodology was later presented for generic MDEs in \cite{KL07} and is referred to as dynamical low-rank approximation (DLRA).

Various TDB ROM schemes have also been developed  to solve  stochastic partial differential equations (SPDEs). Dynamically orthogonal (DO) decomposition \cite{SL09}, bi-orthogonal (BO) decomposition \cite{CHZI13}, dual dynamically orthogonal (DDO) decomposition \cite{MN18}, and dynamically bi-orthogonal decomposition (DBO) \cite{PB20} are all TDB-based low-rank approximation techniques for solving stochastic PDEs (SPDEs). In all of these decompositions (DO, BO, DDO, and DBO), an evolution equation for the mean field is developed, along with evolution equations for the TDB-ROM of the mean-subtracted stochastic fields. Although these decompositions have different forms and constraints, they are all equivalent, i.e., they produce identical low-rank matrices  \cite{CSK14,PB20}, and their differences lie only in their numerical performance. TDB ROMs have also been used in other fields and applications including dynamical systems \cite{Dieci:2006aa,doi:10.1137/19M1257275}, combustion \cite{NBGCL21,RNB21}, linear sensitivity analysis \cite{DCB22}, \cb{dynamical instabilities \cite{Babaee_PRSA, BS19, blanchard2019learning}}, deep learning \cite{SZK23}, and singular vale decomposition (SVD) estimation for matrices that vary smoothly with a parameter \cite{Wright:1992aa}.

Despite the potential of using TDB ROMs to significantly reduce the computational cost of solving massive MDEs, there are still a number of outstanding challenges for most practical problems of interest. We  summarize three key challenges below:
\vspace{-.1cm}
\begin{enumerate}[label=(\roman*)]
    \item \textbf{Computational efficiency:} For specific classes of equations (e.g. homogeneous linear and quadratic nonlinear),  rank-$r$ TDB ROMs can be solved efficiently with operations that scale with $\mathcal{O}(nr)$ and $\mathcal{O}(sr)$ for linear MDEs or  scale with $\mathcal{O}(nr^2)$ and $\mathcal{O}(sr^2)$ for quadratic MDEs. However, this computational efficiency is lost for general nonlinearities, requiring operations that scale with the size of the FOM, i.e., $\mathcal{O}(ns)$. 
    \item \textbf{Intrusiveness:} Even in the  special cases of homogeneous linear and quadratic nonlinear equations, efficient implementation of TDB ROM evolution equations is an intrusive process   \cite[Appendix B]{NB23}. This involves replacing the low-rank approximation in the FOM, projecting the resulting equation onto the tangent manifold, and obtaining low-rank matrices for each term on the right-hand side. The process  requires significant effort to derive, implement, and debug the code. This poses a major obstacle for most practitioners, creating a significant barrier to adopting the methodology.
    \item \textbf{Ill-conditioning:} The TDB ROM evolution equations become numerically unstable when the singular values of the low-rank approximation  become very small. This is particularly problematic because it is often necessary to retain very small singular values in order to have an accurate approximation. Small singular values lead to ill-conditioned matrices that require inversion in all variations of TDB ROM evolution equations \cite{KL07,SL09,CHZI13,MN18,PB20}, resulting in restrictive time step limitations for numerical integration and error amplification.
\end{enumerate}

Although some of these challenges have been tackled, there is currently no methodology that can address all of them. In particular, the problems of ill-conditioning and computational expense must be resolved for practitioners to adopt TDB-based low-rank approximations for MDEs.
To address the issue of ill-conditioning, a projector-splitting time integration was proposed  \cite{LO14}, in which arbitrarily small  singular values can be retained. However, this scheme includes a  backward time integration substep, which is an unstable substep for dissipative problems. To address this issue, an unconventional robust integrator was recently proposed \cite{ceruti2021unconventional}  which retains the robustness with respect to small singular values while avoiding the unstable backward step. The authors also presented an elegant rank adaptive strategy, where the rank of the approximation changes over time to maintain a desired level of accuracy. Despite these advantages, this scheme is first-order in time \cite[Theorem 4]{ceruti2021unconventional}. In \cite{Babaee:2017aa}, a pseudo-inverse methodology was presented  as a remedy to maintain a well-conditioned system. However, in this approach, it is difficult to determine what singular value threshold must be used. Another projection method was presented in \cite{kieri2019projection} that retains robustness with respect to small singular values and can be extended to high-order explicit time discretizations. 

The three time-integration schemes  presented in \cite{LO14,ceruti2021unconventional,kieri2019projection} and the pseudo-inverse methodology presented in \cite{Babaee:2017aa} can retain $\mathcal{O}(n+s)$ cost for linear and quadratic MDEs.  But achieving this speedup comes at the expense  of a highly intrusive implementation. However, for generic nonlinear MDEs, an intrusive implementation is not possible, and  the computational cost of solving the TDB ROMs using methods presented in \cite{LO14,ceruti2021unconventional,kieri2019projection,Babaee:2017aa} scales with $\mathcal{O}(ns)$, which is the same as the cost of solving the FOM. Recently, a sparse interpolation algorithm was presented for solving the TDB ROM evolution equations with a computational complexity that scales with $\mathcal{O}(n+s)$ for  generic nonlinear SPDEs \cite{NB23}.  However, this methodology still lacks robustness when the singular values become small, as it requires the inversion of the matrix of singular values.

In this work, we present a methodology inspired by \cb{interpolation  and hyper-reduction techniques developed to accelerate nonlinear ROMs and finite-element models in vector differential equations \cite{BMNP04,R05,AF08,CS10,FACC14,doi:10.1137/140989169}. In particular,} we present  CUR factorizations of low-rank matrices that address the above challenges, i.e., (i) the computational cost of the methodology scales with $\mathcal{O}(n+s)$ for generic nonlinear SPDEs both in terms of flops and memory costs, (ii) it lends itself to simple implementation in existing codes, and  (iii)  the time-integration is robust in the presence of small  singular values, and high-order explicit time integration can be used.  To this end, the main elements  of the presented methodology are (i) a time-discrete variational principle for minimization of the residual due to low-rank approximation error, and (ii) a CUR factorization based on strategic row and column sampling of the time discrete MDE. 


The remainder of the paper is organized as follows: In Section \ref{sec:methodology}, we first review the time continuous variational principle and its associated challenges. We then present the time discrete variational principle along with the rank-adaptive sparse sampling strategy for solving TDB-ROMs. Finally, we show that the resulting low-rank approximation is equivalent to a CUR factorization and we provide an upper bound on the approximation error. In Section \ref{sec:demonstration}, we demonstrate the method for a toy problem as well as the stochastic Burgers equation and stochastic nonlinear advection-diffusion-reaction equation. In Section \ref{sec:conclusion}, we summarize the present work and discuss its implications.

\section{Methodology}\label{sec:methodology}
\subsection{Setup}
Consider the nonlinear stochastic PDE given by:
\begin{equation}\label{eq:FOM_Cont}
\frac{\partial v}{\partial t} = f(v;x,t,\bs\xi),
\end{equation}
augmented with appropriate initial and boundary conditions. In the above equation, $v=v(x,t;\bs\xi)$, $x$ is the spatial coordinate, $\bs\xi \in \mathbb{R}^d$ are the set of random parameters, $t$ is time, and $f(v;x,t,\bs\xi)$ includes the nonlinear spatial differential operators. We assume generic nonlinear PDEs, where the nonlinearity of $f$ versus $v$ may be non-polynomial, e.g., exponential, fractional, etc.   For the sake of simplicity in the exposition, we consider a collocation/strong-form discretization of  Eq. \ref{eq:FOM_Cont} in $x$ and $\bs\xi$. Because of the simplicity of the resulting discrete system, this choice  facilitates an uncluttered illustration of the main contribution of this paper, which is focused on the efficient low-rank approximation of nonlinear matrix differential equations. However, the presented methodology can also be applied to other types of discretizations, for example, weak form discretizations (finite element, etc). Examples of collocation/strong-form discretizations in the spatial domain are Fourier/polynomial spectral collocation schemes or finite-difference discretizations. Example collocation schemes in the random domain include the probabilistic collocation method (PCM) \cite{xiu2006high} or any Monte-Carlo-type sampling methods \cite{giles2008multilevel, barth2011multi, kuo2012quasi}.  Applying any of the above schemes to  Eq. \ref{eq:FOM_Cont}  leads to the following \emph{nonlinear matrix differential equation}:
\begin{equation}\label{eq:FOM}
\frac{\mathrm d \bm V}{\mathrm d t} = \mathcal F(t,\bm V), \quad   t\in I=[0,T_f],
\end{equation}
where $I=[0,T_f]$ denotes the time interval,  $\bm V(t): I \rightarrow \mathbb{R}^{n \times s}$ is a matrix with $n$ rows corresponding to collocation points in the spatial domain and $s$ columns corresponding to collocation/sampling points of the parameters $\bs\xi$, and  $\mathcal F(t,\bm V): I \times \mathbb{R}^{n\times s} \rightarrow \mathbb{R}^{n\times s}$ is obtained by  discretizing $f(v;x,t,\bs\xi)$ in $x$ and $\bs\xi$. Eq. \ref{eq:FOM} is augmented with appropriate initial conditions, i.e., $\bm{V}(t_0)=\bm{V}_0$. We also assume that boundary conditions are already incorporated into  Eq. \ref{eq:FOM}, which can be accomplished in a number of ways, for example by using weak treatment of the boundary conditions \cite{patil2023reduced}.   

For the remainder of this paper, we will refer to Eq. \ref{eq:FOM} as the FOM, which will be used as the ground truth for evaluating the performance of the proposed methodology. For the problems targeted in this work, we assume $n>s$ without loss of generality.

 \cb{The presented methodology is limited to explicit time integration schemes. For the computational complexity analysis, we consider sparse discretization schemes for spatial discretization, which means that each row is dependent on $p_a$ rows, where $p_a << n$. The majority of discretization schemes, e.g., finite difference, finite volume, finite element, spectral element, result in sparse row dependence. As a result, the computational cost of computing each column of FOM (Eq. \ref{eq:FOM}) is $\mathcal{O}(n)$ and the cost of solving  MDE \ref{eq:FOM} for all $s$ columns scales with $\mathcal{O}(ns)$. We also note that the presented methodology is \emph{not} limited to sparse spatial discretizations and can be applied to dense discretizations as well. See Remark \ref{rm:dense_MDE} for more details.}

\subsection{Preliminaries}
\cb{
 In this section, we present some of the definitions of matrix manifolds, tangent spaces,  orthogonal and oblique projections, and CUR decomposition. 

 \begin{definition}[Low-rank matrix manifolds]\label{def:Mr}
 The low-rank matrix manifold $\mathcal{M}_r$ is defined as the set
\begin{equation*}
\mathcal{M}_r = \{\hat{\bm V} \in \mathbb{R}^{n \times s}: \ \mbox{rank}(\hat{\bm V}) = r \}, 
\end{equation*}
of matrices of fixed rank $r$. Any member of the set $\mathcal{M}_r$ is denoted by a hat symbol $( \hat{ \ \ } )$, e.g., $\hat{ \bm V}$.
\end{definition}
Any member of $\mathcal{M}_r$ may be represented by $\hat {\bm V} = \bm{U} \bs \Sigma \bs{Y}^T$, where $\bm U \in \mathbb{R}^{n \times r}$ and       $\bm Y \in \mathbb{R}^{s \times r}$ are a set of orthonormal columns and $\bs \Sigma \in \mathbb{R}^{r \times r}$ is a rank-$r$ matrix. The rank-$r$ matrix $\hat{\bm V}$ may also be represented via the multiplication of two matrices, i.e.,  $\hat {\bm V} = \bm{U} \bs{Y}^T$, where   $\bm U \in \mathbb{R}^{n \times r}$ and       $\bm Y \in \mathbb{R}^{s \times r}$ have full column rank. 

 \begin{definition}[Tangent space]\label{def:Tan_Spc}
 The tangent space of manifold $\mathcal{M}_r$ at $\hat{\bm V}$, represented with the decomposition of $\hat{\bm V} = \bm{U} \bs \Sigma \bm Y^T$, is  the set of matrices in the form of \cite{KL07}:
\begin{equation*}
\mathcal{T}_{\Vhat} \mathcal{M}_r = \{\delta \bm U \bs \Sigma \bm Y^T +  \bm U \delta\bs \Sigma \bm Y^T +   \bm U \bs \Sigma \delta \bm Y^T: \  \delta \bm U^T \bm U = \bm 0 \ \mbox{and} \  \delta \bm Y^T \bm Y = \bm 0 \},
\end{equation*}
where $\delta \bm U \in \mathbb{R}^{n \times r}$ and $\delta \bm Y \in \mathbb{R}^{s \times r}$.
\end{definition}

\begin{definition}[Orthogonal projection onto the tangent space]\label{def:ortho_Tan_Spc}
 The orthogonal projection of matrix $\bm W \in \mathbb{R}^{n\times s}$ onto the  tangent space of manifold $\mathcal{M}_r$ at $\hat{\bm V}$, represented with the decomposition of $\hat{\bm V} = \bm{U} \bs \Sigma \bm Y^T$, is  given by \cite[Lemma 4.1]{KL07}:
\begin{equation}\label{eq:tan_spc}
    \mathcal{P}_{\mathcal{T}_{\hat{\bm{V}}}} (\bm W) = \bm{U}\bm{U}^T \bm W + \bm W \bm{Y}\bm{Y}^T -  \bm{U}\bm{U}^T \bm W \bm{Y}\bm{Y}^T.
\end{equation}
\end{definition}

In the above projection, $\bm U \bm U^T$ and $\bm Y \bm Y^T$ are orthogonal projections onto spaces spanned by the columns of $\bm U$ and $\bm Y$. We denote these orthogonal projections with:
\begin{equation}
    \bm{P}^{\Langle}_{\bm{U}} = \bm U \bm U^T \quad \mbox{and} \quad  \bm{P}^{\Langle}_{\bm{Y}} = \bm Y \bm Y^T,
\end{equation}
where the symbol $\Langle$ indicates  orthogonal projection.  In the following, we define oblique projectors. We first explain the notation that is  used in this section. Let $\bm U \in \mathbb{R}^{n\times r}$ and $\bm Y \in \mathbb{R}^{s\times r}$ be matrices whose columns are orthonormal and let $\bm p =[p_1,p_2, \dots, p_{r'}]\in\mathbb N^{r'}$ and $\bm s=[s_1,s_2, \dots,  s_{r'}]\in\mathbb N^{r'}$ be vectors containing row and column indices, where the number of indices can be greater than or equal to the dimension of the subspaces spanned by $\bm U$ and $\bm Y $, i.e., $r' \geq r$. Also, $r' \leq n$ for row indices and $r' \leq s$ for column indices.  We use MATLAB indexing where $\bm V(\bm p,:) \in \mathbb{R}^{r' \times s}$ selects all columns at the $\bm p$ rows, and $\bm V(:,\bm s)\in \mathbb{R}^{n \times r'}$ selects all rows at the $\bm s$ columns of the matrix $\bm V$. We also use the indexing matrices, $\bm P = \bm{I}_n(:,\bm p) \in \mathbb{R}^{n \times r'}$ and $\bm S = \bm{I}_s(:,\bm s) \in \mathbb{R}^{s \times r'}$, where $\bm I_n$ and $\bm I_s$ are identity matrices of size $n \times n$ and  $s \times s$, respectively. It is easy to verify that $\bm P^T \bm U \equiv \bm{U}(\bm p,:)$ and $\bm S^T \bm Y \equiv \bm{Y}(\bm s,:)$. Let $( \ )^{\dagger}$ denote the Moore–Penrose pseudoinverse of a matrix, i.e.,  $\bm A^{\dagger} = (\bm A^T \bm A)^{-1} \bm A^T$.  

 \begin{definition}[Oblique projection]\label{def:obl_proj}
Let $\bm {U} \in \mathbb{R}^{n \times r}$ and  $\bm {Y} \in \mathbb{R}^{s \times r}$ be  orthonormal matrices and let $\bm p \in\mathbb N^{r'}$ and $\bm s \in\mathbb N^{r'}$ be  sets of distinct row and column indices, respectively. Oblique projectors  onto Ran($\bm{U}$) and  Ran($\bm{Y}$) are defined as \cite{sorensen2016deim}
\begin{equation}\label{eq:interp_proj}
\bm{P}^{\Oangle}_{\bm{U}} = \bm U (\bm P^T \bm U)^{\dagger} \bm P^T \quad \mbox{and} \quad \bm{P}^{\Oangle}_{\bm{Y}} = \bm S (\bm Y^T \bm S)^{\dagger} \bm Y^T,
\end{equation}
provided $(\bm P^T \bm U)^T(\bm P^T \bm U) \in \mathbb{R}^{r \times r}$ and $(\bm Y^T \bm S)^T (\bm Y^T \bm S) \in \mathbb{R}^{r \times r}$ are invertible. 
 \end{definition}
 In the above definition, the symbol $\Oangle$ distinguishes these projectors from orthogonal projectors.  For a given matrix $\bm A \in \mathbb{R}^{n\times s}$,  $\bm{P}^{\Oangle}_{\bm{U}}$ operates on the left side of the matrix  and  $\bm{P}^{\Oangle}_{\bm{Y}}$ operates on the right side of the matrix. It is easy to verify that $\big (\bm{P}^{\Oangle}_{\bm{U}}\big )^2 =\bm{P}^{\Oangle}_{\bm{U}}$ and $\big (\bm{P}^{\Oangle}_{\bm{Y}}\big )^2 =\bm{P}^{\Oangle}_{\bm{Y}}$.  It is also easy to verify that the oblique projection of $\bm A$ belongs to the manifold of rank-$r$ matrices, i.e. $\bm{P}^{\Oangle}_{\bm{U}} \bm A \bm{P}^{\Oangle}_{\bm{Y}} \in \mathcal{M}_r$.
 
 For the special case of $r'=r$, the oblique projectors $\bm{P}^{\Oangle}_{\bm{U}}$ and $\bm{P}^{\Oangle}_{\bm{Y}}$  are also \emph{interpolatory} projectors. In this case, the oblique projectors become $\bm{P}^{\Oangle}_{\bm{U}} = \bm U (\bm P^T \bm U)^{-1} \bm P^T$ and $\quad \bm{P}^{\Oangle}_{\bm{Y}} = \bm S (\bm Y^T \bm S)^{-1} \bm Y^T$. Unlike orthogonal projection or a general oblique projection, the interpolatory projection is guaranteed to match the original matrix at the selected rows and columns, i.e.,
\begin{equation*}
    \bm P^T\bm{P}^{\Oangle}_{\bm{U}} \bm A = \bm A(\bm p,:) \quad \mbox{and} \quad \bm A \bm{P}^{\Oangle}_{\bm{Y}}\bm S = \bm A(:, \bm s).
\end{equation*}
The other extreme is when all the rows or columns are selected, i.e., $r'=n$ or $r'=s$. Take for example, the projector $\bm{P}^{\Oangle}_{\bm{U}}$ when $r'=n$. In this case, $\bm{P}^{\Oangle}_{\bm{U}}$ becomes the same as the orthogonal projector, i.e., $\bm{P}^{\Oangle}_{\bm{U}} \equiv \bm{P}^{\Langle}_{\bm{U}}$. To show this, first  note that $\bm{P}^{\Oangle}_{\bm{U}}$ is invariant with respect to the ordering of the row indices ($\bm p$) and when $r'=n$, $\bm p$ can be taken to be: $\bm p=[1,2, \dots, n]$. In this case, $\bm P \equiv \bm I_n$. Therefore:
\begin{equation*}
    \bm{P}^{\Oangle}_{\bm{U}} = \bm U (\bm P^T \bm U)^{\dagger} \bm P^T = \bm U (\bm U^T \bm U)^{-1} \bm U^T =  \bm U \bm U^T = \bm{P}^{\Langle}_{\bm{U}},
\end{equation*}
where we have used the orthonormality condition of  $\bm U$, i.e., $\bm U^T \bm U= \bm I_r$, where $\bm{I}_r$ is the $r\times r$ identity matrix.
The analogous relationship exists for $\bm{P}^{\Oangle}_{\bm{Y}}$, when $r'=s$. 
 In the following, we define CUR decompositions, which are closely related to the oblique projections.

\begin{definition}[CUR decomposition]\label{def:CUR}
   A CUR decomposition of matrix $\bm V$ is a rank-$r$ approximation of $\bm V$ in the form of $\bm V \approx \bm{C}\bm{U}\bm{R}$, where $\bm{C} \in \mathbb{R}^{n \times r}$ and $\bm{R}\in \mathbb{R}^{r \times s}$ are actual columns and rows of matrix $\bm V$, i.e., $\bm C = \bm V(:, \bm s)$ and  $\bm R = \bm V(\bm p, :)$.  The matrix $\bm U \in \mathbb{R}^{r \times r}$ is computed such that $\bm C \bm U \bm R$ is a good approximation to $\bm V$.   The CUR of matrix $\bm V$ is denoted with $\mbox{\texttt{CUR}}(\bm V)$.
\end{definition}
Here, the matrices, $\bm C$, $\bm U$, and $\bm R$ are different from the matrices defined in previous sections. Different CUR decompositions can be obtained for the same matrix depending on two factors: (i) the selection of columns and rows, and (ii) the method used to compute the matrix $\bm U$.  For more details on CUR decompositions, we refer the reader to \cite{mahoney2009cur}. Finally, it is easy to verify that  $\mbox{\texttt{CUR}}(\bm V) \in \mathcal{M}_r$. The connection between CUR decomposition and oblique projections is shown in Section \ref{sec:eq-cur}. 
}


\subsection{Time-Continuous Variational Principle}
The central idea behind TDB-based low-rank approximation is that the bases evolve optimally to minimize the residual due to low-rank approximation error. The residual is obtained by substituting an SVD-like low-rank approximation into the FOM so that $\bm V(t)$ is closely approximated by the rank-$r$ matrix
\begin{equation}\label{eqn:low-rank-approx}
    \Vhat(t) = \bU(t)\bs\Sigma(t) \bY(t)^T,
\end{equation}
where $\bU(t)\in\mathbb{R}^{n\times r}$ is a time-dependent orthonormal spatial basis for the column space, $\bY(t)\in\mathbb R^{s\times r}$ is a time-dependent orthonormal parametric basis for the row space, $\bs\Sigma(t)\in\mathbb R^{r\times r}$ is, in general, a full matrix, and $r\ll\mathrm{min}(n,s)$ is the rank of the approximation. 

Because this is a low-rank approximation, it cannot satisfy the FOM exactly and there will be a residual equal to:
\begin{equation}
    \bm R(t) = \frac{\D \left(\bU \bs\Sigma \bY ^T\right)}{\D t}-\mathcal F(t,\bU\bs\Sigma \bY^T).
\end{equation}
This residual is minimized via the first-order optimality conditions of the variational principle given by 
\begin{equation}\label{eq:var-princ}
\mathcal{J}(\dot{\bU}, \dot{\bs\Sigma}, \dot{\bY})=\left\|\frac{\D\left(\bU\bs\Sigma \bY^T\right)}{\D t}-\mathcal{F}(t,\bU\bs\Sigma \bY^T)\right\|_{F}^2,
\end{equation}
subject to orthonormality constraints on $\bU$ and $\bY$. Since the above variational principle involves the time-continuous equation (i.e. no temporal discretization is applied), the idea is to minimize the \emph{instantaneous} residual by optimally updating $\bU$, $\bs\Sigma$, and $\bY$ in time. Therefore, we refer to this as the \emph{time-continuous} variational principle. As indicated in \cite{KL07,RNB21}, the optimality conditions of Eq. \ref{eq:var-princ} lead to closed-form evolution equations for $\bU$, $\bs\Sigma$, and $\bY$:
\begin{subequations}\label{eq:DBO_evol}
\begin{align}
&\dot{\bs\Sigma}=\bU^T  \bm F \bY, \label{eq:DBO_evol_S}\\ 
&\dot{\bU}=\left(\bm I-\bU\bU^T \right) \bm F \bY\bs\Sigma^{-1}, \label{eq:DBO_evol_U} \\  
&\dot{\bY}=\left(\bm I-\bY\bY^T \right) \bm F^T  \bU\bs\Sigma^{-\T},  \label{eq:DBO_evol_Y}
\end{align}
\end{subequations}
where $\bm F \in \mathbb{R}^{n\times s}$ is a matrix defined as $\bm F =\mathcal{F}(t,\bU\bs\Sigma \bY^T)$, and $\bm I$ is the identity matrix of appropriate dimensions. The above variational principle is the same as the Dirac–Frenkel time-dependent variational principle in the quantum chemistry literature \cite{Beck:2000aa} or the  dynamical low-rank approximation (DLRA) \cite{KL07}. \cb{In \cite{KL07},  the geometry of the tangent space, $\mathcal{T}_{\Vhat}\mathcal{M}_r$, was exploited to solve the  constrained residual minimization problem given by Eq. \ref{eq:var-princ}. In this setting, the residual, $\mathcal{I}(\dot{\hat{\bm V}} ) = \|\dot{\hat{\bm V}} - \mathcal{F}(t,\hat{\bm{V}})\|_F $, is minimized with the constraint that $\dot{\hat{\bm{V}}}(t) \in \mathcal{T}_{\Vhat} \mathcal{M}_r$. The solution to the above minimization problem is obtained by 
\begin{equation}\label{eq:tangent}
    \dot{\hat{\bm V}} = \mathcal{P}_{\mathcal{T}_{\hat{\bm{V}}}} (\mathcal{F}(t,\hat{\bm{V}})),
\end{equation}
where $\mathcal{P}_{\mathcal{T}_{\hat{\bm{V}}}}$ is the orthogonal projection onto the tangent space  $\mathcal{T}_{\hat{\bm{V}}}$ at $\hat{\bm{V}}=\bm U \bs \Sigma \bm Y^T$ \cite[Lemma 4.1]{KL07}.
It is easy to show that Eqs. \ref{eq:DBO_evol_S}-\ref{eq:DBO_evol_Y} can be recovered from Eq. \ref{eq:tangent}.  
}
 As it was shown in \cite{B19}, it is possible to derive a similar variational principle for the DO decomposition, $\Vhat(t)=\bm{U}_{DO}(t)\bm{Y}^T_{DO}(t)$, whose optimality conditions are constrained to the orthonormality of the spatial modes, $\bm{U}^T_{DO}\bm{U}_{DO}= \bm{I}$, via the dynamically orthogonal condition, $\dot{\bm{U}}^T_{DO}\bm{U}_{DO}= \bm{0}$. However, for the sake of simplicity and unlike the original DO formulation presented in \cite{SL09}, an evolution equation for the mean field is not derived. Without loss of generality, the low-rank  DO evolution equations become
\begin{subequations}
\begin{align}
    &\dot{\bU}_{DO} = \left( \bm I - \bU_{DO}\bU_{DO}^T \right)\bm{F}\bY_{DO}\bm C^{-1}, \label{eqn:u-do} \\
    &\dot{\bY}_{DO} = \bm F^T\bU_{DO}, \label{eqn:y-do}
\end{align}
\end{subequations}
where $\bm C = \bY_{DO}^T\bY_{DO}$ is the low-rank correlation matrix.
Note that the low-rank approximation based on DO is \emph{equivalent} to Eq. \ref{eqn:low-rank-approx}, i.e.,  $\bm{U}_{DO}\bm{Y}^T_{DO}=\bU\bs\Sigma \bY^T$. Similarly, the BO decomposition, $\Vhat(t)=\bm{U}_{BO}(t)\bm{Y}^T_{BO}(t)$, which is subject to BO conditions, $\bm{U}^T_{BO}\bm{U}_{BO}= \mbox{diag}(\lambda_1, \dots, \lambda_r)$ and $\bm{Y}^T_{BO}\bm{Y}_{BO}= \bm{I}$, is also identical to DO and Eq. \ref{eqn:low-rank-approx}.  As it was shown in \cite{PB20}, one can derive matrix differential equations that transform the factorization $\{\bm{U},\bs{\Sigma},\bm{Y} \}$ to    $ \{\bm{U}_{DO},\bm{Y}_{DO} \}$ or $ \{\bm{U}_{BO},\bm{Y}_{BO} \}$. The equivalence of DO and BO formulations was shown in \cite{CSK14}. Using the DO/BO terminology, Eqs. \ref{eq:DBO_evol_S}-\ref{eq:DBO_evol_Y} have both DO and BO conditions, i.e., the dynamically orthogonal conditions for $\bm{U}$ and $\bm{Y}$: $\dot{\bm{U}}^T\bm{U}= \bm{0}$ and  $\dot{\bm{Y}}^T\bm{Y}= \bm{0}$ as well as bi-orthonormality conditions: $\bm{U}^T\bm{U}= \bm{I}$ and  $\bm{Y}^T\bm{Y}= \bm{I}$. Despite their equivalence, these three factorizations have  different numerical performances in the presence of small singular values. As it was shown in \cite{PB20}, Eqs. \ref{eq:DBO_evol_S}-\ref{eq:DBO_evol_Y} outperform both DO and BO. 


Despite the potential of Eqs. \ref{eq:DBO_evol_S}-\ref{eq:DBO_evol_Y} to significantly reduce the computational cost of solving massive matrix differential equations like Eq. \ref{eq:FOM}, there are still a number of outstanding challenges for most practical problems of interest. As highlighted in the Introduction, computing $\bm F = \mathcal{F}(t,\bm U \bs \Sigma \bm Y^T)$ requires $\mathcal{O}(ns)$ operations that scale with the size of the FOM. This involves applying the nonlinear map ($\mathcal{F}$) on every column of the  matrix $\hat{\bm V}=\bm U \bs \Sigma \bm Y^T$. While it is possible to achieve $\mathcal{O}(n+s)$ for the  special cases of homogeneous linear and quadratic nonlinear $\mathcal{F}$, this comes at the expense of a highly intrusive process, that requires a careful term-by-term treatment of the right side of  Eqs. \ref{eq:DBO_evol_S}-\ref{eq:DBO_evol_Y} \cite[Appendix B]{NB23}. Furthermore, solving Equations \ref{eq:DBO_evol_U} and \ref{eq:DBO_evol_Y} become unstable when $\bs\Sigma$ is singular or near singular. This is particularly problematic because  it is often necessary to retain very small singular values in order to have an accurate approximation.

While the low-rank approximation based on TDBs can be cast in different, yet equivalent formulations, we have chosen Eqs. \ref{eq:DBO_evol_S}-\ref{eq:DBO_evol_Y} over DO/BO/DDO decompositions to highlight the underlying challenges. Since, DO/BO/DDO  decompositions exhibit all of the above  challenges, addressing these challenges in the context of Eqs. \ref{eq:DBO_evol_S}-\ref{eq:DBO_evol_Y} automatically addresses the DO/BO/DDO challenges as well.

\subsection{Time-Discrete Variational Principle}
To address the  challenges of low-rank approximations based on TDB using the time-continuous variational principle, we consider a \emph{time-discrete} variational principle for rank-adaptive matrix approximations, which has recently been applied in \cite{CL23} and also \cite{kieri2019projection, rodgers2022adaptive} in the context of tensors. To this end, consider \cb{an explicit Runge-Kutta temporal discretization of Eq. \ref{eq:FOM}:
\begin{equation}\label{eq:disc_FOM}
\bm V^{k} = \hat{\bm V}^{k-1} + \Delta t \ol{\bm{F}},
\end{equation}
where $\Delta t$ is the step size and $\ol{\bm F}$ is obtained via an explicit Euler or  Runge-Kutta scheme. For example, the first-order explicit Euler method is given by: $\ol{\bm F} = \mathcal{F}(t^{k-1},\hat{\bm V}^{k-1})$. In the above equation, it is important to note that $\bm V^k$ is not the FOM solution since the right hand side is computed using the low-rank state from the previous time step. Despite using the rank-$r$ $\Vhat^{k-1}$ in Eq. \ref{eq:disc_FOM}, $\bm V^{k}$ will not be a rank-$r$ matrix, i.e., $\bm V^{k} \notin \mathcal{M}_r$. Excluding some rare exceptions, taking one step according to Eq. \ref{eq:disc_FOM} will put $\bm V^k$ off the rank-$r$ manifold. Therefore, to solve the MDE while remaining on $\mathcal {M}_r$,  a rank truncation is needed to map  the solution back onto the rank-$r$ manifold at each time step.   In other words, we need to  approximate $\bm V^k$ with a rank-$r$ matrix, $\Vhat^k$, such that
\begin{equation}\label{eq:Residual}
\bm V^k = \Vhat^k+ \bm R^k,
\end{equation}
where $\bm R^k$ is the low-rank approximation error.

The time-discrete variational principle can be stated as finding the best $\Vhat^k\in\mathcal{M}_r$ such that  the Frobenius norm of the  residual is minimized \cite{kieri2019projection}:}
\begin{equation}
\mathcal{Z}(\Vhat^k) = \bigg \| \bm V^k - \Vhat^{k}  \bigg \|_F^2.
\end{equation}
 The solution of the above residual minimization scheme is obtained via
\begin{equation}\label{eq:TD_SVD}
\Vhat_{best}^k = \mbox{\texttt{SVD}}(\bm V^k),
\end{equation}
\cb{
where $\mbox{\texttt{SVD}}(\bm V^k)=\bm{U}_{best}^k\bs \Sigma_{best}^k \bm Y_{best}^{k^T}$ is the rank-$r$ truncated SVD of matrix $\bm V^k$, where  $\bm{U}_{best}^k \in \mathbb{R}^{n\times r}$ and $\bm{Y}_{best}^k \in \mathbb{R}^{s\times r}$ are the matrices of the first $r$ left and right singular vectors of $\bm{V}^k$, respectively and $\bs \Sigma_{best}^k \in \mathbb{R}^{r \times r}$ is the matrix of singular values. }

An important advantage of Eq. \ref{eq:TD_SVD} over Eqs. \ref{eq:DBO_evol_S}-\ref{eq:DBO_evol_Y} is that the time advancement according to Eq. \ref{eq:TD_SVD}  does not become singular in the presence of small singular values. While this solves the issue of ill-conditioning, computing Eq. \ref{eq:TD_SVD} at each iteration of the time stepping scheme is cost prohibitive. This computational cost is due to two sources: (i) computing the nonlinear map \cb{$\mathcal{F}(t^{k-1},\Vhat^{k-1})$ to obtain $\bm{V}^k$} and (ii) computing $\mbox{\texttt{SVD}}(\bm V^k)$. The cost of (i) alone makes the solution of the time-discrete variational principle as expensive as the FOM, i.e., $\mathcal{O}(ns)$. Besides the flops cost associated with computing $\bm V^k$, the memory cost of storing  $\bm V^k$ is prohibitive for most realistic applications. On the other hand, computing the exact SVD of $\bm V^k$ scales with $\mathrm{min} \{ \mathcal{O}(n^3),\mathcal{O}(s^3)\}$. While this cost is potentially alleviated by fast algorithms for approximating the SVD, e.g. randomized SVD \cite{halko2011finding} or incremental QR \cite{sorensen2016deim}, for general nonlinearities in $\mathcal{F}$, (i) is unavoidable. This ultimately leads to a computational cost that exceeds that of the FOM.


\subsection{Low-Rank Approximation via an Oblique Projection}\label{sec:sparse-sampling}
\begin{figure}
    \centering
    \includegraphics[width=\textwidth]{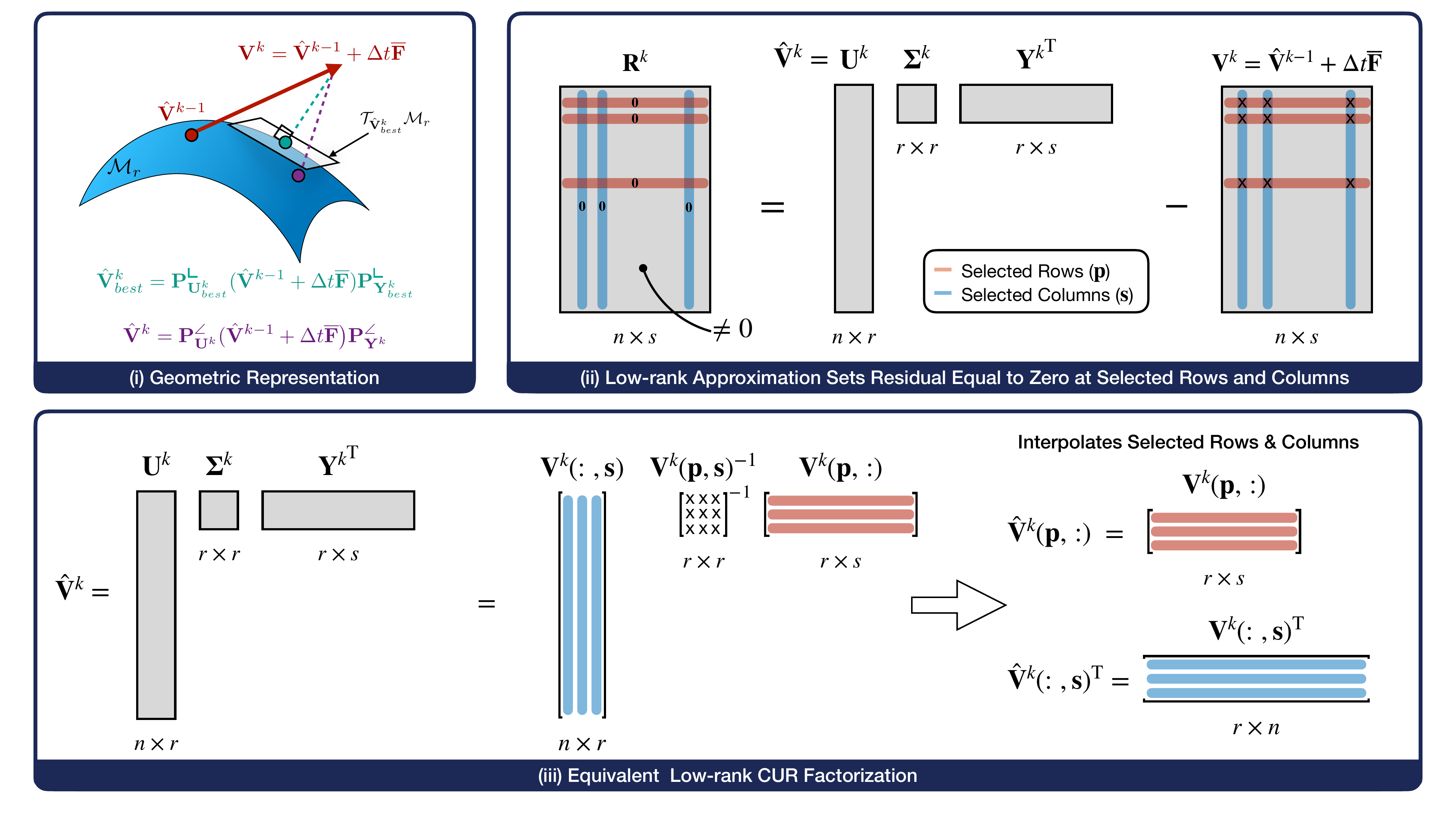}
    \caption{Schematic of the TDB-CUR methodology. \cb{(i) A geometric representation of the methodology depicting the departure from the rank-$r$ manifold at $\bm \Vhat^{k-1}$ to $\bm V^k$ (red). Two possible mappings that truncate $\bm V^k$ back to the rank-$r$ manifold are shown: $\Vhat^k_{best}$ (green) computed via orthogonal projection (SVD) and $\Vhat^k$ (purple) computed via oblique projection (TDB-CUR). The error of $\Vhat^k_{best}$ (dashed green line) is orthogonal to the tangent space $\mathcal{T}_{\Vhat^k_{best}}\mathcal{M}_r$ at $\Vhat^k_{best}$, while the error of $\Vhat^k$ (dashed purple line) is not}.   (ii) \cb{The low-rank approximation,} $\Vhat^k$, can be computed such that the residual at the selected rows (red) and columns (blue) is equal to zero. This is accomplished via sparse interpolation of the selected rows and columns. (iii) \cb{Using interpolatory projection}, the resulting $\Vhat^k$ is an approximation to the rank-$r$ truncated $\texttt{\mbox{SVD}}(\bm V^k)$, and is equivalent to the low-rank CUR factorization that interpolates the selected rows and columns of $\bm V^k$. Although it is equivalent to this CUR factorization, the numerical computation of $\Vhat^k$ is different, as it does not require inverting $\bm V^k(\bm p, \bm s)$.}
    \label{fig:schematic}
\end{figure}

To overcome these challenges, \cb{ we present an oblique projection scheme   that enables a cost-effective approximation to the rank-$r$ $\mbox{\texttt{SVD}}(\bm V^k)$. Before presenting our methodology, we provide a geometric interpretation of   $\mbox{\texttt{SVD}}(\bm V^k)$. In particular, $\mbox{\texttt{SVD}}(\bm V^k)$ can be interpreted as an orthogonal projection onto the manifold $\mathcal{M}_r$ at  $\Vhat_{best}^k$, since $\Vhat_{best}^k$ can be expressed as the orthogonal projection  of $\bm V^k$ onto the tangent space at  $\Vhat_{best}^k$, i.e., $\Vhat_{best}^k = \mathcal{P}_{\mathcal{T}_{\hat{\bm{V}}^k_{best}}} (\bm V^k)$. Therefore, $\Vhat_{best}^k$  can be expressed as: 
\begin{equation}\label{eq:svd_orth}
     \Vhat_{best}^k= \bm{P}^{\Langle}_{\bm{U}_{best}^k} \bm V^k  \bm{P}^{\Langle}_{\bm{Y}_{best}^k} = \bm{P}^{\Langle}_{\bm{U}_{best}^k} ( \hat{\bm V}^{k-1} + \Delta t \ol {\bm F})  \bm{P}^{\Langle}_{\bm{Y}_{best}^k}, 
\end{equation} 
where  $\bm{P}^{\Langle}_{\bm{U}_{best}^k} = \bm{U}_{best}^k\bm{U}_{best}^{k^T}$  and $\bm{P}^{\Langle}_{\bm{Y}_{best}^k} = \bm{Y}_{best}^k\bm{Y}_{best}^{k^T}$ are orthogonal projections onto the column and row space of $\Vhat_{best}^k$, respectively. The approximation $\Vhat_{best}^k$ is the optimal rank-$r$ approximation of $\bm V^k$, however as mentioned in the previous section, computing $\Vhat_{best}^k$ is more expensive than solving the FOM.

In the following, we present a methodology that computes an accurate approximation to  $\Vhat_{best}^k$ in a cost-effective manner. From the geometric perspective, our approach is to use an oblique projection onto a set of rank-$r$ orthonormal column ($\bm{U}^k$)  and row ($\bm{Y}^k$) subspaces:
\begin{equation}\label{eq:oblique}
     \Vhat^k= \bm{P}^{\Oangle}_{\bm{U}^k} \bm V^k  \bm{P}^{\Oangle}_{\bm{Y}^k} = \bm{P}^{\Oangle}_{\bm{U}^k} ( \hat{\bm V}^{k-1} + \Delta t \ol {\bm F})  \bm{P}^{\Oangle}_{\bm{Y}^k}.
\end{equation} 
The above equation is analogous to Eq. \ref{eq:svd_orth} where the orthogonal projections are replaced with oblique projections. While orthogonal projection requires access to the entire  $\hat{\bm V}^{k-1} + \Delta t \ol {\bm F}$ matrix, the oblique projectors can be designed to require the computation of $\hat{\bm V}^{k-1} + \Delta t \ol {\bm F}$ at $\mathcal{O}(r)$ columns and rows. This geometric perspective is depicted in Figure \ref{fig:schematic} panel (i).

From the matrix decomposition point of view,  the above approximation may be represented via a CUR decomposition:
\begin{equation}\label{eq:TD_CUR}
\Vhat^k = \mbox{\texttt{CUR}}(\bm V^k),
\end{equation}
where \texttt{CUR} represents the algorithmic implementation of a CUR decomposition. See Figure \ref{fig:schematic} panel (iii).}

\cb{From the residual minimization perspective, the SVD can be viewed as a Galerkin projection where $\|\bm R^k \|_F$ is minimized. On the other hand, the presented approach based on interpolatory projection (a special case of oblique projection) can be viewed as a \emph{collocated} scheme where the residual is set}  to zero at $r$ strategically selected rows and columns of the residual matrix, $\bm R^k$, \cb{in Eq. \ref{eq:Residual}}. To this end, we present an algorithm to set $\bm R^k(\bm p,:) = \bm 0$ and $\bm R^k(:,\bm s)=\bm 0$, where $\bm p\in\mathbb N^r$ and $\bm s \in\mathbb N^r$ are vectors containing the row and column indices at which the residual is set to zero. This simply requires $\Vhat^k(\bm p,:) = \bm  V^k(\bm p,:)$ and $\Vhat^k(:,\bm s) =  \bm V^k(:,\bm s)$. \cb{See Figure \ref{fig:schematic} panel (ii). In Section \ref{sec:OS}, we consider oblique projection for the general case when $r'>r$ (Definition \ref{def:obl_proj}), where the residual at the selected rows and columns is not guaranteed to be zero.}  

\cb{Although the approach we will present is equivalent to the oblique projection of Eq. \ref{eq:oblique}, $\bm U^k$ and $\bm Y^k$ are the unknown column and row bases of $\Vhat^k$ at the \emph{current} time step. Therefore, Eq. \ref{eq:oblique} cannot be readily used for the computation of $\Vhat^k$. In the following, we present a methodology to compute these bases (and subsequently $\Vhat^k$) by strategically sampling $r$ columns and rows of  $\bm V^k$.}
While there are many possible choices for the indices $\bm p$ and $\bm s$, selecting these points should be done in a principled manner, to ensure the residual at all points remains small. To compute these points, we use the discrete empirical interpolation method (DEIM) \cite{chaturantabut2010nonlinear} which has been shown to provide near optimal sampling points for computing CUR matrix decompositions \cite{sorensen2016deim}. A similar approach was recently applied in \cite{NB23} to accelerate the computation of Eqs. \ref{eq:DBO_evol_S}-\ref{eq:DBO_evol_Y}, by only sampling $\bm F$ at a small number of rows and columns. However, the approach presented in \cite{NB23} still suffers from the issue of ill-conditioning.


To compute the DEIM points, the rank-$r$ SVD (or an approximation) is required \cite{sorensen2016deim}. Since we do not have access to the rank-$r$ SVD at the current time step, $k$, we use the approximation of the SVD from the previous time step, $\Vhat^{k-1}=\bU^{k-1}\bs\Sigma^{k-1}\bY{^{k-1}}^T$, to compute the DEIM points. \cb{The initial approximation is ideally obtained from the FOM initial condition as the rank-$r$ $\mbox{\texttt{SVD}}(\bm V_0)$.} The algorithm for computing $\Vhat^{k}$ \cb{using interpolation} is as follows:
\begin{enumerate}[leftmargin=*, label=(\roman*)]
    \item Compute the sampling indices, $\bm p \gets \texttt{DEIM}(\bU^{k-1})$, and $\bm s \gets \texttt{DEIM}(\bY^{k-1})$, in parallel.
    \item \cb{Compute $\bm V^k(\bm p,:)$ and $\bm V^k(:,\bm s)$ by taking one step according to Eq. \ref{eq:disc_FOM}} at the selected rows and columns, in parallel.
    \item Compute $\bm Q\in\mathbb R^{n\times r}$ as the orthonormal basis for the range of $\bm V^k(:,\bm s)$ by QR decomposition such that $\bm V^k(:,\bm s) = \bm Q \bm R$, where $\bm R \in \mathbb{R}^{r\times r}$.
    \item Interpolate every column of $\bm{V}^k$  onto the orthonormal basis $\bm Q$ at sparse indices $\bm p$: 
    \begin{equation} \label{eqn:oblique-projection}
        \bm Z = \bm Q(\bm p,:)^{-1}\bm V^k(\bm p,:),
    \end{equation}
    where $\bm Z\in\mathbb R^{r\times s}$ is the matrix of interpolation coefficients  such that $\bm Q \bm Z$ \cb{interpolates} $\bm V^k$ onto the basis $\bm Q$ at the interpolation points indexed by $\bm p$. 
    \item Compute the SVD of $\bm Z$ so that
    \begin{equation}
        \bm Z = \bU_{\bm Z}\bs\Sigma^k{\bY^k}^T,
    \end{equation}
    where $\bm U_{\bm Z} \in \mathbb{R}^{r \times r}$,  $\bs\Sigma^k \in \mathbb{R}^{r \times r}$, and $\bm Y^k \in \mathbb{R}^{s \times r}$.
    \item Compute $\bU^k \in\mathbb R^{n\times r}$ as the in-subspace rotation:
    \begin{equation}
        \bU^k = \bm Q\bU_{\bm Z}.
    \end{equation}
\end{enumerate}
In Step (i), the details of the \texttt{DEIM} algorithm can be found in \cite[Algorithm 1]{chaturantabut2010nonlinear}. A DEIM algorithm based on the QR factorization, a.k.a \texttt{QDEIM}, may also be used \cite{DG16,MBKB18}. Both \texttt{DEIM} and \texttt{QDEIM} are sparse selection algorithms  and they perform comparably in the cases considered in this paper. We explain here how the above algorithm addresses the three challenges mentioned in  the Introduction (\ref{sec:Intro}). 
\begin{enumerate}[label=(\roman*)]
    \item \textbf{Computational efficiency:} The above procedure returns the updated low-rank approximation $\Vhat^k = \bm{QZ} = \bU^k\bs\Sigma^k{\bY^k}^T$, and only requires sampling $\bm V^k$ at $r$ rows and columns. This alone significantly reduces both the required number of flops and memory, compared to computing the entire $\bm V^k$. Furthermore, instead of directly computing the SVD of the $n\times s$ matrix $\bm V^k$, we only require computing the QR of the $n\times r$ matrix $\bm V^k(:,\bm s)$, and the SVD of the $r\times s$ matrix $\bm Z$. This reduces the computational cost to $\mathcal{O}(s+n)$ for $r \ll s$ and $r \ll n$. Moreover, in most practical applications, computing  $\bm V^k(:,\bm s)$ is the costliest part of the algorithm, which requires solving $s$ samples of the FOM. However, since these samples are independent of each other, the columns of $\bm V^k(:,\bm s)$ can  be computed in parallel. Similarly, each row of $\bm V^k(\bm p,:)$ can be computed in parallel.
    \item \textbf{Intrusiveness:} While this significantly reduces the computational burden, perhaps an equally important outcome is the minimally intrusive nature of the above approach. For example, when the columns of $\bm V^k$ are independent, e.g. random samples, $\bm V^k(:,\bm s)$ can be computed by directly applying Eq. \ref{eq:disc_FOM} \cb{to the selected columns of the low-rank approximation from the previous time step.} This effectively allows for existing numerical implementations  of Eq. \ref{eq:disc_FOM} to be used as a black box for computing $\bm V^k(:,\bm s)$. \cb{The nonintrusive column sampling in the presented algorithm is the counterpart of  solving Eq. (\ref{eq:DBO_evol_U}) in DLRA and  Eq. (\ref{eqn:u-do}) in DO.  However, Eq. (\ref{eq:DBO_evol_U}) and Eq. (\ref{eqn:u-do})  require deriving and implementing new PDEs, whereas the presented algorithm allows an existing deterministic solver to be used in a black box fashion, in which a suitable column space basis is extracted.}  On the other hand, the rows of $\bm V^k$ are in general \emph{dependent}, based on a known map for the chosen spatial discretization scheme, e.g. sparse discretizations like finite difference, \cb{ finite element, or dense discretization schemes, e.g. global spectral methods}. Therefore, computing $\bm V^k(\bm p,:)$ does require specific knowledge of the governing equations, namely the discretized differential operators. Based on the discretization scheme, one can determine a set of adjacent points, $\bm p_a$, that are required for computing the derivatives at the points specified by $\bm p$. While this introduces an added layer of complexity, this is much less intrusive than deriving and implementing reduced order operators for each term in the governing equations; which we emphasize again, is only feasible for homogeneous linear or quadratic nonlinear equations. In the present work, that bottleneck is removed, regardless of the type of nonlinearity.
    \item \textbf{Ill-conditioning:} The presented algorithm  is robust in the presence of small or zero singular values. First note that the inversion of the matrix of singular values is not required in the presented algorithm. In fact,  the conditioning of the algorithm depends on $\bm Q(\bm p,:)$ and $\bm Y(\bm s,:)$, and the DEIM algorithm ensures that these two matrices are well-conditioned. To illustrate this point, let us consider the case of overapproximation where the rank of $\bm V^k(:,\bm s)$ is $r_1<r$. In this case, Eqs. \ref{eq:DBO_evol_S}-\ref{eq:DBO_evol_Y} and Eqs. \ref{eqn:u-do}-\ref{eqn:y-do} cannot be advanced because $\bs \Sigma \in \mathbb{R}^{r\times r}$ and $\bm C \in \mathbb{R}^{r\times r}$ will be singular, i.e., rank$(\bs \Sigma)=$rank$(\bm C) = r_1<r$. On the other hand, despite $\bm V^k(:,\bm s)$ being rank deficient, $\bm Q$ will still be a full rank matrix in the presented algorithm.   While there is no guarantee that a subset of rows of $\bm Q$, i.e., $\bm Q (\bm p,:)$ is well conditioned, the DEIM  is a greedy algorithm that is designed to keep  $\| \bm{Q}(\bm p, :)^{-1} \|$ as small as possible in a near-optimal fashion. In Section \ref{sec:OS}, we show that oversampling further improves the condition number of the presented algorithm, and in Theorem \ref{thm:err}, we show that $\| \bm{Y}(\bm s, :)^{-1} \|$ plays an equally important role in maintaining a well-conditioned algorithm.
\end{enumerate}

As we will show in Section \ref{sec:eq-cur}, the low-rank approximation computed above is equivalent to a CUR matrix decomposition that interpolates $\bm V^k$ at the selected rows and columns. Therefore, we refer to the above procedure as the TDB-CUR algorithm. 

\subsection{Computing $\bm V^k(\bm p,:)$}\label{sec:compute-G}
Up until this point, we have considered $\bm V^k = \hat{\bm V}^{k-1} + \Delta t \ol{\bm{F}}$ to be an $n\times s$ matrix resulting from an \cb{explicit Runge-Kutta temporal discretization} of Eq. \ref{eq:disc_FOM}. We showed that sparse row and column measurements, $\bm V^k(\bm p,:)$ and $\bm V^k(:,\bm s)$, could be used to efficiently compute an approximation to the rank-$r$ SVD of $\bm V^k$. While $\bm V^k(:,\bm s)$ is straightforward to compute for independent random samples, as discussed in Section \ref{sec:sparse-sampling}, computing $\bm V^k(\bm p,:)$ depends on a set of adjacent points, $\bm p_a$, according to the spatial discretization scheme.  As a result, for higher-order integration schemes, special care must be taken in the computation of $\bm V^k(\bm p,:)$. To demonstrate this, we consider the second-order explicit Runge-Kutta scheme where
\cb{
\begin{equation*}
    \ol{\bm F} = \mathcal{F}\left(t^{k-1}+\frac{1}{2}\Delta t, \, \Vhat^{k-1}+\frac{1}{2}\Delta t\mathcal{F}\left(t^{k-1},\,\Vhat^{k-1}\right)\right)
\end{equation*}
}
After determining the row indices, $\bm p$ and $\bm p_a$, $\bm V^k(\bm p, :)$ can be computed as follows:
\begin{enumerate}[leftmargin=*, label=(\roman*)]
    \item Compute $\Vhat^{k-1}([\bm p,\bm p_a],:)=\bU^{k-1}([\bm p,\bm p_a],:)\bs\Sigma^{k-1}{\bY^{k-1}}^T$.
    \item \label{compute-F1} Compute the first stage $\bm F_1 = \mathcal{F}\left(t^{k-1},\Vhat^{k-1}\right)$ at the $\bm p$ rows as
    \begin{equation*}
        \bm F_1(\bm p,:)=\mathcal F\left(t^{k-1},\Vhat^{k-1}([\bm p,\bm p_a],:)\right). 
    \end{equation*}
    Note, if the explicit Euler method is used, \cb{$\ol{\bm F}=\bm F_1$}, and no additional steps are required. Simply compute    $\bm V^k(\bm p,:) = \Vhat^{k-1} + \Delta t\bm F_1(\bm p, :)$. If a higher-order scheme is used, proceed with the following steps.
    \item The final stage of the second order integration scheme requires taking a half step to evaluate $\mathcal F$ at the midpoint:
    \begin{equation*}
        \bm F_2 = \mathcal F \left(t^{k-1}+\frac{1}{2}\Delta t,  \Vhat^{k-1}+\frac{1}{2}\Delta t\bm F_1\right).
    \end{equation*}
    \cb{Note that for the second order scheme, $\ol{\bm F} = \bm F_2$.} Here, we require $\bm F_2(\bm p,:)$, given by
    \begin{equation*}
        \bm F_2(\bm p, :) = \mathcal F \left(t^{k-1}+\frac{1}{2}\Delta t,  \Vhat^{k-1}([\bm p,\bm p_a],:)+\frac{1}{2}\Delta t\bm F_1([\bm p,\bm p_a],:)\right).
    \end{equation*}
    Notice that we now require $\bm F_1(\bm p_a,:)$ to evaluate the above expression. While this can be computed according to Step \ref{compute-F1}, where $\bm p_a$ will have its own set of adjacent points $\bm p_{aa}$, this process quickly gets out of hand, especially as more stages are added to the integration scheme. As a result, for higher-order schemes, the efficiency afforded by the presented algorithm will deteriorate, and the resulting implementation will become increasingly complex. To overcome these challenges, we instead compute the low-rank approximation $\hat{\bm F}_i\approx\bm F_i$, using the sparse row and column measurements, $\bm F_i(\bm p,:)$ and $\bm F_i(:,\bm s)$, which are already required for computing $\bm V^k(\bm p,:)$ and $\bm V^k(:,\bm s)$. Here, the subscript denotes the $i^{\mathrm{th}}$ stage of the integration scheme. The first step is to compute $\bm U_{\bm F_i}$ as an orthonormal basis for the $\mathrm{Ran}(\bm F_i(:,\bm s))$, using QR. Next, compute the oblique projection of $\bm F_i$ onto $\bm U_{\bm F_i}$, such that 
    \begin{equation*}
        \hat{\bm F}_i=\bm U_{\bm F_i}\bm U_{\bm F_i}(\bm p,:)^{-1}\bm F_i(\bm p,:). 
    \end{equation*}
    Using this low-rank approximation, $\bm F_i(\bm p_a,:)$ is readily approximated by  $\hat{\bm F}_i(\bm p_a,:)=\bm U_{\bm F_i}(\bm p_a,:)\bm U_{\bm F_i}(\bm p,:)^{-1}\bm F_i(\bm p,:)$.
    Although we have considered the second-order Runge-Kutta method in the example above, this approach is easily extended to higher-order Runge-Kutta methods. \cb{It is straightforward to show that the above procedure is equivalent to a CUR decomposition of  matrix $\bm F_i$, similar to our previous work \cite{NB23}. }
        
\end{enumerate}

\subsection{Equivalence to a CUR Decomposition \& Oblique Projection}\label{sec:eq-cur}
\cb{Before presenting the details of our methodology in Section \ref{sec:sparse-sampling}, we discussed how the presented approach can be understood as an oblique projection (Eq. \ref{eq:oblique}) or alternatively as a CUR decomposition (Eq. \ref{eq:TD_CUR}).} In this section, we  show that (i) the presented algorithm is equivalent to a CUR decomposition (Theorem \ref{thm:CUR}), and (ii) the matrix $\Vhat^k$ is obtained via an oblique projection, which  requires access to only the selected rows and columns of $\bm V^k$ (Theorem \ref{thm:obliq}).  \cb{In Theorem \ref{thm:CUR} and Theorem \ref{thm:obliq}, $\bm P$ and $\bm S$ are matrices of size $n\times r$ and $s\times r$, respectively. We drop the time step  $k$ for simplicity. } 
 
\begin{theorem}\label{thm:CUR}
Let $\Vhat=\bm{QZ}$ be the low-rank approximation of $\bm V$ computed according to the TDB-CUR algorithm. Then: (i) $\Vhat=\bm{QZ}$ is equivalent to the CUR factorization given by $(\bm{VS}) (\bm P^T\bm{VS})^{-1} (\bm P^T\bm V)$. (ii) The low-rank approximation is exact at the selected rows and columns, i.e. $\bm P^T\Vhat = \bm P^T\bm V$ and $\Vhat\bm S=\bm{VS}$. 
 \begin{proof}\leavevmode
 \begin{enumerate}[leftmargin=*, align=left, wide=0pt, label=(\roman*)]
 \item According to the TDB-CUR algorithm, $\bm Q$ is a basis for the $\mathrm{Ran}(\bm V(:,\bm s))$. Therefore, $\bm V(:,\bm s)=\bm{VS} = \bm{QQ}^T\bm V\bm S$, and it follows that $\bm P^T\bm{VS}=\bm P^T\bm{QQ}^T\bm{VS}$. Substituting this result into the CUR factorization gives
 \begin{equation*}
    (\bm{VS}) (\bm P^T\bm{VS})^{-1} (\bm P^T\bm V) = \bm{QQ}^T\bm{VS}( \bm P^T\bm Q\bm Q^T\bm{VS} )^{-1}\bm P^T \bm V.
 \end{equation*}
 Rearranging the above expression gives the desired result 
 \begin{equation*}
     (\bm{VS}) (\bm P^T\bm{VS})^{-1} (\bm P^T\bm V) = \bm{QQ}^T\bm{VS}(\bm Q^T\bm{VS})^{-1}(\bm P^T \bm Q)^{-1}\bm P^T \bm V=\bm{QZ}=\Vhat,
 \end{equation*}
 where we have used $\bm Z=(\bm P^T\bm Q)^{-1}\bm P^T \bm V = \bm Q(\bm p,:)^{-1}\bm V(\bm p,:)$, from Eq. \ref{eqn:oblique-projection}.

 \item Using the above result, $\Vhat = (\bm{VS}) (\bm P^T\bm{VS})^{-1} (\bm P^T\bm V)$, we show the selected rows of $\Vhat $ are exact, i.e., $\Vhat(\bm p, :) = \bm V (\bm p,:)$:
 \begin{equation*}
     \bm P^T\Vhat = (\bm P^T\bm{VS}) (\bm P^T\bm{VS})^{-1} (\bm P^T\bm V) = \bm P^T\bm V. 
 \end{equation*}
 Similarly for the columns,
 \begin{equation*}
     \Vhat\bm S = (\bm{VS}) (\bm P^T\bm{VS})^{-1} (\bm P^T\bm{VS}) = \bm{VS}. 
 \end{equation*}
 This completes the proof.
\end{enumerate}
\end{proof}
\end{theorem}
Now we  show that $\Vhat$ is   an oblique projection of $\bm V$ onto the selected columns and rows of $\bm V$. In particular, the oblique projector involved is an  interpolatory projector. For the sake of brevity, we drop the superscript $k$ in the following.  

 \begin{theorem}\label{thm:obliq}
Let $\Vhat=\bm{QZ}$ be the low-rank approximation of $\bm V$ computed according to the TDB-CUR algorithm. Then $\Vhat=\bm{P}^{\Oangle}_{\bm{U}}\bm V\bm{P}^{\Oangle}_{\bm{Y}}$ where $\bm{P}^{\Oangle}_{\bm{U}}$ and $\bm{P}^{\Oangle}_{\bm{Y}}$ are oblique projectors onto Ran($\bm U$) and Ran($\bm Y$), respectively, according to Eq. \ref{eq:interp_proj}.  
\begin{proof}\leavevmode
We first show that  $\bm{P}^{\Oangle}_{\bm{U}}$ can be represented versus   $\bm Q$  as the interpolation basis. To this end, replacing $\bm U = \bm Q \bm U_{\bm Z}$ in the definition of $\bm{P}^{\Oangle}_{\bm{U}}$ results in:
\begin{equation*}
    \bm{P}^{\Oangle}_{\bm{U}} = \bm U (\bm P^T \bm U)^{-1} \bm P^T = \bm Q \bm U_{\bm Z} (\bm P^T \bm  Q \bm U_{\bm Z})^{-1} \bm P^T =  \bm Q \bm U_{\bm Z} \bm U_{\bm Z}^{-1}(\bm P^T \bm  Q )^{-1} \bm P^T =\bm Q \bm (\bm P^T \bm  Q )^{-1} \bm P^T. 
\end{equation*}
where we have used the fact that $\bm U_{\bm Z}$ is a square orthonormal matrix and therefore, $\bm U_{\bm Z} \bm U_{\bm Z}^{-1}=\bm{I}$. Similarly, $\bm{P}^{\Oangle}_{\bm{Y}}$ can be represented versus   $\bm Z^T$  as the interpolation basis by replacing $\bm Y^T = \bs \Sigma^{-1} \bm U_{\bm Z}^{-1}\bm Z$ in $\bm{P}^{\Oangle}_{\bm{Y}}$:
\begin{equation*}
   \bm{P}^{\Oangle}_{\bm{Y}} = \bm S (\bm Y^T \bm S)^{-1} \bm Y^T = \bm S ( \bs \Sigma^{-1} \bm U_{\bm Z}^{-1} \bm Z \bm S)^{-1}  \bs \Sigma^{-1} \bm U_{\bm Z}^{-1} \bm Z = \bm S ( \bm Z \bm S)^{-1}   \bm Z. 
\end{equation*}
Using these projection operators we have
\begin{equation}\label{eqn:aux1}
    \bm{P}^{\Oangle}_{\bm{U}}\bm V\bm{P}^{\Oangle}_{\bm{Y}} = \bm Q \bm (\bm P^T \bm  Q )^{-1} \bm P^T \bm V \bm S ( \bm Z \bm S)^{-1}   \bm Z. 
\end{equation}
Using the results of Theorem \ref{thm:CUR}, Part (ii), we have: $\bm V(\bm p, \bm s) = \hat{\bm V}(\bm p, \bm s)$. Therefore:
\begin{equation*}
    \bm P^T \bm V \bm S = \bm V(\bm p, \bm s) = \hat{\bm V}(\bm p, \bm s) = \bm Q(\bm p,:) \bm Z(:,\bm s)= \bm P^T \bm Q \bm Z \bm S.
\end{equation*}
Using this result in Eq. \ref{eqn:aux1}, yields:
\begin{equation*}
    \bm{P}^{\Oangle}_{\bm{U}}\bm V\bm{P}^{\Oangle}_{\bm{Y}} = \bm Q \bm (\bm P^T \bm  Q )^{-1} \bm P^T \bm Q \bm Z \bm S( \bm Z \bm S)^{-1}   \bm Z = \bm{Q} \bm{Z} = \hat{\bm{V}}. 
\end{equation*}
This result completes the proof.  
\end{proof}
\end{theorem}

In the following theorem, we show that the oblique projection error is bounded by an error factor multiplied by the maximum of orthogonal projection errors onto $\bm U$ or $\bm Y$. We follow a similar procedure that was used in \cite{sorensen2016deim}, however,  in \cite{sorensen2016deim} the CUR is computed based on orthogonal projections onto the selected columns and rows, whereas in the presented TDB-CUR algorithm, \cb{oblique} projectors are used. \cb{Without loss of generality, we consider a generic oblique projection, where the indexing matrices $\bm P$ and $\bm S$ are of size $n\times r'$ and $s\times r'$, respectively, and in general, $r'\ge r$ (see Definition \ref{def:obl_proj} for details).} In the following, we use the second norm ($\| \sim \| \equiv \| \sim \|_2$).

\begin{theorem}\label{thm:err}
Let $\bm{P}^{\Oangle}_{\bm{U}}$ and $\bm{P}^{\Oangle}_{\bm{Y}}$ be oblique projectors according to Definition \ref{def:obl_proj}  and let $\bm U \in \mathbb{R}^{n\times r}$ and $\bm Y \in \mathbb{R}^{s\times r}$ be a set of orthonormal matrices, i.e., $\bm U^T \bm U = \bm I$ and $\bm Y^T \bm Y = \bm I$. \cb{Let   $\epsilon_f \geq 0$ be   given by: $\epsilon_f = \mbox{min}\{ \eta_p(1+\eta_s), \eta_s(1+\eta_p)\}-1$}, where $\eta_p = \| (\bm{P}^T \bm{U})^{\dagger}\|$ and $\eta_s = \| (\bm{S}^T \bm{Y})^{\dagger}\|$ and $\hat{\sigma}_{r+1} = \mbox{max}\{\|\bm V -\bm U \bm U^T \bm V\|,\| \bm V - \bm V\bm Y \bm Y^T \| \}$. Then the error of the oblique projection is bounded by
\begin{equation}\label{eq:err_bnd}
    \|\bm V - \bm{P}^{\Oangle}_{\bm{U}} \bm V \bm{P}^{\Oangle}_{\bm{Y}}\| \leq \cb{(1+\epsilon_f) \hat{\sigma}_{r+1}}.
\end{equation}  
\begin{proof}\leavevmode
\cb{First note that $\epsilon_f \geq 0$ because  $\eta_p\geq 1$ and $\eta_q\geq 1$.} The error matrix can be written as:
\begin{equation*}
    \bm V  - \bm{P}^{\Oangle}_{\bm{U}} \bm V \bm{P}^{\Oangle}_{\bm{Y}} = (\bm I - \bm{P}^{\Oangle}_{\bm{U}})\bm V +  \bm{P}^{\Oangle}_{\bm{U}}\bm V - \bm{P}^{\Oangle}_{\bm{U}} \bm V \bm{P}^{\Oangle}_{\bm{Y}} = (\bm I - \bm{P}^{\Oangle}_{\bm{U}})\bm V +  \bm{P}^{\Oangle}_{\bm{U}}\bm V (\bm I -  \bm{P}^{\Oangle}_{\bm{Y}})
\end{equation*}
where $\bm I$ is the identity matrix of appropriate size. Also, $\mathcal P \bm U = \bm U(\bm P^T \bm U)^{-1} \bm P^T \bm U = \bm U$. Therefore, $(\bm I -\bm{P}^{\Oangle}_{\bm{U}})\bm U = \bm 0$. Similarly,  $\bm Y^T (\bm I -\bm{P}^{\Oangle}_{\bm{Y}})= \bm 0$. Therefore,
\begin{align*}
    \| \bm V  - \bm{P}^{\Oangle}_{\bm{U}} \bm V \bm{P}^{\Oangle}_{\bm{Y}}\| &\leq  \|(\bm I - \bm{P}^{\Oangle}_{\bm{U}})\bm V \|+  \|\bm{P}^{\Oangle}_{\bm{U}}\bm V (\bm I -  \bm{P}^{\Oangle}_{\bm{Y}})\| \\
    & = \|(\bm I - \bm{P}^{\Oangle}_{\bm{U}})(\bm V - \bm U \bm U^T \bm V )\|+  \|\bm{P}^{\Oangle}_{\bm{U}}(\bm V - \bm V \bm Y \bm Y^T) (\bm I -  \bm{P}^{\Oangle}_{\bm{Y}})\|\\
    & \leq  \big( \|(\bm I - \bm{P}^{\Oangle}_{\bm{U}})\|+  \|\bm{P}^{\Oangle}_{\bm{U}} \|  \|(\bm I -  \bm{P}^{\Oangle}_{\bm{Y}}) \| \big) \hat{\sigma}_{r+1}\\
    & = \eta_p(1 + \eta_s)\hat{\sigma}_{r+1}.
\end{align*}
In the above inequality, we have made use of the fact that $\|\bm I- \bm{P}^{\Oangle}_{\bm{U}}\| = \| \bm{P}^{\Oangle}_{\bm{U}}\| = \eta_p$ and $\|\bm I- \bm{P}^{\Oangle}_{\bm{Y}}\| = \| \bm{P}^{\Oangle}_{\bm{Y}}\| = \eta_s$ as long as the projectors are neither null nor the identity \cite{S06}. In the second line of the above inequality, we have made use of $(\bm I -\bm{P}^{\Oangle}_{\bm{U}})\bm U = \bm 0$ and $\bm Y^T (\bm I -\bm{P}^{\Oangle}_{\bm{Y}})= \bm 0$.  Similarly, it is possible to express the error matrix as:
\begin{equation*}
    \bm V  - \bm{P}^{\Oangle}_{\bm{U}} \bm V \bm{P}^{\Oangle}_{\bm{Y}} = \bm V(\bm I - \bm{P}^{\Oangle}_{\bm{Y}}) + \bm V \bm{P}^{\Oangle}_{\bm{Y}}  - \bm{P}^{\Oangle}_{\bm{U}} \bm V \bm{P}^{\Oangle}_{\bm{Y}} = \bm V(\bm I - \bm{P}^{\Oangle}_{\bm{Y}}) +  (\bm I - \bm{P}^{\Oangle}_{\bm{U}})\bm V \bm{P}^{\Oangle}_{\bm{Y}}.
\end{equation*}
Therefore, another error bound can be obtained as
\begin{align*}
    \| \bm V  - \bm{P}^{\Oangle}_{\bm{U}} \bm V \bm{P}^{\Oangle}_{\bm{Y}}\| 
    &\leq  \| \bm V(\bm I - \bm{P}^{\Oangle}_{\bm{Y}})\| +  \|(\bm I - \bm{P}^{\Oangle}_{\bm{U}})\bm V \bm{P}^{\Oangle}_{\bm{Y}}\|\\
    & =  \| (\bm V - \bm V \bm Y \bm Y^T)(\bm I - \bm{P}^{\Oangle}_{\bm{Y}})\| +  \|(\bm I - \bm{P}^{\Oangle}_{\bm{U}})(\bm V - \bm U \bm U^T \bm V) \bm{P}^{\Oangle}_{\bm{Y}}\|\\
    & \leq  \big( \|(\bm I - \bm{P}^{\Oangle}_{\bm{Y}})\|+  \|\bm I - \bm{P}^{\Oangle}_{\bm{U}} \|  \| \bm{P}^{\Oangle}_{\bm{Y}} \| \big) \hat{\sigma}_{r+1}\\
    & = \eta_s(1 + \eta_p)\hat{\sigma}_{r+1}.
\end{align*}
where $\|\bm I- \bm{P}^{\Oangle}_{\bm{Y}}\|=\eta_s$ is used. Combining the above two inequalities yields inequality  \ref{eq:err_bnd}. 
\end{proof}
\end{theorem}
In the above error bound, when $\bm U$ and  $\bm Y$ are the $r$ most dominant  exact left and right singular vectors of $\bm V$, then $\hat{\sigma}_{r+1}=\sigma_{r+1}$, where $\sigma_{r+1}$ is the $r+1$-th singular value of $\bm V$, since 
\begin{equation}\label{eq:aux_svd}
    \|\bm V -\bm U \bm U^T \bm V\|=\| \bm V - \bm V\bm Y \bm Y^T \| = \sigma_{r+1}. 
\end{equation}
In that case, \cb{$1+\epsilon_f$} is the error factor of the CUR decomposition when compared against the optimal rank-$r$ reduction error obtained by SVD. As demonstrated in our numerical experiments, the TDB-CUR algorithm closely approximates the rank-$r$ SVD approximation of $\bm V$.

\subsection{Oversampling for Improved Condition Number}\label{sec:OS}
The above error analysis shows that the CUR rank-$r$ approximation can be bounded by an error factor $\epsilon_f$ times the maximum error obtained from the orthogonal projection of $\bm V$ onto $\bm U$ or $\bm Y$. This analysis reveals that  better conditioned $\bm P^T \bm U$ and $\bm S^T \bm Y$ matrices result in smaller $\eta_p$ and $\eta_s$, which  then results in smaller error factor $\epsilon_f$. In the context of DEIM interpolation, it was shown that  \emph{oversampling} can improve the condition number of oblique projections \cite{ZW16}. The authors demonstrated that augmenting the original DEIM algorithm with an additional $m=\mathcal{O}(r)$ sampling points can reduce the value of $\eta_p$, leading to smaller approximation errors. This procedure of sampling more rows than the number of basis vectors leads to an overdetermined system where an approximate solution can be found via a least-square regression rather than interpolation. Additionally, it was shown in \cite{anderson2015spectral} that for matrices with rapidly decaying singular values (as targeted in this work), oversampling improves the accuracy of CUR decompositions.


 In the following, we extend the TDB-CUR algorithm for row oversampling. As a direct result of the oversampling procedure, the oblique projection of $\bm V$ onto the range of the orthonormal basis $\bm Q$ becomes:
 \begin{equation}\label{eq:interp}
     \bm Z = \bm Q(\bm p, :)^{\dagger}\bm V(\bm p,:), \quad \mbox{where} \quad \bm Q(\bm p, :)^{\dagger}=(\bm Q(\bm p, :)^T\bm Q(\bm p, :))^{-1}\bm Q(\bm p, :)^T,
 \end{equation}
 and $\bm p \in \mathbb N^{r'}$ contains the $r'=r+m \ll n$ row indices. \cb{Note that  $\bm Q(\bm p, :)^{\dagger}$ is the pseudo-inverse of $\bm Q(\bm p, :)$, however, we do not apply any singular value threshold cutoff to compute $\bm Q(\bm p, :)^{\dagger}$ and exact inversion of $\bm Q(\bm p, :)^T\bm Q(\bm p, :)$ is used.}   Therefore, the oblique projection becomes a least squares best-fit solution. Also, increasing the number of oversampling points decreases $\eta_p = \| \bm Q(\bm p, :)^{\dagger}\|$ and it  follows that for the maximum number of oversampling points, i.e., when all the rows are sampled, the orthogonal projection of every column of $\bm V$ onto $\mathrm{Ran}(\bm V(:,\bm s))$ is recovered, where $\eta_p$ attains its smallest value, which is $\eta_p=1$.  
 Note that, unlike the interpolatory projector, $\bm P^T\bm{P}^{\Oangle}_{\bm{U}} \bm A \neq \bm A(\bm p,:)$. \cb{ The oversampling is also applied analogously to the CUR decomposition of $\bm F$:}
 \begin{equation*}
        \hat{\bm F}_i=\bm U_{\bm F_i}\bm U_{\bm F_i}(\bm p,:)^{\dagger}\bm F_i(\bm p,:), 
    \end{equation*}
    where $\bm U_{\bm F_i}(\bm p,:)^{\dagger}=(\bm U_{\bm F_i}(\bm p,:)^T\bm U_{\bm F_i}(\bm p,:))^{-1}\bm U_{\bm F_i}(\bm p,:)^T$.
The CUR approximation of $\bm F$ is presented in Section \ref{sec:compute-G}.

 We refer to the above sampling procedure as \texttt{OS-DEIM}, where OS refers to the oversampling algorithm. Since the DEIM only provides sampling points equal to the number of basis vectors, we use the GappyPOD+E algorithm from \cite{peherstorfer2020stability}  to sample a total of $r'$ rows. For convenience, the algorithm is provided in Listing \ref{list:gpode}. While any sparse selection procedure can be used, the GappyPOD+E was shown to outperform other common choices like random sampling or leverage scores \cite{mahoney2009cur}. Finally,  it is possible to oversample the columns in an analogous manner to decrease $\eta_s$. In all of the examples considered in this paper, we apply row oversampling, but ultimately the decision for row oversampling, column oversampling, or both may be made by requiring that $\eta_p$ and $\eta_s$ be smaller than some threshold values.

\begin{minipage}{\textwidth}
\renewcommand*\footnoterule{}
\begin{savenotes}
\begin{algorithm}[H]
\caption{Rank-Adaptive TDB-CUR Algorithm\label{alg:S-TDB-ROM}}
\hspace*{\algorithmicindent} \textbf{Input}: $\tilde{\mathbf{U}} \in \mathbb{R}^{n \times \rt}$, $\tilde{\boldsymbol{\Sigma}} \in \mathbb{R}^{\rt \times \rt}$, $\tilde{\mathbf{Y}}\in \mathbb{R}^{s \times \rt}$, $\rt$, $m$ ($\sim$ indicates quantities from previous time step) \\
\hspace*{\algorithmicindent} \textbf{Output}: $\mathbf{U}\in\mathbb R^{n\times r}$, $\boldsymbol{\Sigma}\in\mathbb R^{r\times r}$, $\mathbf{Y}\in\mathbb R^{s\times r}$, $r$
\begin{algorithmic}[1]
\State $\epsilon = \tilde{\bs\Sigma}(\rt,\rt)/\Vert \texttt{diag}(\tilde{\bs\Sigma}) \Vert_2$
\Comment{Compute error proxy for adaptive rank criteria.}
\If{$\epsilon > \epsilon_u$} \Comment{Increase rank if $\epsilon$ exceeds the upper threshold, $\epsilon_u$.}
    \State $r=\rt+1$ 
\ElsIf{$\epsilon < \epsilon_l$}\Comment{Decrease rank if $\epsilon$ falls below the lower threshold, $\epsilon_l$.}
    \State $r = \rt-1$
    \State $\tilde{\bU} = \tilde{\bU}(:,1:r)$;  $\tilde{\bs\Sigma} = \tilde{\bs\Sigma}(1:r,1:r)$; $\tilde{\bY} = \tilde{\bY}(:,1:r)$ \Comment{Truncate TDB matrices.}
\Else\Comment{Keep rank the same.}
    \State $r=\rt$
\EndIf
\State ${\bm s} \gets$ \texttt{sparse\_selection}($\tilde{\mathbf{Y}},r$) \footnote{While the present work uses the GappyPOD+E for \texttt{sparse\_selection}, the user is free to choose their favorite sparse selection algorithm.}\Comment{Compute $r$ column indices.}
\State ${\bm p} \gets$ \texttt{sparse\_selection}($\tilde{\bU},r+m$)
\Comment{Compute $r+m$ row indices.}
\State $\bm p_a\gets$ \texttt{find\_adjacent}($\bm p$)\Comment{Find adjacent points required to compute $\bm V(\bm p,:)$.}
\State $\Vhat(:,\bm s) = \tilde{\mathbf{U}} \tilde{\boldsymbol{\Sigma}} \tilde{\mathbf{Y}}(\bm s,:)^{T}$ \Comment{Construct low-rank approximation of columns in $\bm s$.}
\State $\bm V(:,\bm s) = \cb{\Vhat(:,\bm s)+\Delta t\ol{\bm F}(:,\bm s)}$ \Comment{\cb{Take one step} at the selected columns.}
\State $\Vhat([\bm{p},\bm{p}_a],:) = \tilde{\mathbf{U}}([\bm{p},\bm{p}_a],:)\tilde{\boldsymbol{\Sigma}}\tilde{\mathbf{Y}}^T$\Comment{Construct low-rank approximation of rows in $[\bm p,\bm p_a]$.}
\State $\bm V(\bm{p},:) = \cb{\Vhat(\bm p,:)+\Delta t\ol{\bm F}(\bm p,:)}$\Comment{\cb{Take one step} at the selected rows.}
\State $\bm{Q}\bm{R} = \texttt{QR}(\bm V(:,\bm s),\texttt{`econ'})$\Comment{Compute the economy QR of $\bm V(:,\bm s)$.}
\State $\mathbf{Z}=\mathbf{Q}(\mathbf{p},:)^{\dagger} \bm V(\bm p,:)$\Comment{Compute $\bm Z$ as an oblique projection of $\bm V$ onto $\bm Q$.} 
\State $\mathbf{U}_{\mathbf{Z}} \boldsymbol{\Sigma} \mathbf{Y}^T = \texttt{SVD}(\mathbf{Z},\texttt{`econ'}$)\Comment{Compute the economy SVD of $\bm Z$.}
\State $\mathbf{U}=\mathbf{Q U}_{\mathbf{Z}}$ \Comment{In-subspace rotation of the orthonormal basis, $\bm Q$.}
\end{algorithmic}
\end{algorithm}
\end{savenotes}
\end{minipage}

\subsection{Rank Adaptivity}\label{sec:adaptive-rank}
In order to control the error while avoiding unnecessary computations, the rank of the TDB must be able to adapt on the fly. \cb{The importance of rank adaptivity for low-rank approximation with TDB has been recognized and several algorithms have been proposed recently. See for example \cite{yang2020time, lubich2021rankadaptive, dektor2021rank}}. We show  that it is easy to incorporate mode adaptivity into the TDB-CUR algorithm.
In the case of rank reduction, once the new rank is chosen, such that $r^k < r^{k-1}$, the low-rank matrices are simply truncated to retain only the first $r^k$ components, i.e. $\bU(:,1:r^k)$, $\bs\Sigma(1:r^k,1:r^k)$, and $\bY(:,1:r^k)$. On the other hand, the rank can be increased, such that $r^k>r^{k-1}$, by sampling more columns ($r^k$) than the number of basis vectors ($r^{k-1}$), i.e. oversampling. Similar to the procedure used for oversampling the rows, the column indices are determined via the GappyPOD+E algorithm. While this provides a straightforward approach for \emph{how} to adapt the rank, it does not address \emph{when} the rank should be adapted, or \emph{what} that new rank should be.

Informed by the error analysis from the preceding section, we devise a suitable criterion for controlling the error via rank addition and removal. Since it is not possible to know the true error without solving the expensive FOM, we devise a proxy for estimating the low-rank approximation error: 
\begin{equation}
    \epsilon(t) = \frac{\hat{\sigma}_r(t)}{\left(\sum_{i=1}^r\hat{\sigma}_i(t)^2\right)^{1/2}},
\end{equation}
where $\hat{\sigma}_i$ are the singular values of the low-rank approximation from the previous time step. Assuming the low-rank approximation is near-optimal in its initial condition, we can use the trailing singular value as a proxy for the low-rank approximation error.

To make the error proxy more robust for problems of varying scale and magnitude, we divide by the Frobenius norm of $\Vhat$, where it is well-known that $\Vert \Vhat \Vert_F=(\sum_{i=1}^r\hat{\sigma}_i^2)^{1/2}$. Rather than set a hard threshold, we add/remove modes to maintain $\epsilon$ within a desired range, $\epsilon_l\le\epsilon\le\epsilon_u$, where $\epsilon_l$ and $\epsilon_u$ are user-specified lower and upper bounds, respectively. If $\epsilon>\epsilon_u$ we increase the rank to $r+1$, and if $\epsilon<\epsilon_l$ we decrease the rank to $r-1$. As a result, this approach avoids the undesirable behavior of repeated mode addition and removal, which is observed by setting a hard threshold. The rank-adaptive TDB-CUR algorithm is detailed in Algorithm \ref{alg:S-TDB-ROM}.

It is important to note that this isn't the only criterion for mode addition and removal, and one can devise a number of strategies based on the problem at hand. However, from our numerical experiments, this approach has proved to be simple and effective, and it does a good job at capturing the trend of the true error. For more details on estimating rank and selection criteria, we refer the reader to \cite[Section 2.3]{vidal2005generalized}. \cb{Finally, it is possible to increase the rank by more than one in Algorithm \ref{alg:S-TDB-ROM}, if required. This can be determined by applying the singular value threshold check after executing Line 19. If $\epsilon>\epsilon_u$ is still true, one more column can be sampled. This requires executing Line 10 to find the new column index, Line 11 to update $\bm p$, evaluating $\bm V^k$ only for the new column using Line 14, and following Lines 15-20. These iterations can be carried out many times, until $\epsilon$ falls below $\epsilon_u$}.

\begin{remark}\label{rm:gen_MDE}
    \cb{Algorithm \ref{alg:S-TDB-ROM} is presented for solving MDEs that arise from discretizing PDEs with parametric uncertainties, where the rows are dependent on each other but columns can be solved independently. However, Algorithm \ref{alg:S-TDB-ROM}, with minor modification, can be applied to MDEs where the columns are also dependent on each other. In such cases, evaluating $\bm V^k(:,\bm s)$ requires providing  $\hat{\bm{V}}^{k-1}(:,[\bm s , \bm s_a])= \bm U^{k-1} \bs \Sigma^{k-1} {\bm{Y}^{k-1}}([\bm s , \bm s_a],:)^T$, where $\bm s_a$ is the set of column indices, whose values are needed to compute $\bm V^k(:,\bm s)$. } 
\end{remark}



\begin{remark}\label{rm:dense_MDE}
    \cb{Algorithm \ref{alg:S-TDB-ROM} can be applied to  problems with dense spatial discretizations, where $\bm p_a =[1,2 \dots, n]$. These MDEs can arise, for example, from global discretization methods such as spectral methods. In the most generic form,  the computational complexity of computing each entry of the right-hand side matrix $\bm{F}$ can be $\mathcal{O}(n^{\alpha} s^{\beta} )$ for some $\alpha, \beta \geq 0$. The computational complexity of solving FOM is $ns$ times the computational complexity of each entry, i.e.,  $\mathcal{O}(ns n^{\alpha} s^{\beta} )$ or $\mathcal{O}(n^{(\alpha+1)} s^{(\beta+1)}$. The presented algorithm reduces the cost of evaluating the FOM for this generic setting to $\mathcal{O}((n+s) n^{\alpha} s^{\beta} )$.   } 
\end{remark}

\cb{For MDEs arising from the discretization of PDEs with parametric uncertainties, $\beta=0$, and when sparse discretization schemes are used for the spatial discretization, $\alpha=0$. However, when global discretization methods are used,  $\alpha=1$.  Take for example, $\mathcal{F}(t,\bm V) = \bm D \bm V$, where $\bm D \in \mathbb{R}^{n \times n}$ is a full matrix obtained from discretization of linear differential operators. For this problem, the cost of solving FOM scales with $\mathcal{O}(n^2s)$, while the cost of solving TDB-CUR scales with $\mathcal{O}(n^2+ns)$. The toy problem presented in Section \ref{sec:toy_prob} is a demonstration of a case where there is a dense coupling between both columns and rows.}

\section{Demonstrations}\label{sec:demonstration}

\subsection{Toy Problem}\label{sec:toy_prob}
\begin{figure}[!t]
     \centering
     \subfigure[Full Rank $\bm V(t)$]{
         \centering
%
%
\definecolor{myblue}{rgb}{0.00000,0.44700,0.74100}%
\definecolor{myred}{rgb}{0.85000,0.32500,0.09800}%
\definecolor{myblack}{rgb}{0,0,0}%
\definecolor{mygrey}{rgb}{.7 .7 .7}
\begin{tikzpicture}[scale=0.8]

\begin{axis}[%
width=2.55in,
height=2.35in,
at={(1.046in,0.878in)},
scale only axis,
xmode=log,
xmin=1,
xmax=100,
xminorticks=true,
xlabel style={font=\color{white!15!black}},
xlabel={$r$},
ymode=log,
ymin=1e-16,
ymax=10^2,
yminorticks=true,
ylabel style={font=\color{white!15!black}},
ylabel={$\mathcal{E}(T_f)$},
axis background/.style={fill=white},
legend style={at={(0.01,0.01)}, anchor=south west, legend cell align=left, align=left,font=\scriptsize,draw=white!15!black, font=\scriptsize,fill opacity=0.8, draw opacity=1,
  text opacity=1}
]
\addplot [color=myblue, line width=1.0pt, mark=+, mark options={solid, myblue}]
  table[row sep=crcr]{%
1   0.5\\
2	0.250014622063584\\
4	0.0625745986945715\\
6	0.0159446442370342\\
8	0.00498244100794277\\
10	0.00338125605087416\\
12	0.00319636740587373\\
14	0.00278967392655518\\
16	0.00261271831369576\\
18	0.00307884173791493\\
20	0.00221021325150604\\
22	0.00246926437025475\\
24	0.0107144027633218\\
26	0.00245394711890699\\
28	0.00229832067288488\\
30	0.0022825688345937\\
32	0.00227016381767944\\
34	0.00226882281744376\\
36	0.00226571211835052\\
38	0.00226846933156996\\
40	0.00227083205713417\\
42	0.00226835696707088\\
44	0.00226603654849796\\
46	0.00226560854792322\\
48	0.00226606307180994\\
50	0.00226593695388088\\
52	0.00226604217317113\\
54	0.0022660142567017\\
56	0.00226600913179534\\
58	0.00226605701194188\\
60	0.00226601567172416\\
};
\addlegendentry{TDB-CUR (DEIM) $\mid$ $\Delta t_1$}

\addplot [color=myblue, line width=1.0pt, mark=o, mark options={solid, myblue}]
  table[row sep=crcr]{%
1   0.5\\  
2	0.250008831989557\\
4	0.0625372580964803\\
6	0.0157786201431957\\
8	0.00456880361271387\\
10	0.00257577052949978\\
12	0.00231734952518242\\
14	0.00225624898015126\\
16	0.00221871366116188\\
18	0.00220705399314932\\
20	0.00222566407915653\\
22	0.0022435658607635\\
24	0.00226282210415369\\
26	0.00226315317493618\\
28	0.0022637854859269\\
30	0.00226502700506716\\
32	0.00226521049134609\\
34	0.00226597911295283\\
36	0.00226596735771769\\
38	0.0022659336538412\\
40	0.00226594750500344\\
42	0.00226609524011182\\
44	0.00226595749766555\\
46	0.00226599364472786\\
48	0.00226597984358981\\
50	0.0022660095098956\\
52	0.0022659804436432\\
54	0.00226598787273141\\
56	0.0022659877590048\\
58	0.00226598789078447\\
60	0.00226598733544655\\
};
\addlegendentry{TDB-CUR (OS-DEIM) $\mid$ $\Delta t_1$}

\addplot [color=myblue, line width=1.0pt, mark=triangle, mark options={solid, myblue}]
  table[row sep=crcr]{%
1	0.500000696446737\\
2	0.25000163292868\\
3	0.125003401643888\\
4	0.0625081791301673\\
5	0.0312778898091301\\
6	0.0156783056663922\\
7	0.0679960229016973\\
8	1e5\\
};
\addlegendentry{DLRA Std. Int. $\mid$ $\Delta t_1$ \cite{KL07}}

\addplot [color=myblue, line width=1.0pt, mark=asterisk, mark options={solid, myblue}]
  table[row sep=crcr]{%
1	0.500000078637091\\
2	0.250000263855318\\
3	0.125000599454631\\
4	0.0625018554880956\\
5	0.0312542975765653\\
6	0.0157639187736372\\
7	0.011999735474146\\
8	0.0540332362792798\\
9	0.0869739922552131\\
10	0.521347456398839\\
11	1e5\\
};
\addlegendentry{DO $\mid$ $\Delta t_1$ \cite{SL09}}

\addplot [color=myblue, line width=1.0pt, mark=square, mark options={solid, myblue}]
  table[row sep=crcr]{%
1	0.5\\
2	0.250000269467373\\
4	0.0625015386247087\\
6	0.0156323375244914\\
8	0.00394018826845718\\
10	0.00111981756520385\\
12	0.00068899804086379\\
14	0.00092985320124468\\
16	0.00276372611218331\\
18	0.00356091835850867\\
20	0.00377766946874024\\
22	0.00381680613816122\\
24	0.00382700714539204\\
26	0.00381771809513268\\
28	0.00381693986514617\\
30	0.0038235363265421\\
32	0.00381437782009753\\
34	0.00381891828049093\\
36	0.00381547512093207\\
38	0.00381986260207407\\
40	0.00382668420924098\\
42	0.00383015791765392\\
44	0.00383239347113405\\
46	0.00383510034057348\\
48	0.00384357196711891\\
50	0.003846890568762\\
52	0.00384683109968176\\
54	0.00385117496376006\\
56	0.00385404096868183\\
58	0.00386187723266484\\
60	0.00386287228857146\\
};
\addlegendentry{PS $\mid$ $\Delta t_1$ \cite{LO14}}

\addplot [color=myred, dotted, line width=1.0pt, mark=+, mark options={solid, myred}]
  table[row sep=crcr]{%
1   0.5\\
2	0.25\\
4	0.0624999999999999\\
6	0.015625\\
8	0.00390624999999999\\
10	0.000976562499999994\\
12	0.000244140624999998\\
14	6.10351562500024e-05\\
16	1.52587890625253e-05\\
18	3.81469726571722e-06\\
20	9.53674316807772e-07\\
22	2.3841858070871e-07\\
24	5.9604651109395e-08\\
26	1.49011871123958e-08\\
28	3.72539350192344e-09\\
30	9.3173231853196e-10\\
32	2.34502905334195e-10\\
34	6.45629935358266e-11\\
36	3.139983256644e-11\\
38	2.79347499518271e-11\\
40	2.71796504052288e-11\\
42	2.66686414089951e-11\\
44	2.44794480240596e-11\\
46	2.48827278956319e-11\\
48	2.49185390333839e-11\\
50	2.50866526836541e-11\\
52	2.4927980967831e-11\\
54	2.49778331005293e-11\\
56	2.48548608053878e-11\\
58	2.49494767970496e-11\\
60	2.48722221124103e-11\\
};
\addlegendentry{TDB-CUR (DEIM) $\mid$ $\Delta t_2$}

\addplot [color=myred, dotted, line width=1.0pt, mark=o, mark options={solid, myred}]
  table[row sep=crcr]{%
1   0.5\\
2	0.25\\
4	0.0625\\
6	0.015625\\
8	0.00390625\\
10	0.000976562499999993\\
12	0.000244140625000001\\
14	6.10351562500027e-05\\
16	1.52587890625229e-05\\
18	3.8146972657069e-06\\
20	9.53674316732502e-07\\
22	2.38418580421351e-07\\
24	5.9604650157917e-08\\
26	1.4901183319316e-08\\
28	3.7253793772055e-09\\
30	9.31675117666945e-10\\
32	2.34236219935131e-10\\
34	6.36375794478095e-11\\
36	2.96155279258619e-11\\
38	2.60000933075703e-11\\
40	2.55978199215875e-11\\
42	2.52962336841639e-11\\
44	2.45862866144045e-11\\
46	2.48711941080655e-11\\
48	2.49640445529742e-11\\
50	2.50499661476805e-11\\
52	2.50007337864359e-11\\
54	2.49839613509822e-11\\
56	2.48919244951643e-11\\
58	2.48908036195595e-11\\
60	2.48351282123606e-11\\
};
\addlegendentry{TDB-CUR (OS-DEIM) $\mid$ $\Delta t_2$}

\addplot [color=myred, dotted, line width=1.0pt, mark=triangle, mark options={solid, myred}]
  table[row sep=crcr]{%
1	0.5\\
2	0.25\\
3	0.125\\
4	0.0625\\
5	0.03125\\
6	0.015625\\
7	0.00781249999999999\\
8	0.00390624999999999\\
9	0.00195312500000005\\
10	0.000976562500000151\\
11	0.000488281250018764\\
12	0.00024414062623655\\
13	0.000122070324892482\\
14	6.10351881486726e-05\\
15	3.100506318982e-05\\
16	2.80418636045631e-05\\
17	0.000160596552782842\\
18	0.00635366412880544\\
19	0.0745419909990391\\
20	0.555406042484047\\
21	1.27481150305247\\
22	1.28290541722276\\
23	1.20812999596956\\
24	1e5\\
};
\addlegendentry{DLRA Std. Int. $\mid$ $\Delta t_2$ \cite{KL07}}

\addplot [color=myred, dotted, line width=1.0pt, mark=asterisk, mark options={solid, myred}]
  table[row sep=crcr]{%
1	0.5\\
2	0.25\\
3	0.125\\
4	0.0625\\
5	0.03125\\
6	0.015625\\
7	0.00781249999999998\\
8	0.00390625\\
9	0.00195312500000003\\
10	0.000976562500000047\\
11	0.000488281250090073\\
12	0.000244140627212022\\
13	0.000122070528048094\\
14	6.10357336957267e-05\\
15	3.06300307315672e-05\\
16	1.57866302910535e-05\\
17	7.77374448117332e-05\\
18	0.000856562487694364\\
19	0.0506045744818511\\
20	0.638848672182523\\
21	0.171903003852237\\
22	0.307764022593519\\
23	1.04841547587368\\
24	1.07476940185158\\
25	1.25896475114826\\
26	0.977380596168923\\
27	0.940763299576873\\
28	1.44712557655873\\
29	1.42010646373574\\
30	1.3311442141583\\
31	1.36510046719813\\
32	1.42995900812233\\
33	1.3140103696187\\
34	1.45134081247101\\
35	1.3356451601476\\
36	1.53566933827227\\
37	1.49918706544435\\
38	1.32662043085303\\
39	1.32982733402042\\
40	1.32415477201334\\
41	1.39110554936323\\
42	1e5\\
};
\addlegendentry{DO $\mid$ $\Delta t_2$ \cite{SL09}}

\addplot [color=myred, dotted,line width=1.0pt, mark=square, mark options={solid, myred}]
  table[row sep=crcr]{%
1	0.5\\
2	0.25\\
4	0.0624999999999999\\
6	0.015625\\
8	0.00390624999999999\\
10	0.000976562499999995\\
12	0.000244140625000002\\
14	6.10351562499967e-05\\
16	1.52587890625065e-05\\
18	3.81469726563683e-06\\
20	9.53674316454316e-07\\
22	2.38418579315985e-07\\
24	5.96046458114969e-08\\
26	1.49011654900203e-08\\
28	3.72530870393992e-09\\
30	9.31401076400462e-10\\
32	2.33199920783727e-10\\
34	5.97767603475955e-11\\
36	2.02705934756652e-11\\
38	1.48276575071819e-11\\
40	1.72771370967761e-11\\
42	2.02425801490246e-11\\
44	2.33016423527791e-11\\
46	2.42894253743222e-11\\
48	2.45949632625638e-11\\
50	2.47067600299958e-11\\
52	2.47507458967689e-11\\
54	2.47631850602666e-11\\
56	2.47741649491975e-11\\
58	2.47828391734066e-11\\
60	2.47797166578047e-11\\
};
\addlegendentry{PS $\mid$ $\Delta t_2$ \cite{LO14}}

\addplot [color=black, line width=1.0pt]
  table[row sep=crcr]{%
1   0.5\\
2	0.25\\
4	0.0625\\
6	0.015625\\
8	0.00390625\\
10	0.000976562499999992\\
12	0.000244140624999997\\
14	6.10351562499962e-05\\
16	1.52587890625035e-05\\
18	3.81469726562688e-06\\
20	9.53674316403435e-07\\
22	2.38418579100398e-07\\
24	5.96046447736548e-08\\
26	1.49011611974322e-08\\
28	3.72529030499035e-09\\
30	9.31322574218637e-10\\
32	2.32830644280305e-10\\
34	5.8207657118119e-11\\
36	1.45519185365298e-11\\
38	3.63798548449969e-12\\
40	9.09499084927673e-13\\
42	2.27392606650626e-13\\
44	5.69113874393991e-14\\
46	1.44858964307182e-14\\
48	4.5290281899017e-15\\
50	2.96095138133268e-15\\
52	2.90967194029906e-15\\
54	2.88498981542425e-15\\
56	2.88681931726308e-15\\
58	2.88686308774316e-15\\
60	2.88807662323677e-15\\
};
\addlegendentry{\texttt{SVD}$(\bm V(T_f))$}
\end{axis}
\end{tikzpicture}%
         \label{fig:Error_r}
     }
     \subfigure[Rank-Deficient $\bm V(t)$]{
         \centering
%
%
\definecolor{myblue}{rgb}{0.00000,0.44700,0.74100}%
\definecolor{myred}{rgb}{0.85000,0.32500,0.09800}%
\definecolor{myblack}{rgb}{0,0,0}%
\definecolor{mygrey}{rgb}{.7 .7 .7}
\begin{tikzpicture}[scale=0.8]

\begin{axis}[%
width=2.55in,
height=2.35in,
at={(1.136in,0.878in)},
scale only axis,
xmode=log,
xmin=1,
xmax=100,
xminorticks=true,
xlabel style={font=\color{white!15!black}},
xlabel={$r$},
ymode=log,
ymin=1e-16,
ymax=1e2,
yminorticks=true,
ylabel style={font=\color{white!15!black}},
ylabel={$\mathcal{E}(T_f)$},
axis background/.style={fill=white},
legend style={at={(0.01,0.01)}, anchor=south west, legend cell align=left, align=left,font=\scriptsize,draw=white!15!black, font=\scriptsize,fill opacity=0.8, draw opacity=1,
  text opacity=1}
]

\addplot [color=myblue, line width=1.0pt, mark=triangle, mark options={solid, myblue}]
  table[row sep=crcr]{%
1	0.49926702208453\\
2	0.248162033997406\\
3	0.12109338543896\\
4	0.0541624857960875\\
5	0.00132210457039386\\
6	1e5\\
};

\addplot [color=myblue, line width=1.0pt, mark=asterisk, mark options={solid, myblue}]
  table[row sep=crcr]{%
1	0.499266402705131\\
2	0.248160653114511\\
3	0.121090489248597\\
4	0.0541551794310088\\
5	0.000518505023722816\\
6	1e5\\
};
\addplot [color=myblue, line width=1.0pt, mark=+, mark options={solid, myblue}]
  table[row sep=crcr]{%
1    0.499269642586741
2	0.248175131200938\\
4	0.0542391717823277\\
6	0.00270544575701949\\
8	0.00222355422825159\\
10	0.00227306032028547\\
12	0.00227467345758884\\
14	0.00225814458058512\\
16	0.00226384890296631\\
18	0.00226157817177852\\
20	0.00226370147738314\\
22	0.00226144369398602\\
24	0.0022615074311466\\
26	0.00226149150953558\\
28	0.00226149123456083\\
30	0.00226149132643113\\
32	0.00226149120296967\\
34	0.00226149121097904\\
36	0.00226149120830836\\
38	0.00226149120847581\\
40	0.00226149120847295\\
42	0.00226149120847842\\
44	0.00226149120847983\\
46	0.00226149120847639\\
48	0.00226149120847684\\
50	0.00226149120847593\\
52	0.00226149120847653\\
54	0.00226149120847627\\
56	0.00226149120847702\\
58	0.00226149120847702\\
60	0.00226149120847673\\
};

\addplot [color=myblue, line width=1.0pt, mark=o, mark options={solid, myblue}]
  table[row sep=crcr]{%
1   0.499270111211453 
2	0.248169290744527\\
4	0.0541960548097126\\
6	0.00213633439822848\\
8	0.00223310801771008\\
10	0.00225448633080751\\
12	0.00226133778248204\\
14	0.00225833516910828\\
16	0.00226231913913508\\
18	0.00226199080215297\\
20	0.00226156025137857\\
22	0.00226143959560398\\
24	0.00226148219735061\\
26	0.00226149124510039\\
28	0.00226149114537793\\
30	0.00226149119718375\\
32	0.00226149120306432\\
34	0.00226149120727232\\
36	0.00226149120851055\\
38	0.00226149120850158\\
40	0.00226149120847489\\
42	0.00226149120847502\\
44	0.00226149120847712\\
46	0.00226149120847658\\
48	0.00226149120847698\\
50	0.00226149120847573\\
52	0.00226149120847734\\
54	0.00226149120847626\\
56	0.00226149120847617\\
58	0.00226149120847594\\
60	0.00226149120847627\\
};

\addplot [color=myblue, line width=1.0pt, mark=square, mark options={solid, myblue}]
  table[row sep=crcr]{%
1	0.499267\\
2	0.248160658941626\\
4	0.0541548142194864\\
6	0.00168574373452664\\
8	0.00364599331381577\\
10	0.00390912013892031\\
12	0.00392378918375546\\
14	0.00386349921659771\\
16	0.00381079038664793\\
18	0.00380556750599535\\
20	0.00380627464285453\\
22	0.00380549436337858\\
24	0.00380339394480853\\
26	0.00380212571932838\\
28	0.00380103724026982\\
30	0.00380166233039213\\
32	0.00380449651178172\\
34	0.00380717059418838\\
36	0.0038106924006401\\
38	0.00381324651259358\\
40	0.00381572153419904\\
42	0.00381952448961759\\
44	0.00382365164631498\\
46	0.00382777853985452\\
48	0.00383487295166003\\
50	0.00384103600608518\\
52	0.0038421173777326\\
54	0.00384786303131014\\
56	0.00385046023023474\\
58	0.00385786260702947\\
60	0.00385884210800649\\
};
\addplot [color=myred, dotted, line width=1.0pt, mark=triangle, mark options={solid, myred}]
  table[row sep=crcr]{%
1	0.499266323889453\\
2	0.248160387073783\\
3	0.121089869924121\\
4	0.0541530361073882\\
5	1.24411507266825e-11\\
6	9264.94054690419\\
7	558395.436452626\\
8	108517.099417521\\
};

\addplot [color=myred, dotted, line width=1.0pt, mark=asterisk, mark options={solid, myred}]
  table[row sep=crcr]{%
1	0.499266323889453\\ 
2	0.248160387073783\\
3	0.121089869924121\\
4	0.0541530361073882\\
5	5.10873175405064e-12\\
6	0.910184622848508\\
7	1.25238226010956\\
8	1.30314109943758\\
9	1.30731745735103\\
10	1.36867119042206\\
11	10000.36867119042206\\
};

\addplot [color=myred, dotted, line width=1.0pt, mark=+, mark options={solid, myred}]
  table[row sep=crcr]{%
1   0.499266323889453
2	0.248160387073783\\
4	0.0541530361073882\\
6	2.34408121695341e-11\\
8	2.43330948673613e-11\\
10	2.47058807943826e-11\\
12	2.47792685840225e-11\\
14	2.47270654451993e-11\\
16	2.47551638146133e-11\\
18	2.47602837060423e-11\\
20	2.48031255729571e-11\\
22	2.47491135243931e-11\\
24	2.4739901261674e-11\\
26	2.47455718250266e-11\\
28	2.47533838695569e-11\\
30	2.47546377046957e-11\\
32	2.47537397497779e-11\\
34	2.47534394859275e-11\\
36	2.47514514837317e-11\\
38	2.47535849887357e-11\\
40	2.47552872734603e-11\\
42	2.47525036646773e-11\\
44	2.47616760214551e-11\\
46	2.47544792080542e-11\\
48	2.4759319073132e-11\\
50	2.47534809651454e-11\\
52	2.47512826576384e-11\\
54	2.47571942737807e-11\\
56	2.47577038121125e-11\\
58	2.47615168311593e-11\\
60	2.47495510782971e-11\\
};

\addplot [color=myred, dotted, line width=1.0pt, mark=o, mark options={solid, myred}]
  table[row sep=crcr]{%
1   0.499266323889453
2	0.248160387073783\\
4	0.0541530361073882\\
6	2.16535100784877e-11\\
8	2.43403012696605e-11\\
10	2.46425802522662e-11\\
12	2.47491741489857e-11\\
14	2.47321841622044e-11\\
16	2.47396555314821e-11\\
18	2.47484801075246e-11\\
20	2.47501698723893e-11\\
22	2.4736687196361e-11\\
24	2.47489386006761e-11\\
26	2.47466863461052e-11\\
28	2.47585567826684e-11\\
30	2.47554130249564e-11\\
32	2.47481822165157e-11\\
34	2.47551210387519e-11\\
36	2.47439715263516e-11\\
38	2.47532424528891e-11\\
40	2.47549689537045e-11\\
42	2.4758168286253e-11\\
44	2.47534308815266e-11\\
46	2.47614866432038e-11\\
48	2.47559239889562e-11\\
50	2.47675317384081e-11\\
52	2.47536452432394e-11\\
54	2.47571955673479e-11\\
56	2.47511368359062e-11\\
58	2.47599737698078e-11\\
60	2.47547276248087e-11\\
};
\addplot [color=myred, dotted, line width=1.0pt, mark=square, mark options={solid, myred}]
  table[row sep=crcr]{%
1	0.499266\\
2	0.248160387073783\\
4	0.0541530361073882\\
6	1.68275947306279e-11\\
8	2.36739421280324e-11\\
10	2.45242097601603e-11\\
12	2.46631777611775e-11\\
14	2.47043396542112e-11\\
16	2.47171301575667e-11\\
18	2.47264621639209e-11\\
20	2.47259756392734e-11\\
22	2.47290583765102e-11\\
24	2.47283167611167e-11\\
26	2.47292338107307e-11\\
28	2.47292973348994e-11\\
30	2.47268142012761e-11\\
32	2.47284775253655e-11\\
34	2.47291275351944e-11\\
36	2.472592533718e-11\\
38	2.47301892260742e-11\\
40	2.47290958263794e-11\\
42	2.47282842497492e-11\\
44	2.47291722341242e-11\\
46	2.47287749422163e-11\\
48	2.47291456512603e-11\\
50	2.4727744945451e-11\\
52	2.47283854661177e-11\\
54	2.47265880842177e-11\\
56	2.47296089276576e-11\\
58	2.47281617203641e-11\\
60	2.47297878354662e-11\\
};

\addplot [color=myblack, line width=1.0pt]
  table[row sep=crcr]{%
1	0.499266323889453\\
2	0.248160387073783\\
3	0.121089869924121\\
4	0.0541530361073883\\
5	8.44327280934686e-16\\
6	9.4842881140901e-16\\
7	1.03700147204602e-15\\
8	1.10584596002282e-15\\
9	1.15228359012857e-15\\
10	1.19252119990134e-15\\
11	1.2186255052603e-15\\
12	1.24035679219235e-15\\
13	1.25630301289187e-15\\
14	1.26697818627568e-15\\
15	1.27415990318291e-15\\
16	1.27990113128385e-15\\
17	1.28211440559709e-15\\
18	1.28378752416426e-15\\
19	1.28540654102983e-15\\
20	1.28662825331042e-15\\
21	1.28847022286999e-15\\
22	1.29032017848579e-15\\
23	1.29165254442079e-15\\
24	1.29322991565485e-15\\
25	1.29496104131002e-15\\
26	1.29622103005832e-15\\
27	1.29795994192233e-15\\
28	1.29941159776484e-15\\
29	1.30093842091625e-15\\
30	1.30260570414428e-15\\
31	1.30412986108336e-15\\
32	1.30598938261325e-15\\
33	1.30736516679461e-15\\
34	1.30904941162344e-15\\
35	1.31042404159512e-15\\
36	1.31190783062718e-15\\
37	1.3135520702423e-15\\
38	1.31499102461573e-15\\
39	1.31648211232327e-15\\
40	1.31846512152763e-15\\
41	1.32027799096689e-15\\
42	1.32157666572787e-15\\
43	1.32317251049739e-15\\
44	1.32481960894249e-15\\
45	1.3265335878612e-15\\
46	1.32798326602698e-15\\
47	1.3291689997008e-15\\
48	1.33060529080315e-15\\
49	1.3320498804117e-15\\
50	1.33385965237793e-15\\
51	1.33548758954421e-15\\
52	1.33727754732359e-15\\
53	1.33876127607062e-15\\
54	1.34022978417617e-15\\
55	1.3418891648639e-15\\
56	1.34316808624404e-15\\
57	1.3448564971266e-15\\
58	1.3464130500214e-15\\
59	1.34781696311115e-15\\
60	1.34941297136284e-15\\
};
\end{axis}

\end{tikzpicture}%
         \label{fig:Error_r_RD}
    
    }
     \subfigure[Full Rank $\bm V(t)$]{
         \centering
%
%
\definecolor{myblue}{rgb}{0.00000,0.44700,0.74100}%
\definecolor{myred}{rgb}{0.85000,0.32500,0.09800}%
\definecolor{myblack}{rgb}{0,0,0}%
\definecolor{mygrey}{rgb}{.7 .7 .7}

\begin{tikzpicture}[scale=0.8]

\begin{axis}[%
width=2.55in,
height=2.35in,
at={(1.046in,0.794in)},
scale only axis,
xmin=0,
xmax=90,
xlabel style={font=\color{white!15!black}},
xlabel={$r$},
ymode=log,
ymin=1,
ymax=1e+10,
yminorticks=true,
ylabel style={font=\color{white!15!black}},
ylabel={Norm of inverse},
axis background/.style={fill=white},
legend style={at={(0.1,0.65)}, anchor=south west, legend cell align=left, align=left, draw=white!15!black,font=\scriptsize}
]
\addplot [color=mygrey, line width=1.0pt, mark=o, mark options={solid, mygrey}]
  table[row sep=crcr]{%
1	0.73603318146115\\
2	1.47186853878068\\
3	2.94347565764296\\
4	5.88422195079267\\
5	11.7629514462931\\
6	23.5501493719503\\
7	51.1139985857867\\
8	1.25052278851982e+19\\
};
\addlegendentry{DLRA Std. Int. \cite{KL07} ($\|\bs \Sigma^{-1}\|$)}

\addplot [color=myblack, line width=1.0pt, mark=o, mark options={solid, myblack}]
  table[row sep=crcr]{%
1	0.541504293941878\\
2	2.16595024506351\\
3	8.66283858906197\\
4	34.6594733912378\\
5	138.615674815279\\
6	543.693447776471\\
7	1926.36320078083\\
8	8011.49138982699\\
9	8784.88557648143\\
10	19507.2791720504\\
11	5.12442939369175e+18\\
};
\addlegendentry{DO \cite{SL09} ($\|\bm C^{-1}\|$)}

\addplot [color=myred, line width=1.0pt, mark=o, mark options={solid, myred}]
  table[row sep=crcr]{%
2	3.53816688278307\\
4	4.64627275781365\\
6	5.14873704231453\\
8	6.10430933986175\\
10	5.94537424328696\\
12	8.27854096383219\\
14	6.24427094565775\\
16	9.83514048967725\\
18	7.5537615875938\\
20	10.1311878547535\\
22	14.1199380216532\\
24	9.58132195321378\\
26	10.7293584148472\\
28	9.35375025012918\\
30	10.9949469833999\\
32	11.1292172402793\\
34	8.35207084267975\\
36	10.2157382580018\\
38	11.3786343292542\\
40	8.17543907267701\\
42	13.284458044755\\
44	12.6769954416536\\
46	10.4752051446123\\
48	17.3506466742657\\
50	13.0582460674774\\
52	13.716757265101\\
54	11.1285589941674\\
56	19.8175932748377\\
58	12.1638666093191\\
60	10.4057670479674\\
62	19.3271686289627\\
64	16.2750980344182\\
66	11.0369522875218\\
68	10.7846063707609\\
70	7.59619990161686\\
72	19.4805044976579\\
74	9.69231522251791\\
76	10.8180746567843\\
78	8.71295450347796\\
80	9.34464204425723\\
82	11.0545064266121\\
84	11.1843403757053\\
86	7.82128127049525\\
88	6.56252098483527\\
90	8.37131213763082\\
92	11.1847539534581\\
94	5.71987249247283\\
96	6.71608686464704\\
98	2.68543448833929\\
100	1\\
};
\addlegendentry{TDB-CUR (DEIM) ($\|\bm Q(\bm p,:)^{-1}\|$)}

\addplot [color=myblue, line width=1.0pt, mark=o, mark options={solid, myblue}]
  table[row sep=crcr]{%
2	1.86575028329244\\
4	2.11665553452682\\
6	2.15385971594223\\
8	2.2931001926474\\
10	2.57379800407065\\
12	2.96453444286989\\
14	2.78337922497049\\
16	2.96901740434083\\
18	2.71468420203663\\
20	2.75996306101009\\
22	2.97370539776708\\
24	2.89912145164454\\
26	3.29422794290857\\
28	3.13296219674791\\
30	3.2735059555372\\
32	2.98629749631875\\
34	2.90896112645501\\
36	3.04322064873622\\
38	2.92507652438221\\
40	3.33471284886125\\
42	3.21127805568218\\
44	2.95538749045638\\
46	2.95372116702701\\
48	3.29378007102947\\
50	3.11664416307735\\
52	3.00817134945009\\
54	2.93876846392022\\
56	3.16916994074479\\
58	3.07778028979453\\
60	2.94030238643\\
62	2.97614780107035\\
64	3.04435740438166\\
66	3.14832215392614\\
68	2.80398083023845\\
70	2.79154498691963\\
72	2.73069403757398\\
74	2.78776203780961\\
76	2.88856143767799\\
78	2.83995473224367\\
80	2.5498410019842\\
82	2.5059247698496\\
84	2.35767501191358\\
86	2.37409042227793\\
88	1.87799415847592\\
90	1\\
};
\addlegendentry{TDB-CUR (OS-DEIM) ($\|\bm Q(\bm p,:)^{\dagger}\|$)}
\end{axis}
\end{tikzpicture}%
         \label{fig:CN_r}
     }
     \subfigure[Full Rank $\bm V(t)$]{
         \centering
%
%
\definecolor{myblue}{rgb}{0.00000,0.44700,0.74100}%
\definecolor{myred}{rgb}{0.85000,0.32500,0.09800}%
\definecolor{myblack}{rgb}{0,0,0}%
\definecolor{mygrey}{rgb}{.7 .7 .7}
\begin{tikzpicture}[scale=0.8]

\begin{axis}[%
width=2.55in,
height=2.35in,
at={(1.011in,0.642in)},
scale only axis,
xmode=log,
xmin=0.0001,
xmax=1,
xminorticks=true,
xlabel style={font=\color{white!15!black}},
xlabel={$\Delta t$},
ymode=log,
ymin=1.0e-17,
ymax=100,
yminorticks=true,
ylabel style={font=\color{white!15!black}},
ylabel={$\mathcal{E}(T_f)$},
axis background/.style={fill=white},
legend style={at={(0.1,0.0001)}, anchor=south west, legend cell align=left, align=left,font=\scriptsize, fill opacity=0.8, draw opacity=1,
  text opacity=1}
]

\addplot [color=myblue, line width=1.0pt, mark=+, mark options={solid, myblue}]
  table[row sep=crcr]{%
0.5	0.783474529336291\\
0.25	0.614646526128554\\
0.125	0.397596249088415\\
0.0625	0.224288105536789\\
0.03125	0.119527052728295\\
0.015625	0.0618328209133949\\
0.0078125	0.0316109351107923\\
0.00390625	0.0162906440975639\\
0.001953125	0.00884956444939056\\
0.0009765625	0.00557573999939009\\
0.00048828125	0.00438458725705676\\
0.000244140625	0.00403128789218865\\
0.0001220703125	0.00393789820986087\\
};
\addlegendentry{DLRA Unconv. Integ. $\mid$ $r = 8$ \cite{ceruti2021unconventional}}

\addplot [color=myblue, line width=1.0pt, mark=o, mark options={solid, myblue}]
  table[row sep=crcr]{%
0.5	1.46262507794608\\
0.25	0.0734387793930567\\
0.125	0.00709675899098937\\
0.0625	0.00392242505158827\\
0.03125	0.00390630692668625\\
0.015625	0.00390625021700442\\
0.0078125	0.00390625000083653\\
0.00390625	0.00390625000000325\\
0.001953125	0.00390625000000001\\
0.0009765625	0.00390624999999999\\
0.00048828125	0.00390624999999999\\
0.000244140625	0.00390624999999999\\
0.0001220703125	0.00390624999999999\\
};
\addlegendentry{TDB-CUR (OS-DEIM) $\mid$ $r = 8$}

\addplot [color=myred, line width=1.0pt, mark=+, mark options={solid, myred}]
  table[row sep=crcr]{%
0.5	1.79475823571256\\
0.25	0.212320623425001\\
0.125	0.179810607983951\\
0.0625	0.178106045084809\\
0.03125	0.112041082934767\\
0.015625	0.0579006784979872\\
0.0078125	0.0294339595561399\\
0.00390625	0.0148404004190946\\
0.001953125	0.00745139042698502\\
0.0009765625	0.00373355436624847\\
0.00048828125	0.00186878801078693\\
0.000244140625	0.000934979056688481\\
0.0001220703125	0.000467799193363757\\
};
\addlegendentry{DLRA Unconv. Integ. $\mid$ $r = 16$ \cite{ceruti2021unconventional}}

\addplot [color=myred, line width=1.0pt, mark=o, mark options={solid, myred}]
  table[row sep=crcr]{%
0.5	1.58376688284777\\
0.25	0.0689678078488772\\
0.125	0.00522608984961341\\
0.0625	0.000353255605720409\\
0.03125	2.71755319476723e-05\\
0.015625	1.5321358228045e-05\\
0.0078125	1.52590322299146e-05\\
0.00390625	1.52587900091391e-05\\
0.001953125	1.52587890662025e-05\\
0.0009765625	1.52587890625151e-05\\
0.00048828125	1.52587890625023e-05\\
0.000244140625	1.52587890625133e-05\\
0.0001220703125	1.52587890625414e-05\\
};
\addlegendentry{TDB-CUR (OS-DEIM) $\mid$ $r = 16$}

\addplot [color=myblack, line width=1.0pt, mark=+, mark options={solid, myblack}]
  table[row sep=crcr]{%
0.5	1.59455902202364\\
0.25	0.0791366571179818\\
0.125	0.0263258133373095\\
0.0625	0.0192811308766461\\
0.03125	0.0141228404498425\\
0.015625	0.0115469947060962\\
0.0078125	0.0103672680239765\\
0.00390625	0.00961146158160766\\
0.001953125	0.00602989052653142\\
0.0009765625	0.00303778071884244\\
0.00048828125	0.00152082299481097\\
0.000244140625	0.000760790079907968\\
0.0001220703125	0.000380486486907731\\
};
\addlegendentry{DLRA Unconv. Integ. $\mid$ $r = 32$ \cite{ceruti2021unconventional}}

\addplot [color=myblack, line width=1.0pt, mark=o, mark options={solid, myblack}]
  table[row sep=crcr]{%
0.5	1.57867893858231\\
0.25	0.0691928622459301\\
0.125	0.0053721054453534\\
0.0625	0.00035881121480422\\
0.03125	2.30306142040398e-05\\
0.015625	1.44783230302827e-06\\
0.0078125	9.00176140284007e-08\\
0.00390625	5.80767042397349e-09\\
0.001953125	4.36831678602252e-10\\
0.0009765625	2.33999803860202e-10\\
0.00048828125	2.32897896416795e-10\\
0.00048828125	2.32898083165861e-10\\
0.000244140625	2.33079840349673e-10\\
0.0001220703125	2.33822480825768e-10\\
};
\addplot [color=mygrey, forget plot]
  table[row sep=crcr]{%
0.0625	0.000358993935814926\\
0.03125	0.000358993935814926\\
0.03125	2.30515005798464e-05\\
};
\node[right, align=left, inner sep=0]
at (axis cs:0.004,9.0e-5) {Slope = 4};

\addplot [color=mygrey, forget plot]
  table[row sep=crcr]{%
0.015625	0.0618328209133949\\
0.0078125	0.0618328209133949\\
0.0078125	0.0316109351107923\\
};
\addlegendentry{TDB-CUR (OS-DEIM) $\mid$ $r = 32$}
\node[right, align=left, inner sep=0]
at (axis cs:0.003,0.5) {Slope = 1};
\end{axis}
\end{tikzpicture}%
         \label{fig:Error_dt}
     }
        \caption{Toy Problem: (a)-(b) Error at the final time $\mathcal{E}(T_f)$ versus reduction order $r$ for step-sizes $\Delta t_1 = 10^{-1} $ and $\Delta t_2 = 10^{-3} $ (c) $l^2$ norm of the inverse matrix versus reduction order $r$ (d) $\mathcal{E}(T_f)$ versus step size $\Delta t$ for various reduction orders $r = 8, 16, 32$. }
        \label{fig:Toy_problem}
\end{figure}
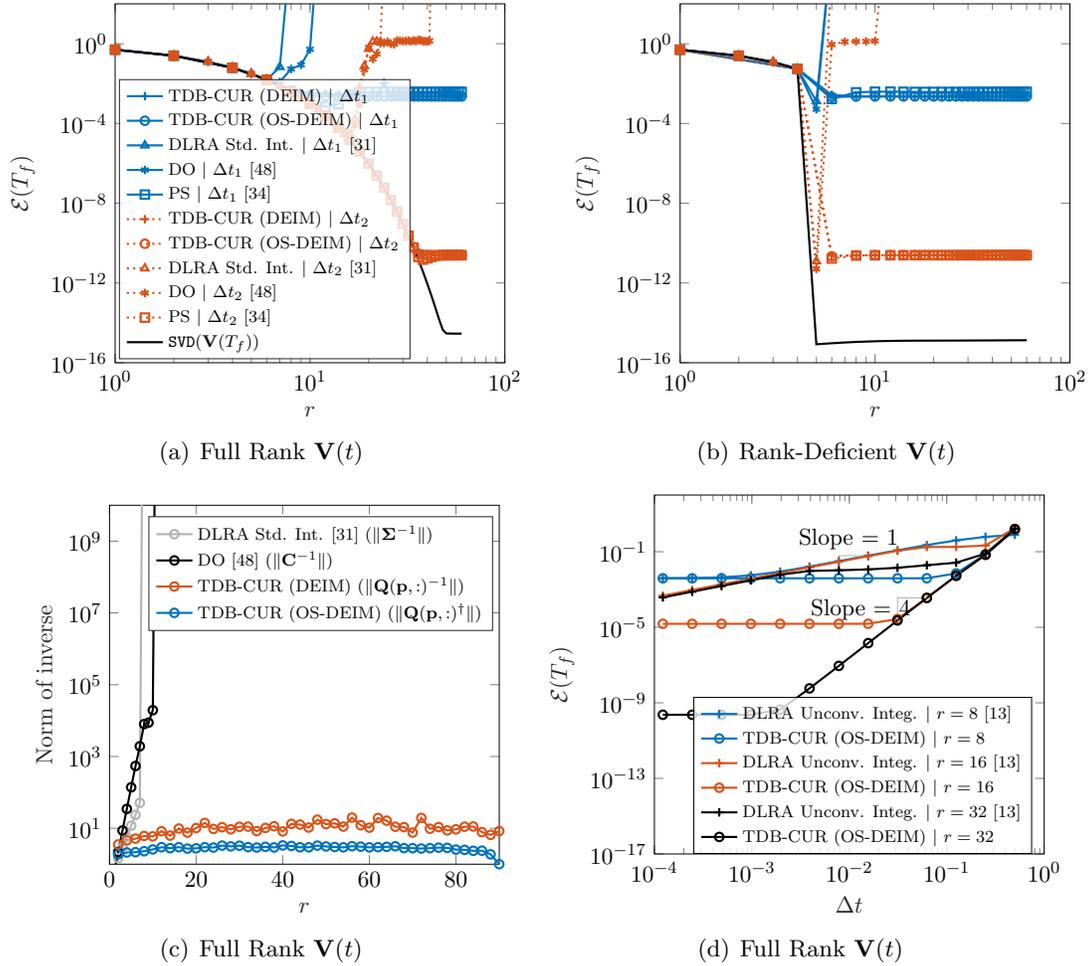
As our first example, we compare the accuracy of the presented algorithm against DLRA \cb{using standard integrator} \cite{KL07}, DO \cite{SL09}, \cb{the projector splitting time integrator (PS) \cite{LO14}} and the recently proposed unconventional robust integrator \cite{ceruti2021unconventional}. \cb{We emphasize that it is already established that the standard integrator, for example, Runge Kutta schemes, are unstable for solving Eqs. \ref{eq:DBO_evol_S} - \ref{eq:DBO_evol_Y}, which has motivated the development of new time integration techniques \cite{LO14,ceruti2021unconventional}.} We consider the time-dependent matrix from \cite{ceruti2021unconventional} given explicitly as
\begin{equation*}
    \bm V(t) = (e^{t\bm W_1})e^t\bm D(e^{t\bm W_2})^T,  \;\; 0\leq t\leq 1.
\end{equation*}
The matrices $\bm W_1\in \mathbb{R}^{n\times n}$ and $ \bm W_2\in \mathbb{R}^{n\times n}$ are randomly generated skew-symmetric matrices as follows: $\bm W_1 = (\tilde{\bm W}_1-\tilde{\bm W}_1^T)/2$ and $\bm W_2 = (\tilde{\bm W}_2-\tilde{\bm W}_2^T)/2$, where $\tilde{\bm W}_1 \in \mathbb{R}^{n \times n}$ and $\tilde{\bm W}_2 \in \mathbb{R}^{n \times n}$ are uniformly distributed random matrices. The matrix $\bm D\in \mathbb{R}^{n\times n}$ is diagonal with diagonal entries $d_i = 1/2^i$ for $i\in \{1,2,...,n\}$. We choose $n=100$ and final time $T_f = 1$. We create a linear MDE: $\D \bm V/\D t=\bm W_1\bm V + \bm V + \bm V\bm W_2^T$, where the right-hand side is linearly dependent on $\bm V$. We use the explicit fourth-order Runge-Kutta integrator for all of the methods including the substeps of the unconventional robust integrator \cite{ceruti2021unconventional}.  We use relative Frobenius error $\mathcal{E}^k = \| \hat{\bm V}^k - \bm V^k\|_{F}/\| \bm V^k\|_{F}$ in our analysis.

In Figure \ref{fig:Error_r}, we plot the error at the final time $\mathcal{E}(T_f)$ versus rank for two different step-sizes: $\Delta t_1 = 10^{-1} $ and $\Delta t_2 = 10^{-3} $. We consider the following cases:  TDB-CUR (DEIM), TDB-CUR (OS-DEIM) with $m=10$, DO \cite{SL09} (Eqs. \ref{eqn:u-do}-\ref{eqn:y-do}), DLRA \cb{using standard integrator} \cite{KL07} (Eqs. \ref{eq:DBO_evol_S}-\ref{eq:DBO_evol_Y}), \cb{ and PS \cite{LO14}}. For reference, we also show the optimal error that is obtained via the  rank-$r$ SVD of the exact  solution at the final time $T_f$, denoted by \texttt{SVD}($\bm V(T_f)$).
For $\Delta t_1$, both DO  and DLRA diverge before reaching $r=10$. This is because the matrix $\bs\Sigma$ in DLRA and the matrix $\bm C$ in DO become poorly conditioned as $r$ increases. However, \cb{PS}, TDB-CUR (OS-DEIM), and TDB-CUR (OS-DEIM) follow the optimal error until the temporal integration error dominates, at which point the error cannot be reduced further by increasing $r$. It is worth noting that without oversampling, TDB-CUR has a sudden increase in error at around $r=25$. However, this undesirable behavior is eliminated by oversampling. As the time step is reduced to $\Delta t_2$, we observe a corresponding decrease in the \cb{PS}, TDB-CUR (DEIM) and TDB-CUR (OS-DEIM) errors. Although DO and DLRA still diverge for the smaller time step, this occurs at a much larger value of $r$. Thus, Figure \ref{fig:Error_r} also highlights the severe time step restrictions for the stability of DO and DLRA in the presence of small singular values.

In Figure \ref{fig:Error_r_RD}, we consider a rank-deficient matrix and overapproximation using different low-rank techniques. Specifically, we consider the matrix $\bm D\in \mathbb{R}^{n\times n}$, where the diagonal entries are given by $d_i = 1/2^i$ for $i\in {1,\ldots,5}$, and all remaining entries are zero. Our results show that the error drops to the optimal temporal error when $r=5$. This is because the rank of the matrix is $5$. Furthermore, we observe that even in the case of rank deficiency, TDB-CUR remains stable for $r>5$. This finding supports the observation made in Section \ref{sec:sparse-sampling} regarding the conditioning of the presented algorithm. Specifically, even when $r>5$, the matrix $\bm Q(\bm p,:)$ is well-conditioned, and the TDB-CUR scheme remains stable, while both DO and DLRA \cb{with standard integrator} diverge. \cb{Similarly, PS remains stable and converges to the optimal error as it does not require inverting $\bs \Sigma$.}

In Figure \ref{fig:CN_r}, we plot the $l^2$ norm of the inverse matrix versus $r$ for all four methods used:  $\bs\Sigma^{-1}$ for DLRA, $\bm C^{-1}$ for DO, $\bm Q(\bm p,:)^{-1}$ for TDB-CUR (DEIM) and $\bm Q(\bm p,:)^{\dagger}$ for TDB-CUR (OS-DEIM). As $r$ increases, the matrices $\bs\Sigma^{-1}$ and $\bm C^{-1}$ become ill-conditioned, hence the condition numbers for DLRA and DO become unbounded. On the other hand, the condition numbers for TDB-CUR (DEIM) and TDB-CUR (OS-DEIM) remain nearly constant since the matrix $\bm Q(\bm p,:)$ is well-conditioned. We also observe that the condition number for TDB-CUR (DEIM) can be improved by oversampling as seen in the plot for TDB-CUR (OS-DEIM). 

In the Figure \ref{fig:Error_dt}, we compare  
$\mathcal{E}(T_f)$ versus step size $\Delta t$ for various reduction orders $r = 8, 16, 32$. We observe that TDB-CUR (OS-DEIM) saturates to the optimal low-rank error for each $r$ much quicker than using the unconventional robust integrator \cite{ceruti2021unconventional}. Furthermore, the TDB-CUR method retains the fourth-order accuracy of the Runge-Kutta scheme, whereas the order of accuracy for the unconventional robust integrator is first order, despite using fourth-order Runge-Kutta  for each substep of the algorithm. This confirms the first-order temporal accuracy of the unconventional integrator \cite[Section 3.1]{ceruti2021unconventional}.

\subsection{Stochastic Burgers Equation}
For the second test case, we consider the one-dimensional Burgers equation subject to random initial and boundary conditions as follows: 
\begin{equation*}
\begin{aligned}
&\frac{\partial v}{\partial t}+ \frac{1}{2}\frac{\partial v^2}{\partial x} = \nu \frac{\partial^{2} v}{\partial x^{2}}, && x \in[0,1], t \in[0,5], \\
&v(x, 0 ; \bs\xi) = \sin(2\pi x)\left[0.5\left(e^{\cos(2\pi x)}-1.5\right)+\sigma \sum_{i=1}^{d} \sqrt{\lambda_{x_{i}}} \psi_{i}(x) \xi_{i}\right], &&   x \in[0,1], \xi_{i}\sim \mathcal{N}(\mu,\sigma^2), \\
&v(0, t; \bs\xi) = -\sin (2 \pi t)+\sigma \sum_{i=1}^{d} \lambda_{t_{i}}\varphi_{i}(t) \xi_{i}, && x=0, \xi_{i}\sim \mathcal{N}(\mu,\sigma^2),
\end{aligned}
\end{equation*}
where $\nu=2.5\times 10^{-3}$. The stochastic boundary at $x=0$ is specified above and the boundary at $x=1$ is $v(x=1,t;\bs\xi)=0$. We use weak treatment of the boundary conditions for both the FOM and TDB \cite{patil2023reduced}. The random space is taken to be $d = 17$ dimensional and $\xi_{i}$'s are sampled from a normal distribution with mean $\mu=0$, standard deviation $\sigma=0.001$, and $s= 256$. In the stochastic boundary specification, we take $\varphi_{i}(t)=\sin(i\pi t)$ and  $\lambda_{t_{i}}=i^{-2}$. In the stochastic initial condition, $\lambda_{x_{i}}$ and $\psi_{i}(x)$ are the eigenvalues and eigenvectors of the spatial squared-exponential kernel, respectively. The fourth-order explicit Runge-Kutta method is used for time integration \cb{of the FOM and TDB-CUR} with $\Delta t=2.5 \times 10^{-4}$. For discretization of the spatial domain, we use a second-order finite difference scheme on a uniform grid with $n=401$. \cb{This leads to the following MDE of the form $\D\bm V/\D t = \mathcal{F}(t,\bm V)$:
\begin{equation*}
\frac{\D\bm V(t)}{\D t} = -\frac{1}{2}\bm D_1 (\bm V(t)\odot \bm V(t)) + \nu\bm D_2 \bm V(t) + \bm B(t), \quad \bm V(0)=\bm V_0,
\end{equation*}
where $\bm D_1$ and $\bm D_2$ are $n\times n$ sparse matrices defining the first and second spatial derivatives of the discretized system. The first and last row of $\bm D_1$ and $\bm D_2$ are equal to zero. The $n\times s$ matrix $\bm B(t)$ enforces the stochastic boundary at $x=0$ by setting each element in its first row equal to $\D v(0,t;\bs\xi)/\D t$, for $s$ independent samples of the random variables $\bs \xi$. All other entries of $\bm B(t)$ are equal to zero. The columns of $\bm V_0$ are the initial conditions for $s$ samples of the random variables.
}

\begin{figure}[!t]
    \centering
     \subfigure[Error]{
         \centering
         \input{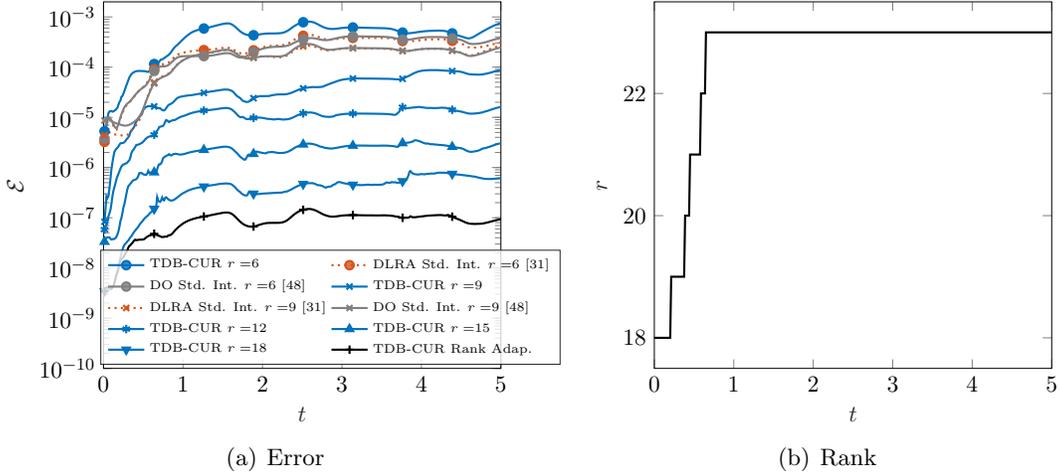}
         \label{fig:Burger-error}
     }
     \subfigure[Rank]{
%
%
\definecolor{mycolor1}{rgb}{0.00000,0.44700,0.74100}%
\definecolor{mycolor2}{rgb}{0.85000,0.32500,0.09800}%
\begin{tikzpicture}[scale=0.8]

\begin{axis}[%
width=2.6in,
height=2.4in,
at={(1.648in,0.642in)},
scale only axis,
xmin=0,
xmax=5,
xlabel style={font=\color{white!15!black}},
xlabel={$t$},
ylabel style={font=\color{white!15!black}},
ylabel={$r$},
axis background/.style={fill=white},
legend style={at={(0.97,0.03)}, anchor=south east, legend cell align=left, align=left, draw=white!15!black,font=\scriptsize}
]
\addplot [color=black, line width=1.0pt]
  table[row sep=crcr]{%
0	18\\
0.0125	18\\
0.025	18\\
0.0375	18\\
0.05	18\\
0.0625	18\\
0.075	18\\
0.0875	18\\
0.1	18\\
0.1125	18\\
0.125	18\\
0.1375	18\\
0.15	18\\
0.1625	18\\
0.175	18\\
0.1875	18\\
0.2	18\\
0.2125	19\\
0.225	19\\
0.2375	19\\
0.25	19\\
0.2625	19\\
0.275	19\\
0.2875	19\\
0.3	19\\
0.3125	19\\
0.325	19\\
0.3375	19\\
0.35	19\\
0.3625	19\\
0.375	19\\
0.3875	20\\
0.4	20\\
0.4125	20\\
0.425	20\\
0.4375	20\\
0.45	21\\
0.4625	21\\
0.475	21\\
0.4875	21\\
0.5	21\\
0.5125	21\\
0.525	21\\
0.5375	21\\
0.55	21\\
0.5625	21\\
0.575	21\\
0.5875	22\\
0.6	22\\
0.6125	22\\
0.625	22\\
0.6375	22\\
0.65	23\\
0.6625	23\\
0.675	23\\
0.6875	23\\
0.7	23\\
0.7125	23\\
0.725	23\\
0.7375	23\\
0.75	23\\
0.7625	23\\
0.775	23\\
0.7875	23\\
0.8	23\\
0.8125	23\\
0.825	23\\
0.8375	23\\
0.85	23\\
0.8625	23\\
0.875	23\\
0.8875	23\\
0.9	23\\
0.9125	23\\
0.925	23\\
0.9375	23\\
0.95	23\\
0.9625	23\\
0.975	23\\
0.9875	23\\
1	23\\
1.0125	23\\
1.025	23\\
1.0375	23\\
1.05	23\\
1.0625	23\\
1.075	23\\
1.0875	23\\
1.1	23\\
1.1125	23\\
1.125	23\\
1.1375	23\\
1.15	23\\
1.1625	23\\
1.175	23\\
1.1875	23\\
1.2	23\\
1.2125	23\\
1.225	23\\
1.2375	23\\
1.25	23\\
1.2625	23\\
1.275	23\\
1.2875	23\\
1.3	23\\
1.3125	23\\
1.325	23\\
1.3375	23\\
1.35	23\\
1.3625	23\\
1.375	23\\
1.3875	23\\
1.4	23\\
1.4125	23\\
1.425	23\\
1.4375	23\\
1.45	23\\
1.4625	23\\
1.475	23\\
1.4875	23\\
1.5	23\\
1.5125	23\\
1.525	23\\
1.5375	23\\
1.55	23\\
1.5625	23\\
1.575	23\\
1.5875	23\\
1.6	23\\
1.6125	23\\
1.625	23\\
1.6375	23\\
1.65	23\\
1.6625	23\\
1.675	23\\
1.6875	23\\
1.7	23\\
1.7125	23\\
1.725	23\\
1.7375	23\\
1.75	23\\
1.7625	23\\
1.775	23\\
1.7875	23\\
1.8	23\\
1.8125	23\\
1.825	23\\
1.8375	23\\
1.85	23\\
1.8625	23\\
1.875	23\\
1.8875	23\\
1.9	23\\
1.9125	23\\
1.925	23\\
1.9375	23\\
1.95	23\\
1.9625	23\\
1.975	23\\
1.9875	23\\
2	23\\
2.0125	23\\
2.025	23\\
2.0375	23\\
2.05	23\\
2.0625	23\\
2.075	23\\
2.0875	23\\
2.1	23\\
2.1125	23\\
2.125	23\\
2.1375	23\\
2.15	23\\
2.1625	23\\
2.175	23\\
2.1875	23\\
2.2	23\\
2.2125	23\\
2.225	23\\
2.2375	23\\
2.25	23\\
2.2625	23\\
2.275	23\\
2.2875	23\\
2.3	23\\
2.3125	23\\
2.325	23\\
2.3375	23\\
2.35	23\\
2.3625	23\\
2.375	23\\
2.3875	23\\
2.4	23\\
2.4125	23\\
2.425	23\\
2.4375	23\\
2.45	23\\
2.4625	23\\
2.475	23\\
2.4875	23\\
2.5	23\\
2.5125	23\\
2.525	23\\
2.5375	23\\
2.55	23\\
2.5625	23\\
2.575	23\\
2.5875	23\\
2.6	23\\
2.6125	23\\
2.625	23\\
2.6375	23\\
2.65	23\\
2.6625	23\\
2.675	23\\
2.6875	23\\
2.7	23\\
2.7125	23\\
2.725	23\\
2.7375	23\\
2.75	23\\
2.7625	23\\
2.775	23\\
2.7875	23\\
2.8	23\\
2.8125	23\\
2.825	23\\
2.8375	23\\
2.85	23\\
2.8625	23\\
2.875	23\\
2.8875	23\\
2.9	23\\
2.9125	23\\
2.925	23\\
2.9375	23\\
2.95	23\\
2.9625	23\\
2.975	23\\
2.9875	23\\
3	23\\
3.0125	23\\
3.025	23\\
3.0375	23\\
3.05	23\\
3.0625	23\\
3.075	23\\
3.0875	23\\
3.1	23\\
3.1125	23\\
3.125	23\\
3.1375	23\\
3.15	23\\
3.1625	23\\
3.175	23\\
3.1875	23\\
3.2	23\\
3.2125	23\\
3.225	23\\
3.2375	23\\
3.25	23\\
3.2625	23\\
3.275	23\\
3.2875	23\\
3.3	23\\
3.3125	23\\
3.325	23\\
3.3375	23\\
3.35	23\\
3.3625	23\\
3.375	23\\
3.3875	23\\
3.4	23\\
3.4125	23\\
3.425	23\\
3.4375	23\\
3.45	23\\
3.4625	23\\
3.475	23\\
3.4875	23\\
3.5	23\\
3.5125	23\\
3.525	23\\
3.5375	23\\
3.55	23\\
3.5625	23\\
3.575	23\\
3.5875	23\\
3.6	23\\
3.6125	23\\
3.625	23\\
3.6375	23\\
3.65	23\\
3.6625	23\\
3.675	23\\
3.6875	23\\
3.7	23\\
3.7125	23\\
3.725	23\\
3.7375	23\\
3.75	23\\
3.7625	23\\
3.775	23\\
3.7875	23\\
3.8	23\\
3.8125	23\\
3.825	23\\
3.8375	23\\
3.85	23\\
3.8625	23\\
3.875	23\\
3.8875	23\\
3.9	23\\
3.9125	23\\
3.925	23\\
3.9375	23\\
3.95	23\\
3.9625	23\\
3.975	23\\
3.9875	23\\
4	23\\
4.0125	23\\
4.025	23\\
4.0375	23\\
4.05	23\\
4.0625	23\\
4.075	23\\
4.0875	23\\
4.1	23\\
4.1125	23\\
4.125	23\\
4.1375	23\\
4.15	23\\
4.1625	23\\
4.175	23\\
4.1875	23\\
4.2	23\\
4.2125	23\\
4.225	23\\
4.2375	23\\
4.25	23\\
4.2625	23\\
4.275	23\\
4.2875	23\\
4.3	23\\
4.3125	23\\
4.325	23\\
4.3375	23\\
4.35	23\\
4.3625	23\\
4.375	23\\
4.3875	23\\
4.4	23\\
4.4125	23\\
4.425	23\\
4.4375	23\\
4.45	23\\
4.4625	23\\
4.475	23\\
4.4875	23\\
4.5	23\\
4.5125	23\\
4.525	23\\
4.5375	23\\
4.55	23\\
4.5625	23\\
4.575	23\\
4.5875	23\\
4.6	23\\
4.6125	23\\
4.625	23\\
4.6375	23\\
4.65	23\\
4.6625	23\\
4.675	23\\
4.6875	23\\
4.7	23\\
4.7125	23\\
4.725	23\\
4.7375	23\\
4.75	23\\
4.7625	23\\
4.775	23\\
4.7875	23\\
4.8	23\\
4.8125	23\\
4.825	23\\
4.8375	23\\
4.85	23\\
4.8625	23\\
4.875	23\\
4.8875	23\\
4.9	23\\
4.9125	23\\
4.925	23\\
4.9375	23\\
4.95	23\\
4.9625	23\\
4.975	23\\
4.9875	23\\
5	23\\
};

\end{axis}

\end{tikzpicture}%
        \label{fig:Burger-rank}
    }
    \caption{Stochastic Burgers equation with constant diffusion: (a) Relative error evolution for fixed and adaptive rank. (b) Rank evolution using upper threshold $\epsilon_u=10^{-8}$ for mode addition.}
    \label{fig:Burgers-error-rank}
\end{figure}

We first solve the system using TDB-CUR with fixed rank and compare the results against the DLRA and DO by solving Eqs. \ref{eq:DBO_evol_S}-\ref{eq:DBO_evol_Y} \cb{and Eqs. \ref{eqn:u-do}-\ref{eqn:y-do}, respectively. The fourth-order explicit Runge-Kutta method (a standard integrator) is used to solve both the DLRA and DO equations}. For TDB-CUR, the rows are oversampled with $m=5$. No sparse sampling strategy is used for DLRA or DO, and Eqs. \ref{eq:DBO_evol_S}-\ref{eq:DBO_evol_Y} and Eqs. \ref{eqn:u-do}-\ref{eqn:y-do} are solved as is. In Figure \ref{fig:Burger-error}, we compare the error of TDB-CUR, DLRA, and DO for different values of $r$. For $r=6$, TDB-CUR has larger error compared to \cb{both DLRA and DO}. This result is expected since TDB-CUR has an additional source of error from the sparse sampling procedure. However, as the rank is increased to $r=9$, the conditioning of the DLRA \cb{ and DO with standard integrator} starts to deteriorate and the error of TDB-CUR is actually lower than DLRA and DO. In fact, for $r>9$, DLRA \cb{and DO with standard integrators} are unstable and cannot be integrated beyond the first time step. On the other hand, TDB-CUR remains stable, and the error decays as the rank is increased to a maximum value of $r=18$. While it is reasonable to expect that the error can be reduced further by increasing the rank to values of $r>18$, it is important to note that the rank of the initial condition is exactly $r=18$. Therefore, in order to increase the rank of the system beyond $r=18$ in a principled manner, we employ the rank adaptive strategy from Section \ref{sec:adaptive-rank}. To this end, we initialize the system with rank $r_0=18$, and use an upper threshold of $\epsilon_u=10^{-8}$ for mode addition. As observed in Figure \ref{fig:Burger-rank}, the rank is increased in time to a maximum of 23, leading to a further reduction in the error. 

\begin{figure}[!t]
    \centering
    \includegraphics[width=.75\textwidth]{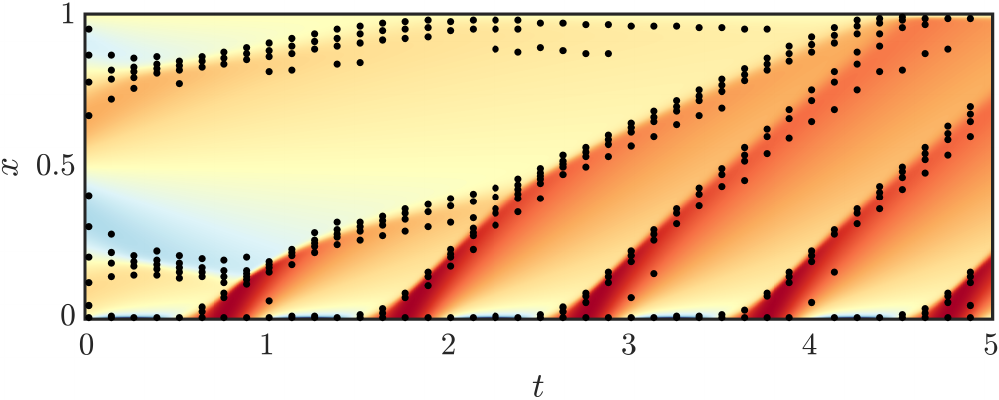}
    \caption{Stochastic Burgers equation with constant diffusion: mean solution with the first 10 QDEIM points (black dots).}
    \label{fig:Burger-mean-sln}
\end{figure}

\begin{figure}
    \centering
    \input{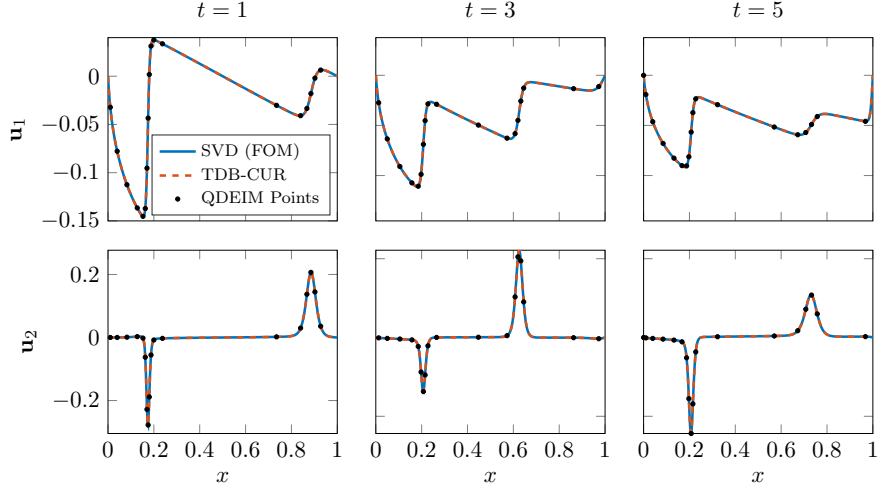}
    \caption{Stochastic Burgers equation with constant diffusion: Evolution of first two spatial modes, $\bm u_1$ and $\bm u_2$, and QDEIM points. Excellent agreement between the FOM and TDB-CUR is observed.}
    \label{fig:Burgers-u-modes}
\end{figure}

The mean solution is shown in Figure \ref{fig:Burger-mean-sln} along with the first 10 QDEIM sampling points. We observe that the sampling points are concentrated near the stochastic boundary at $x=0$ and also at points in the domain where shocks develop. Figure \ref{fig:Burgers-u-modes} shows the evolution of the first two spatial modes, $\bm u_1$ (top) and $\bm u_2$ (bottom), where we observe excellent agreement between the FOM and TDB-CUR. It is important to note that these modes are energetically ranked according to the first and second singular values shown in Figure \ref{fig:Burger-sval}. Therefore, we observe that $\bm u_1$ captures the large scale energy containing structure, while $\bm u_2$ captures the small scale structure that is highly localized in space. 

In Figure \ref{fig:Burger-sval}, we show that TDB-CUR accurately captures the leading singular values of the FOM solution, despite the large gap between the first and last resolved singular values. Finally, Figure \ref{fig:Burger-scaling-n=s} compares the CPU time of the FOM, DLRA \cb{with the unconventional integrator}, and TDB-CUR as the number of rows and columns of the matrix are increased simultaneously. We take $n=s$ and observe that the FOM scales quadratically ($\mathcal{O}(ns)$) while TDB-CUR and \cb{DLRA with the unconventional integrator} scale linearly ($\mathcal{O}(n+s)$). As the matrix size is increased, the disparity in CPU time becomes even more apparent, making the case for solving massive MDEs using \cb{low-rank approximation}. \cb{Despite this result, it is important to note that linear scaling for the unconventional integrator is only possible since the nonlinear term in the Burgers equation is limited to quadratic. For higher-order polynomial and general nonlinearities, the unconventional integrator will scale with $\mathcal{O}(ns)$, and exceed the cost of solving the FOM. Nevertheless, given the factored rank-$r$ approximation, $\Vhat=\bU\bs\Sigma\bY^T$, the quadratic term can be computed efficiently, resulting in a factorization that has a maximum rank of $(r^2+r)/2$.}


\cb
{To demonstrate the true power of the TDB-CUR method, we modify the right hand side of the MDE by making the diffusion term nonlinear, $\nu(1+\tanh(\bm V)\bs \xi')\odot(\bm D_2\bm V)$. To clarify, $\tanh$ is evaluated element-wise on its argument, and $\bs \xi'$ is an $s\times s$ diagonal matrix with elements drawn from $\mathcal N(\mu=0,\sigma^2=0.01^2)$. Furthermore, we verify that $\bs \xi'$ does not result in a negative diffusion. As a result of this simple modification, the cost of directly computing $\mathcal{F}(t,\Vhat)$ will scale with $\mathcal O(ns)$, even for $\Vhat$ of low-rank. Therefore, efficient computation of DLRA with a standard integrator \cite{KL07}, unconventional integrator \cite{CL21}, or projection method \cite{kieri2019projection} is not possible. To highlight this, we compare the error versus cost (time to solution) for TDB-CUR, DLRA with the unconventional integrator, and DLRA using projection methods. We use the projected fourth-order Runge-Kutta method (PRK4) presented in \cite{kieri2019projection} along with fourth-order Runge Kutta for both TDB-CUR and the substeps of the unconventional integrator. For the PRK4 method, the unfactored $n\times s$ matrix, $\bm F_i$, is computed at the $i^{\mathrm{th}}$ stage of the integration scheme. The $n\times s$ matrix is then projected to the tangent space of the rank-$r$ manifold at the $i^{\mathrm{th}}$ stage as, $\mathcal{P}_{\mathcal{T}_{\Vhat_i}}(\bm F_i)$, resulting in a matrix whose rank is at most $2r$. Note that using the subscript to denote the stages results in $\Vhat_1=\Vhat^{k-1}.$ To limit rank growth during the internal steps of the RK4 method, the economy size SVD is applied after each sub-step to obtain the rank-$r$ $\Vhat_i=\bU_i\bs\Sigma_i\bY_i^T$. This is only necessary for $i>1$, since $\Vhat_1=\Vhat^{k-1}$ (see above).  The orthonormal column and row bases, $\bm U_i$ and $\bm Y_i$, are then used for the tangent space projections at each stage. One final economy size SVD is applied so that the updated low-rank matrix, $\Vhat^{k}$, remains on the rank-$r$ manifold. Although the nonlinear diffusion requires forming $\bm F_i$ of size $n\times s$, our implementation does not require computing the SVD of matrices larger than $n\times 8r$. To our knowledge, this represents an efficient implementation of PRK4 when forming the full $\bm F_i$ cannot be avoided.

Figure \ref{fig:Burgers-nonlin-diff-error-vs-cost} shows the error versus cost for $r=6,9,12,15,18$.  For TDB-CUR, we observe a rapid decrease in error for a modest increase in cost. Similar behavior is observed for PRK4, however, both the error and cost exceed those of TDB-CUR. Finally, the unconventional integrator exhibits a larger error than both TDB-CUR and projection for a given $\Delta t$. Due to the first-order accuracy of the unconventional integrator, the error does not monotonically decrease as the cost (rank) is increased. To verify this, we decrease $\Delta t$ by an order of magnitude and rerun. Despite the error dropping by an order of magnitude, the same non-monotonic behavior in the error is observed. This confirms the error in the unconventional integrator is still dominated by the temporal error and not the low-rank approximation error. Finally, we plot the error vs time for TDB-CUR and PRK4 in Figure \ref{fig:Burgers-nonlin-diff-error-vs-time}. As the rank is increased, we observe a corresponding decrease in the error for both methods. However, TDB-CUR ultimately  achieves lower error than PRK4 as the rank is increased.

\begin{figure}
\centering
\subfigure[Singular Values]{
         \centering
         \input{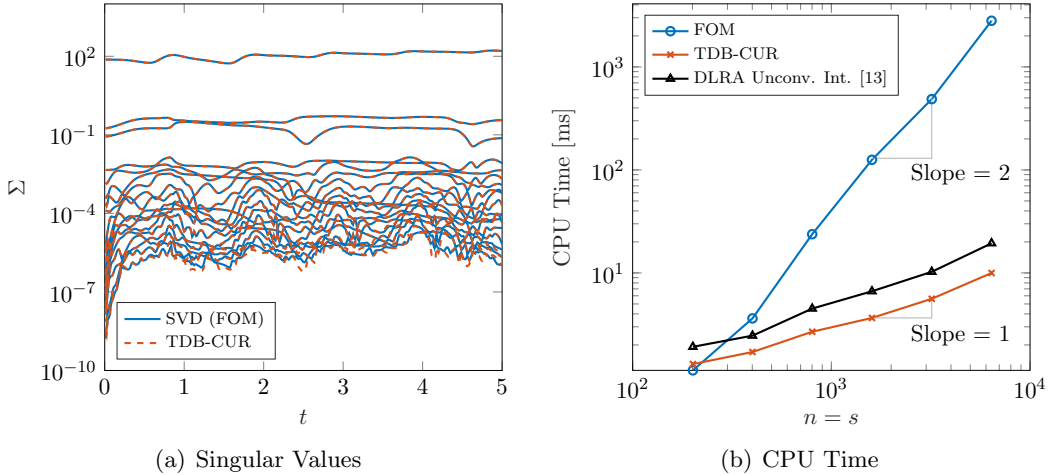}
         \label{fig:Burger-sval}
     }
\subfigure[CPU Time]{
%
%
\definecolor{mycolor1}{rgb}{0.00000,0.44700,0.74100}%
\definecolor{mycolor2}{rgb}{0.85000,0.32500,0.09800}%
\definecolor{mycolor3}{rgb}{0.92900,0.69400,0.12500}%
\definecolor{mygrey}{rgb}{.7 .7 .7}%
\begin{tikzpicture}[scale=0.8]

\begin{axis}[%
width=2.6in,
height=2.4in,
at={(1.648in,0.642in)},
scale only axis,
xmode=log,
xmin=100,
xmax=10000,
xminorticks=true,
xlabel style={font=\color{white!15!black}},
xlabel={$n=s$},
ymode=log,
ymin=1.13279892,
ymax=4097.2801,
yminorticks=true,
ylabel style={font=\color{white!15!black}},
ylabel={CPU Time [ms]},
axis background/.style={fill=white},
legend style={at={(0.03,0.97)}, anchor=north west, legend cell align=left, align=left, draw=white!15!black,,font=\scriptsize}
]
\addplot [color=mycolor1, line width=1.0pt, mark size=2.0pt, mark=o, mark options={solid, mycolor1}]
  table[row sep=crcr]{%
201	1.13279892\\
401	3.62076176\\
801	23.75729965\\
1601	125.39726892\\
3201	487.03110159\\
6401	2805.93614055\\
};
\addlegendentry{FOM}


\addplot [color=mycolor2, line width=1.0pt, mark size=2.0pt, mark=x, mark options={solid, mycolor2}]
  table[row sep=crcr]{%
201	1.30391866\\
401	1.70772596\\
801	2.6859933\\
1601	3.65645136\\
3201	5.61082046\\
6401	9.98523239\\
};
\addlegendentry{TDB-CUR}



\addplot [color=black, line width=1.0pt, mark size=2pt, mark=triangle, mark options={black}]
  table[row sep=crcr]{%
201	1.918331818\\
401	2.460074145\\
801	4.505679654\\
1601	6.638554972\\
3201	10.257292828\\
6401	19.393312501\\
};
\addlegendentry{DLRA Unconv. Int. \cite{ceruti2021unconventional}}

\addplot [color=mygrey, forget plot]
  table[row sep=crcr]{%
1601	129.861230485\\
3201	129.861230485\\
3201	521.79370747\\
};
\node[right, align=left, inner sep=0]
at (axis cs:2500,90) {Slope = 2};

\addplot [color=mygrey, forget plot]
  table[row sep=crcr]{%
1601	3.65645136\\
3201    3.65645136\\
3201	5.61082046\\
};
\node[right, align=left, inner sep=0]
at (axis cs:2500,2.5) {Slope = 1};

\end{axis}
\end{tikzpicture}%
    \label{fig:Burger-scaling-n=s}
}
\caption{Stochastic Burgers equation with constant diffusion: (a) First $18$ singular values of FOM vs TDB-CUR. (b) CPU time for scaling $n=s$ for FOM, TDB-CUR, \cb{and the unconventional integrator} with fixed $r=6$.}
\label{fig:Burgers-sval-scaling}
\end{figure}

\begin{figure}
    \centering
    \subfigure[]{
%
%
\definecolor{mycolor1}{rgb}{0.00000,0.44700,0.74100}%
\definecolor{mycolor2}{rgb}{0.85000,0.32500,0.09800}%
\begin{tikzpicture}[scale=0.8]

\begin{axis}[%
width=2.6in,
height=2.4in,
at={(1.648in,0.642in)},
scale only axis,
xmode=log,
xmin=0,
xmax=6000,
xlabel style={font=\color{white!15!black}},
xlabel={Cost [s]},
ymode=log,
ymin=1e-08,
ymax=0.1,
yminorticks=true,
ylabel style={font=\color{white!15!black}},
ylabel={$\mathcal{E}(t=5)$},
axis background/.style={fill=white},
legend style={at={(1.0,0.001)}, anchor=south east, legend cell align=left, align=left, draw=white!15!black, font=\scriptsize}
]
\addplot [color=mycolor1, line width=1.0pt, mark=o, mark options={solid, mycolor1}]
  table[row sep=crcr]{%
32.292299176	0.00134561195173738\\
46.554076395	0.000191610558938797\\
52.136482276	2.17618780030836e-05\\
65.186446751	6.44947425678869e-06\\
74.688002065	1.43314489884038e-06\\
};
\addlegendentry{TDB-CUR ($\Delta t$)}

\addplot [color=mycolor2, line width=1.0pt, mark=triangle, mark options={solid, mycolor2}]
  table[row sep=crcr]{%
156.968185222	0.00118012726834675\\
188.889361241	0.000290332856922459\\
219.261269075	0.000196491344035316\\
280.739958994	7.92265511587797e-05\\
296.901759136	3.19585289941803e-05\\
};
\addlegendentry{PRK4 ($\Delta t$) \cite{kieri2019projection}}

\addplot [color=black, line width=1.0pt, mark=x, mark options={solid, black}]
  table[row sep=crcr]{%
323.938853955	0.0144670365482364\\
323.536588414	0.00281983073784363\\
325.78167257	0.00927711027718688\\
329.611553543	0.0124355520550417\\
334.361745292	0.00578480487387787\\
};
\addlegendentry{DLRA Unconv. Int. ($\Delta t$) \cite{ceruti2021unconventional}}

\addplot [color=gray, line width=1.0pt, mark=square, mark options={solid, gray}]
  table[row sep=crcr]{%
3236.021958347	0.00188219312228112\\
3255.015563846	0.000565733733721273\\
3241.371883646	0.00104886292012336\\
3286.978552544	0.00145382374041178\\
3307.894076045	0.000620723018193936\\
};
\addlegendentry{DLRA Unconv. Int. ($\Delta t/10)$ \cite{ceruti2021unconventional}}

\end{axis}

\end{tikzpicture}%
        \label{fig:Burgers-nonlin-diff-error-vs-cost}
    }
    \subfigure[]{
        \input{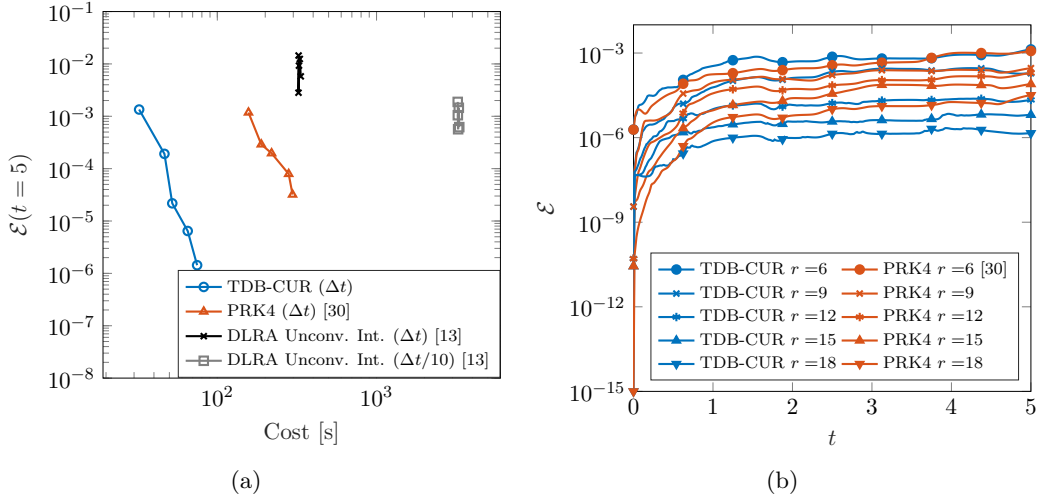}
        \label{fig:Burgers-nonlin-diff-error-vs-time}
    }
    \caption{Stochastic Burgers equation with nonlinear diffusion: (a) Error versus time. (b) Error at $t=5$ versus total cost in seconds. The cost and error data points are obtained by considering $r=6,9,12,15,18$.}
    \label{fig:Burgers-error-cost-entire-fig}
\end{figure}

Although it is not entirely obvious why PRK4 has a larger error as the rank is increased, we propose one possible explanation based on the curvature of the manifold. To this end, it is well known  that the curvature of the manifold is inversely proportional to the smallest singular value in the low-rank solution \cite{LO14, rodgers2022adaptive}. Therefore, as the rank is increased in the above example,  the curvature of the manifold at the low-rank solution increases rapidly. As the curvature increases, the tangent space will no longer provide a good approximation for small deviations (e.g. $\mathcal{O}(\Delta t)$) from the low-rank solution at that point. Since PRK4  takes noninfinitesimal time steps off the rank-$r$ manifold, the subsequent tangent space projections may induce errors that can be large for points on the manifold with high-order curvature. Therefore, one possible explanation for the above result is that the tangent space projection incurs a larger error since its accuracy relies heavily on the curvature of the manifold at that point. On the other hand, the TDB-CUR method does not use tangent space projections, and does not suffer from the high-order curvature of the manifold. For more details, we refer the reader to \cite{CL23} for an excellent discussion on the error induced by the tangent space projection.

}

While Figures \ref{fig:Burgers-error-rank}-\ref{fig:Burgers-error-cost-entire-fig} demonstrate the accuracy, efficiency, rank-adaptivity, and favorable numerical performance of the TDB-CUR method, they do not convey the minimally intrusive nature of its implementation. To give  a better perspective on the implementation efforts, the MATLAB code for solving the stochastic Burgers equation using the TDB-CUR method is provided in Appendix \ref{sec:matlab-code} (Listings \ref{list:matlab} and \ref{list:gpode}). While the code contains lines specific to the TDB-CUR method, after reviewing the entire code, it will become apparent that many of the included lines are already required for solving the FOM Burgers equation. Furthermore, there is no term-by-term implementation required to preserve efficiency and the FOM implementation of the Burgers equation (\texttt{function f}) is used to compute the sparse row and column samples. Therefore, given an existing FOM implementation, the code required to implement the TDB-CUR method is  minimal. The code blocks required for implementing the method are labeled with \texttt{\%\% TDB-CUR} in the attached code.

\subsection{Stochastic Advection-Diffusion-Reaction Equation}

\begin{figure}
    \centering
    \includegraphics[width=.8\textwidth]{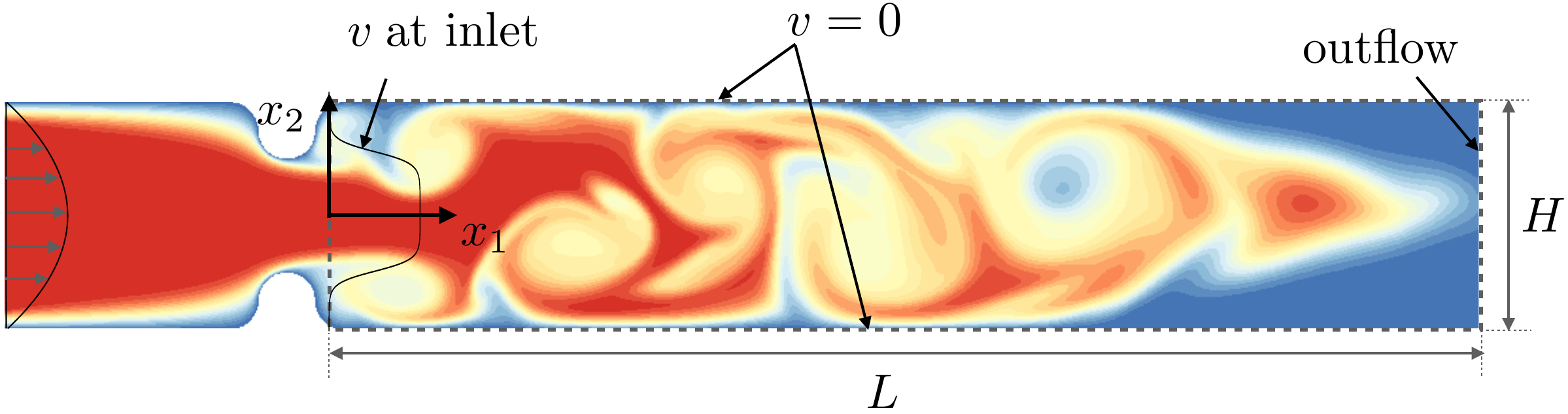}
    \caption{Schematic of the flow visualized with a passive scalar.}
    \label{fig:rxn_setup}
\end{figure}

In this section, we aim to solve the 2D advection-diffusion-reaction (ADR) equation subject to random diffusion coefficient ($\alpha$), with deterministic initial condition:

\begin{equation*}
\begin{aligned}
&\frac{\partial v}{\partial t} + (u \cdot \nabla) v = \nabla \cdot(\nabla v\alpha) + \frac{v^{2}}{10+v}, \quad x_1\in[0,10], \, x_2\in[0,2], \, t \in[0,5], \\
&v(x_1, x_2, 0) = \frac{1}{2}\left(\tanh\left(\frac{x_2+0.5}{0.1}\right) - \tanh\left(\frac{x_2-0.5}{0.1}\right)\right),
\end{aligned}
\end{equation*}
where $v(x_1,x_2,t)$ is the species concentration and $u(x_1, x_2, t)$ is the velocity vector. It is worth noting that the nonlinearity of the equation is non-polynomial, implying that the computational expense of DO or DLRA is comparable to that of the FOM.  The schematic of the problem is shown in Figure \ref{fig:rxn_setup}. The velocity field is obtained by solving the incompressible Navier-Stokes equations  and is independent of the species transport equation. The conditions are identical to those used in previous studies \cite{RNB21,DCB22}.  In particular, we solved the velocity field in the entire domain using the spectral/hp element method on an unstructured mesh with 4008 quadrilateral elements and polynomial order 5.  For more details on the spectral element method see for example \cite{KS05,Babaee:2013aa}. At the inlet, a parabolic velocity is prescribed, with an average velocity of $\overline{u}$. The outflow condition is imposed at the right boundary and the no-slip boundary condition is imposed at the remaining boundaries. The Reynolds number with reference length $H/2$, and kinematic viscosity $\nu$, is given by $Re=\overline{u}H/2\nu=1000$,    

We solved the ADR and TDB-CUR equations using a collocated spectral element method within the rectangular domain indicated by dashed lines in Figure \ref{fig:rxn_setup}. In particular, we use a uniform quadrilateral mesh with 50 elements in the $x_1$ direction and 15 elements in the $x_2$ direction, and a spectral polynomial of order 5 in each direction  within the rectangular domain.  This results in  $n=19076$ degrees of freedom in the spatial domain. We  interpolated the velocity field from the unstructured mesh onto the structured mesh. 
The fourth-order explicit Runge-Kutta method is utilized for time integration with $\Delta t=5 \times 10^{-4}$ for advancing the ADR and TDB-CUR equations. 

Unlike the previous example, we use a deterministic initial condition. Therefore, the  rank at $t=0$ is exactly equal to one, i.e. $\mathrm{rank}(\bm V(0))=1$. While the low-rank approximation with $r=1$ will be exact in its initial condition, the rank of the system will quickly increase due to the nonlinearity. Therefore, to maintain an acceptable level of error, the rank of the approximation must increase in time. While it is possible to initialize TDB-CUR with $r>1$, we opt to use the rank-adaptive strategy from Algorithm \ref{alg:S-TDB-ROM},  starting with the initial rank of $r_0=1$. \cb{Similarly, DLRA using the unconventional integrator can also be initialized with $r>1$, however, several rank-adaptive integrators have been proposed \cite{yang2020time, lubich2021rankadaptive, dektor2021rank}.} On the other hand, initializing DLRA or DO \cb{for a standard integrator} with $r>1$ is not possible, as $\bs \Sigma$ and $\bm C$ will be singular. Therefore, this problem setup emphasizes the need for rank adaptivity for TDB-based low-rank matrix approximation. 

\begin{figure}[!t]
    \centering
    \includegraphics[width=0.95\textwidth]{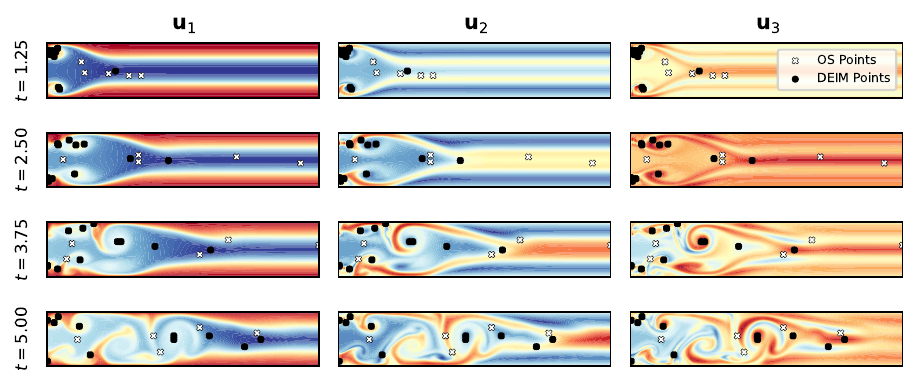}
    \caption{Stochastic Advection-Diffusion-Reaction equation: First three spatial modes of TDB-CUR at different time-steps and
the selected points of DEIM with $m=$5 oversampled (OS) points. ($s=50$)}
    \label{fig:ADR_Solution}
\end{figure}

\begin{figure}[!t]
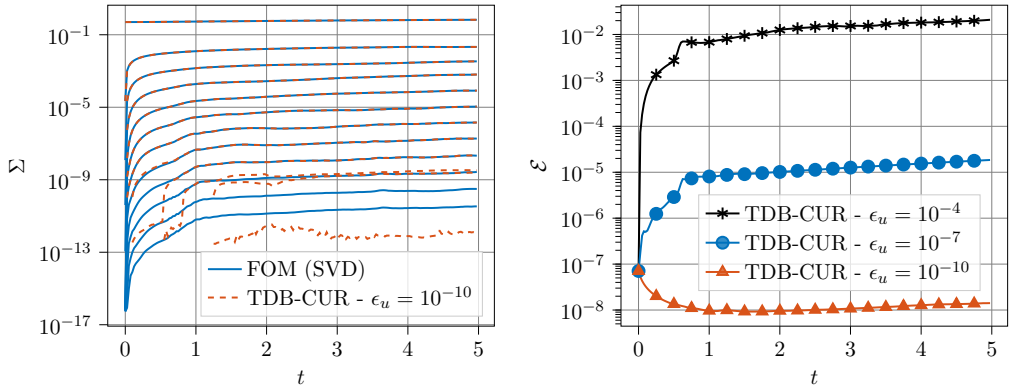

     \centering
    \begin{subfigure}
         \centering
         \input{Figures/ADR_Sigma}
         \label{fig:ADR_Sigma}
     \end{subfigure}
     \begin{subfigure}
         \centering
         \input{Figures/ADR_Error}
         \label{fig:ADR_Error}
     \end{subfigure}
        \caption{Stochastic Advection-Diffusion-Reaction equation: Left: Comparison of first 10 singular values of FOM vs TDB-CUR with 5 oversampled points ($s =50$). Right: Error ($\mathcal{E}$) of TDB-CUR versus time. The result is presented for different values of the upper error bound ($\epsilon_u$) for mode addition.}
        \label{fig:ADR_Error_Sigma}
\end{figure}

For the first case, we consider a random diffusion coefficient according to  $\alpha=\frac{1}{\xi}$, where $\xi$ is a Gaussian random variable with a mean of $\mu=100$ and standard deviation of $\sigma= 25$. Since validating the performance of TDB-CUR requires solving the FOM, we do not consider a large number of samples for the first case. We draw $s=50$ samples of the diffusion coefficient, which allows us to compute the error and compare the singular values with the FOM in a reasonable amount of time. Figure \ref{fig:ADR_Solution} shows the evolution of the first three spatial modes, along with the sparse sampling points. As the simulation evolves in time, the points also evolve as the flow is advected from left to right. In Figure \ref{fig:ADR_Error_Sigma}, the instantaneous singular values from TDB-CUR and the $r$ largest singular values of the FOM solution (SVD singular values) are shown on the left. The discrepancy between the trailing singular values of the TDB-CUR and FOM stems from the effect of unresolved modes in the time integration of the low-rank approximation. However, as observed in the error on the right, accurately resolving the leading singular values results in very small errors. Additionally, the error is  controlled by lowering the error threshold for rank addition, leading to improved accuracy of the TDB-CUR approximation.

\begin{figure}[!t]
     \centering
     \begin{subfigure}
         \centering
         \input{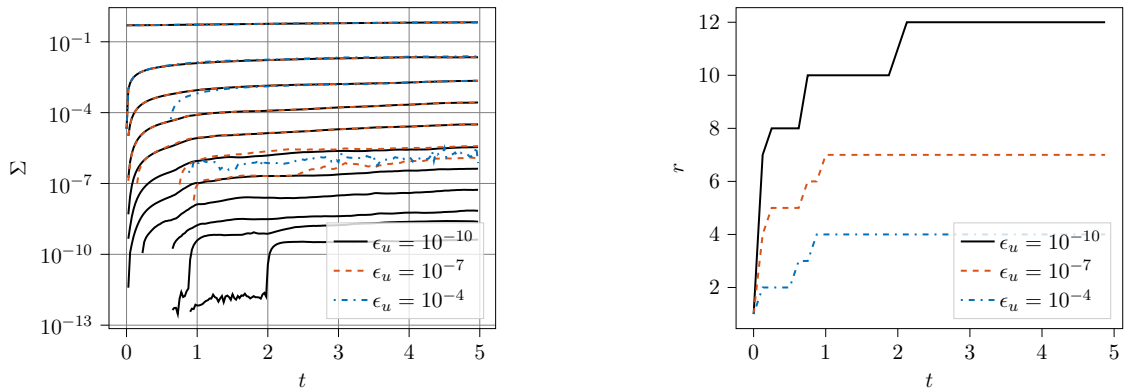}
         \label{fig:ADR_Error}
     \end{subfigure}
     \hfill
     \begin{subfigure}
         \centering
\begin{tikzpicture}[scale=0.75]

  \definecolor{darkgray176}{RGB}{176,176,176}
  \definecolor{lightgray204}{RGB}{204,204,204}
  \definecolor{myblue}{rgb}{0.00000,0.44700,0.74100}%
  \definecolor{myred}{rgb}{0.85000,0.32500,0.09800}%
  \definecolor{myblack}{rgb}{0,0,0}%
  \definecolor{mygrey}{rgb}{.7 .7 .7}

\begin{axis}[
legend cell align={left},
legend style={
  fill opacity=0.8,
  draw opacity=1,
  text opacity=1,
  at={(0.97,0.03)},
  anchor=south east,
  draw=lightgray204
},
tick align=outside,
tick pos=left,
x grid style={darkgray176},
xlabel={\(\displaystyle t\)},
xmin=-0.24325, xmax=5.11925,
xtick style={color=myblack},
y grid style={darkgray176},
ylabel={\(\displaystyle r\)},
ymin=0.45, ymax=12.55,
ytick style={color=myblack}
]
\addplot [line width=1.pt,  myblack]
table {%
0.0005 1
0.1255 7
0.2505 8
0.3755 8
0.5005 8
0.6255 8
0.7505 10
0.8755 10
1.0005 10
1.1255 10
1.2505 10
1.3755 10
1.5005 10
1.6255 10
1.7505 10
1.8755 10
2.0005 11
2.1255 12
2.2505 12
2.3755 12
2.5005 12
2.6255 12
2.7505 12
2.8755 12
3.0005 12
3.1255 12
3.2505 12
3.3755 12
3.5005 12
3.6255 12
3.7505 12
3.8755 12
4.0005 12
4.1255 12
4.2505 12
4.3755 12
4.5005 12
4.6255 12
4.7505 12
4.8755 12
};
\addlegendentry{$\epsilon_u = 10^{-10}$}
\addplot [line width=1.pt,  myred, dashed]
table {%
0.0005 1
0.1255 4
0.2505 5
0.3755 5
0.5005 5
0.6255 5
0.7505 6
0.8755 6
1.0005 7
1.1255 7
1.2505 7
1.3755 7
1.5005 7
1.6255 7
1.7505 7
1.8755 7
2.0005 7
2.1255 7
2.2505 7
2.3755 7
2.5005 7
2.6255 7
2.7505 7
2.8755 7
3.0005 7
3.1255 7
3.2505 7
3.3755 7
3.5005 7
3.6255 7
3.7505 7
3.8755 7
4.0005 7
4.1255 7
4.2505 7
4.3755 7
4.5005 7
4.6255 7
4.7505 7
4.8755 7
};
\addlegendentry{$\epsilon_u = 10^{-7}$}
\addplot [line width=1.pt,  myblue, dash pattern=on 1pt off 3pt on 3pt off 3pt]
table {%
0.0005 1
0.1255 2
0.2505 2
0.3755 2
0.5005 2
0.6255 3
0.7505 3
0.8755 4
1.0005 4
1.1255 4
1.2505 4
1.3755 4
1.5005 4
1.6255 4
1.7505 4
1.8755 4
2.0005 4
2.1255 4
2.2505 4
2.3755 4
2.5005 4
2.6255 4
2.7505 4
2.8755 4
3.0005 4
3.1255 4
3.2505 4
3.3755 4
3.5005 4
3.6255 4
3.7505 4
3.8755 4
4.0005 4
4.1255 4
4.2505 4
4.3755 4
4.5005 4
4.6255 4
4.7505 4
4.8755 4
};
\addlegendentry{$\epsilon_u = 10^{-4}$}
\end{axis}

\end{tikzpicture}
         \label{fig:ADR_r}
     \end{subfigure}
        \caption{Stochastic Advection-Diffusion-Reaction equation: Left: Singular values of TDB-CUR with 100,000 samples and 5 over-sampled points. Right: Rank versus time. The result is presented for different values of the upper error bound ($\epsilon_u$) for mode addition.}
        \label{fig:ADR_100000}
\end{figure}

In the second case, we take $s=100,000$ samples of the random diffusion coefficient. This case is of particular interest, as it demonstrates the true potential of the TDB-CUR in cases where the FOM is too costly to run. In order to execute this scenario with the FOM, not only would we require sufficient memory to store the solution matrix of size $n\times s$, we would have to compute the nonlinear map of this massive matrix at each time step. On the other hand, with the new methodology, we never require storing a  matrix larger than $n\times r$ or $s\times (m+r)$. To demonstrate this capability, we solve the TDB-CUR with $s=100,000$ on a laptop computer. Since the available computational resources (our laptop) did not have sufficient memory to store the FOM solution matrix, we could not solve the FOM for comparison. Instead, we performed a convergence study by decreasing the threshold for rank addition, $\epsilon_u$. By decreasing $\epsilon_u$, we observe two things: (i) the rank is increased more rapidly and (ii) the maximum rank is increased. Figure \ref{fig:ADR_100000} depicts the singular values of the TDB-CUR method versus time (left), and the rank at each time step for different values of $\epsilon_u$ (right). As $\epsilon_u$ is decreased, the rank is increased, and we observe convergence in the leading singular values.

\section{Conclusion}\label{sec:conclusion}
The objective of this work was to develop a method to solve nonlinear matrix differential equations (MDEs) that is accurate, well-conditioned, computationally efficient, and minimally intrusive. To this end, we presented the TDB-CUR
algorithm for solving MDEs via low-rank approximation. The algorithm is based on a time-discrete variational principle that leverages sparse sampling to efficiently compute a low-rank matrix approximation at each iteration of the time-stepping scheme. Numerical experiments illustrate that the TDB-CUR algorithm provides a near-optimal low-rank approximation to the solution of MDEs, while significantly reducing the computational cost. Moreover, we showed the method is robust in the presence of small singular values, and significantly outperforms DLRA based on the time continuous variational principle, \cb{unconventional integrator, and PRK4}. Although not investigated in the present work, the TDB-CUR algorithm is also highly parallelizable, making it an attractive option for high-performance computing tasks.

While the presented approach is minimally intrusive and can be applied to systems containing general nonlinearities, the goal of future work should be to make this method fully non-intrusive, allowing the FOM to be leveraged as a black box. This will allow the method to be applied to proprietary solvers while reducing the overall implementation efforts, making this powerful methodology more accessible to researchers and practitioners, alike. 

\section*{Acknowledgments}
The authors thank Dr. Gianluca Ceruti for numerous insightful and stimulating discussions that led to a number of improvements. This work is supported  by the Air Force Office of Scientific Research award  FA9550-21-1-0247 and funding from Transformational Tools and Technology (TTT), NASA Grant No. 80NSSC22M0282. Computational resources are provided by the Center for Research Computing (CRC) at the University of Pittsburgh.

\clearpage
\appendix
\section{Example Matlab Code}\label{sec:matlab-code}
\begin{lstlisting}[style=Matlab-editor,
frame=single,
numbers=none,
basicstyle=\ttfamily\tiny,
caption=Matlab code to solve the stochastic Burgers equation using TDB-CUR with oversampling,
label={list:matlab}]
close all; clear all; clc
global nu xi d N

rng default
d = 17; sigma = 0.001; % random dimension; standard deviation
Ns = 256; xi = sigma*randn(Ns, d); % # of samples; random parameter samples

dt = 2.5e-4; t0 = 0.0; tf = 5.0; iter_max = round((tf-t0)/dt); save_iter = 50; 
nu = 2.5e-3;
N = 401; xmin = 0; xmax = 1; % # of grid points
x = linspace(xmin, xmax, N)'; dx = x(2) - x(1);
e = ones(N,1);
D1 = spdiags([-e, e]/(2*dx), [-1, 1], N, N); D1([1,end], :) = 0; % d/dx
D2 = spdiags([e, -2*e, e]/dx^2, -1:1, N, N); D2([1,end], :) = 0; % d^2/dx^2

lc = 0.6; K = exp(-(x-x').^2/(2*lc^2)); % squared exponential kernel
[Ux,Lx,~] = svd(K); Ux=Ux(:,1:d); Lx=sqrt(diag(Lx(1:d,1:d)))';
% [Ux,Lx,~] = svds(K,d); Lx=sqrt(diag(Lx))';
ub = sin(2*pi*x) .* (0.5*(exp(cos(2*pi*x))-1.5)); % base flow IC
Ux = sin(2*pi*x) .* (Lx.*Ux);                     % IC perturbations

%% TDB-CUR: initial condition
r = 18; % TDB rank; max=18
m = 5;  % # of rows to oversample
U = [ub, Ux];
Y = [ones(Ns,1), xi];
[U ,R1] = qr(U, 0);
[Y, R2] = qr(Y, 0);
[RU,S,RY] = svd(R1*R2');
U = U*RU; Y = Y*RY;
U=U(:,1:r); S=S(1:r,1:r); Y=Y(:,1:r);

t = t0;
umean = U*S*sum(Y)'/Ns; tt = 0.0; Sig_TDB = diag(S)';
for iter=1:iter_max % time integration loop
    if mod(iter,100)==0; disp(iter); end
    
    %% TDB-CUR: compute row and column indices
    s = gpode(Y, r);    % column indices
    p = gpode(U, r+m); % rows indices
    [~,pa] = find(D2(p,:)); pa = unique(pa);  % adjacent points
    
    %% TDB-CUR: compute rank-r approximation of G
    u_s = U*S*Y(s,:)'; u_p = U(p,:)*S*Y'; u_pa = U(pa,:)*S*Y';
    [F1_s, F1_p, F1_pa] = f_ss(t, u_s, u_pa, s, p, pa, D1, D2);
    [F2_s, F2_p, F2_pa] = f_ss(t+0.5*dt, u_s+0.5*dt*F1_s, u_pa+0.5*dt*F1_pa, s, p, pa, D1, D2);
    [F3_s, F3_p, F3_pa] = f_ss(t+0.5*dt, u_s+0.5*dt*F2_s, u_pa+0.5*dt*F2_pa, s, p, pa, D1, D2);
    [F4_s, F4_p,   ~   ] = f_ss(t+dt, u_s+dt*F3_s, u_pa+dt*F3_pa, s, p, pa, D1, D2);
    G_s = u_s + dt*(F1_s+2.0*F2_s+2.0*F3_s+F4_s)/6;
    G_p = u_p + dt*(F1_p+2.0*F2_p+2.0*F3_p+F4_p)/6;
    [U,~] = qr(G_s,0);
    Y = (inv(U(p,:)'*U(p,:))*U(p,:)'*G_p)';
    [Y,S,RU] = svd(Y,0);
    U = U*RU;

    t = iter*dt;
    if mod(iter, save_iter)==0
        tt = [tt t]; umean = [umean, U*S*sum(Y)'/Ns]; Sig_TDB = [Sig_TDB; diag(S)'];
    end
end
subplot(1,2,1), surf(x,tt,umean'), shading interp, view([0 0 1]), hold on; % mean solution (x,t)
subplot(1,2,2), semilogy(tt,Sig_TDB,'LineWidth',2,'Color','k')             % singular values vs time

function dudt = f(t, u, D1, D2)
global nu
dudt = -0.5*D1*u.^2 + nu*D2*u; % Burgers equation
end

function out = gdot(t, xi)
global d
out = -2*pi*cos(2*pi*t) + (pi*cos(pi*(1:d).*t)./(1:d))*xi';
end

%% TDB-CUR: compute Burgers rhs (f) at specified rows and columns
function [F_s, F_p, F_pa]=f_ss(t, u_s, u_pa, s, p, pa, D1, D2)
global xi N
F_s = f(t, u_s, D1, D2);
F_s(1,:) = gdot(t,xi(s,:)); F_s(end,:) = 0;   % apply boudnary conditions
F_p = f(t, u_pa, D1(p,pa), D2(p,pa));
if ismember(1, p);  F_p(1,:) = gdot(t,xi); end % check left boundary  (x=0)
if ismember(N, p); F_p(end,:) = 0; end         % check right boundary (x=1)
[U_F,~] = qr(F_s,0);
Z_F = inv(U_F(p,:)'*U_F(p,:))*U_F(p,:)'*F_p;
F_pa = U_F(pa,:)*Z_F;
end





\end{lstlisting}

\clearpage
\begin{lstlisting}[style=Matlab-editor,
frame=single,
numbers=none,
basicstyle=\ttfamily\tiny,
caption=GappyPOD+E algorithm  adapted from \cite{peherstorfer2020stability},
label={list:gpode}]
%% TDB-CUR: GappyPOD+E Algorithm
function [ p ] = gpode( U, np )
[~,~,p] = qr(U', 'vector'); % QDEIM (or DEIM)
p = p(1:size(U,2));         % take points equal to number of basis 
for i=length(p)+1:np
    [~, S, W] = svd(U(p, :), 0);
    g = S(end-1, end-1)^2 - S(end, end)^2; 
    Ub = W'*U';
    r = g + sum(Ub.^2, 1);
    r = r-sqrt((g+sum(Ub.^2,1)).^2-4*g*Ub(end, :).^2);
    [~, I] = sort(r, 'descend');
    e = 1;
    while any(I(e) == p)
        e = e + 1;
    end
    p(end + 1) = I(e);
end
p = sort(p);
end
\end{lstlisting}

\bibliographystyle{abbrvnat}
\bibliography{Hessam, MD}

\begin{thebibliography}{56}
\providecommand{\natexlab}[1]{#1}
\providecommand{\url}[1]{\texttt{#1}}
\expandafter\ifx\csname urlstyle\endcsname\relax
  \providecommand{\doi}[1]{doi: #1}\else
  \providecommand{\doi}{doi: \begingroup \urlstyle{rm}\Url}\fi

\bibitem[Amsallem and Farhat(2008)]{AF08}
D.~Amsallem and C.~Farhat.
\newblock Interpolation method for adapting reduced-order models and
  application to aeroelasticity.
\newblock \emph{AIAA Journal}, 46\penalty0 (7):\penalty0 1803--1813, 2023/08/15
  2008.
\newblock \doi{10.2514/1.35374}.
\newblock URL \url{https://doi.org/10.2514/1.35374}.

\bibitem[Anderson et~al.(2015)Anderson, Du, Mahoney, Melgaard, Wu, and
  Gu]{anderson2015spectral}
D.~Anderson, S.~Du, M.~Mahoney, C.~Melgaard, K.~Wu, and M.~Gu.
\newblock Spectral gap error bounds for improving {CUR} matrix decomposition
  and the {N}ystr{\"o}m method.
\newblock In \emph{Artificial Intelligence and Statistics}, pages 19--27. PMLR,
  2015.

\bibitem[Babaee(2019)]{B19}
H.~Babaee.
\newblock An observation-driven time-dependent basis for a reduced description
  of transient stochastic systems.
\newblock \emph{Proceedings of the Royal Society A: Mathematical, Physical and
  Engineering Sciences}, 475\penalty0 (2231):\penalty0 20190506, 2019.
\newblock \doi{10.1098/rspa.2019.0506}.
\newblock URL \url{https://doi.org/10.1098/rspa.2019.0506}.

\bibitem[Babaee and Sapsis(2016)]{Babaee_PRSA}
H.~Babaee and T.~P. Sapsis.
\newblock A minimization principle for the description of modes associated with
  finite-time instabilities.
\newblock \emph{Proceedings of the Royal Society of London A: Mathematical,
  Physical and Engineering Sciences}, 472\penalty0 (2186):\penalty0 20150779,
  2016.
\newblock URL \url{http://dx.doi.org/10.1098/rspa.2015.0779}.

\bibitem[Babaee et~al.(2013)Babaee, Wan, and Acharya]{Babaee:2013aa}
H.~Babaee, X.~Wan, and S.~Acharya.
\newblock Effect of uncertainty in blowing ratio on film cooling effectiveness.
\newblock \emph{Journal of Heat Transfer}, 136\penalty0 (3):\penalty0
  031701--031701, 11 2013.
\newblock URL \url{http://dx.doi.org/10.1115/1.4025562}.

\bibitem[Babaee et~al.(2017)Babaee, Choi, Sapsis, and
  Karniadakis]{Babaee:2017aa}
H.~Babaee, M.~Choi, T.~P. Sapsis, and G.~E. Karniadakis.
\newblock A robust bi-orthogonal/dynamically-orthogonal method using the
  covariance pseudo-inverse with application to stochastic flow problems.
\newblock \emph{Journal of Computational Physics}, 344:\penalty0 303--319, 9
  2017.
\newblock \doi{https://doi.org/10.1016/j.jcp.2017.04.057}.
\newblock URL
  \url{http://www.sciencedirect.com/science/article/pii/S0021999117303364}.

\bibitem[Barrault et~al.(2004)Barrault, Maday, Nguyen, and Patera]{BMNP04}
M.~Barrault, Y.~Maday, N.~C. Nguyen, and A.~T. Patera.
\newblock An `empirical interpolation' method: application to efficient
  reduced-basis discretization of partial differential equations.
\newblock \emph{Comptes Rendus Mathematique}, 339\penalty0 (9):\penalty0
  667--672, 2004.
\newblock \doi{https://doi.org/10.1016/j.crma.2004.08.006}.
\newblock URL
  \url{https://www.sciencedirect.com/science/article/pii/S1631073X04004248}.

\bibitem[Barth et~al.(2011)Barth, Schwab, and Zollinger]{barth2011multi}
A.~Barth, C.~Schwab, and N.~Zollinger.
\newblock Multi-level {M}onte {C}arlo finite element method for elliptic {PDE}s
  with stochastic coefficients.
\newblock \emph{Numerische Mathematik}, 119:\penalty0 123--161, 2011.

\bibitem[Beck et~al.(2000)Beck, J{\"a}ckle, Worth, and Meyer]{Beck:2000aa}
M.~H. Beck, A.~J{\"a}ckle, G.~A. Worth, and H.~D. Meyer.
\newblock The multiconfiguration time-dependent {H}artree ({MCTDH}) method: a
  highly efficient algorithm for propagating wavepackets.
\newblock \emph{Physics Reports}, 324\penalty0 (1):\penalty0 1--105, 1 2000.
\newblock \doi{http://dx.doi.org/10.1016/S0370-1573(99)00047-2}.
\newblock URL
  \url{http://www.sciencedirect.com/science/article/pii/S0370157399000472}.

\bibitem[Blanchard and Sapsis(2019{\natexlab{a}})]{BS19}
A.~Blanchard and T.~P. Sapsis.
\newblock Analytical description of optimally time-dependent modes for
  reduced-order modeling of transient instabilities.
\newblock \emph{SIAM Journal on Applied Dynamical Systems}, 18\penalty0
  (2):\penalty0 1143--1162, 2019{\natexlab{a}}.

\bibitem[Blanchard and Sapsis(2019{\natexlab{b}})]{blanchard2019learning}
A.~Blanchard and T.~P. Sapsis.
\newblock Learning the tangent space of dynamical instabilities from data.
\newblock \emph{Chaos: An Interdisciplinary Journal of Nonlinear Science},
  29\penalty0 (11), 2019{\natexlab{b}}.

\bibitem[Ceruti and Lubich(2021{\natexlab{a}})]{CL21}
G.~Ceruti and C.~Lubich.
\newblock An unconventional robust integrator for dynamical low-rank
  approximation.
\newblock \emph{BIT Numerical Mathematics}, 2021{\natexlab{a}}.
\newblock \doi{10.1007/s10543-021-00873-0}.
\newblock URL \url{https://doi.org/10.1007/s10543-021-00873-0}.

\bibitem[Ceruti and Lubich(2021{\natexlab{b}})]{ceruti2021unconventional}
G.~Ceruti and C.~Lubich.
\newblock An unconventional robust integrator for dynamical low-rank
  approximation.
\newblock \emph{BIT Numerical Mathematics}, pages 1--22, 2021{\natexlab{b}}.

\bibitem[Ceruti et~al.(2021)Ceruti, Kusch, and Lubich]{lubich2021rankadaptive}
G.~Ceruti, J.~Kusch, and C.~Lubich.
\newblock A rank-adaptive robust integrator for dynamical low-rank
  approximation.
\newblock \emph{arXiv preprint arXiv:2104.05247}, 2021.

\bibitem[Charous and Lermusiaux(2023)]{CL23}
A.~Charous and P.~F. Lermusiaux.
\newblock Dynamically orthogonal runge--kutta schemes with perturbative
  retractions for the dynamical low-rank approximation.
\newblock \emph{SIAM Journal on Scientific Computing}, 45\penalty0
  (2):\penalty0 A872--A897, 2023.

\bibitem[Chaturantabut and Sorensen(2010{\natexlab{a}})]{CS10}
S.~Chaturantabut and D.~C. Sorensen.
\newblock Nonlinear model reduction via discrete empirical interpolation.
\newblock \emph{SIAM Journal on Scientific Computing}, 32\penalty0
  (5):\penalty0 2737--2764, 2020/12/11 2010{\natexlab{a}}.
\newblock \doi{10.1137/090766498}.
\newblock URL \url{https://doi.org/10.1137/090766498}.

\bibitem[Chaturantabut and
  Sorensen(2010{\natexlab{b}})]{chaturantabut2010nonlinear}
S.~Chaturantabut and D.~C. Sorensen.
\newblock Nonlinear model reduction via discrete empirical interpolation.
\newblock \emph{SIAM Journal on Scientific Computing}, 32\penalty0
  (5):\penalty0 2737--2764, 2010{\natexlab{b}}.

\bibitem[Cheng et~al.(2013)Cheng, Hou, and Zhang]{CHZI13}
M.~Cheng, T.~Y. Hou, and Z.~Zhang.
\newblock A dynamically bi-orthogonal method for time-dependent stochastic
  partial differential equations i: Derivation and algorithms.
\newblock \emph{Journal of Computational Physics}, 242\penalty0 (0):\penalty0
  843 -- 868, 2013.
\newblock ISSN 0021-9991.
\newblock \doi{http://dx.doi.org/10.1016/j.jcp.2013.02.033}.
\newblock URL
  \url{http://www.sciencedirect.com/science/article/pii/S0021999113001526}.

\bibitem[Choi et~al.(2014)Choi, Sapsis, and Karniadakis]{CSK14}
M.~Choi, T.~P. Sapsis, and G.~E. Karniadakis.
\newblock On the equivalence of dynamically orthogonal and bi-orthogonal
  methods: Theory and numerical simulations.
\newblock \emph{Journal of Computational Physics}, 270:\penalty0 1 -- 20, 2014.
\newblock ISSN 0021-9991.
\newblock \doi{http://dx.doi.org/10.1016/j.jcp.2014.03.050}.
\newblock URL
  \url{http://www.sciencedirect.com/science/article/pii/S002199911400237X}.

\bibitem[Dektor et~al.(2021)Dektor, Rodgers, and Venturi]{dektor2021rank}
A.~Dektor, A.~Rodgers, and D.~Venturi.
\newblock Rank-adaptive tensor methods for high-dimensional nonlinear pdes.
\newblock \emph{Journal of Scientific Computing}, 88\penalty0 (2):\penalty0
  1--27, 2021.

\bibitem[Dieci and Elia(2006)]{Dieci:2006aa}
L.~Dieci and C.~Elia.
\newblock The singular value decomposition to approximate spectra of dynamical
  systems. theoretical aspects.
\newblock \emph{Journal of Differential Equations}, 230\penalty0 (2):\penalty0
  502--531, 2006.
\newblock \doi{http://dx.doi.org/10.1016/j.jde.2006.08.007}.
\newblock URL
  \url{http://www.sciencedirect.com/science/article/pii/S0022039606003263}.

\bibitem[Donello et~al.(2022)Donello, Carpenter, and Babaee]{DCB22}
M.~Donello, M.~H. Carpenter, and H.~Babaee.
\newblock Computing sensitivities in evolutionary systems: A real-time reduced
  order modeling strategy.
\newblock \emph{SIAM Journal on Scientific Computing}, pages A128--A149,
  2022/01/19 2022.
\newblock \doi{10.1137/20M1388565}.
\newblock URL \url{https://doi.org/10.1137/20M1388565}.

\bibitem[Drma{\v c} and Gugercin(2016)]{DG16}
Z.~Drma{\v c} and S.~Gugercin.
\newblock A new selection operator for the discrete empirical interpolation
  method---improved a priori error bound and extensions.
\newblock \emph{SIAM Journal on Scientific Computing}, 38\penalty0
  (2):\penalty0 A631--A648, 2016.
\newblock \doi{10.1137/15M1019271}.
\newblock URL \url{https://doi.org/10.1137/15M1019271}.

\bibitem[Einkemmer and Lubich(2018)]{EL18}
L.~Einkemmer and C.~Lubich.
\newblock A low-rank projector-splitting integrator for the vlasov--poisson
  equation.
\newblock \emph{SIAM Journal on Scientific Computing}, 40\penalty0
  (5):\penalty0 B1330--B1360, 2023/08/15 2018.
\newblock \doi{10.1137/18M116383X}.
\newblock URL \url{https://doi.org/10.1137/18M116383X}.

\bibitem[Farhat et~al.(2014)Farhat, Avery, Chapman, and Cortial]{FACC14}
C.~Farhat, P.~Avery, T.~Chapman, and J.~Cortial.
\newblock Dimensional reduction of nonlinear finite element dynamic models with
  finite rotations and energy-based mesh sampling and weighting for
  computational efficiency.
\newblock \emph{International Journal for Numerical Methods in Engineering},
  98\penalty0 (9):\penalty0 625--662, 2023/08/15 2014.
\newblock \doi{https://doi.org/10.1002/nme.4668}.
\newblock URL \url{https://doi.org/10.1002/nme.4668}.

\bibitem[Giles(2008)]{giles2008multilevel}
M.~B. Giles.
\newblock Multilevel {M}onte {C}arlo path simulation.
\newblock \emph{Operations research}, 56\penalty0 (3):\penalty0 607--617, 2008.

\bibitem[Halko et~al.(2011)Halko, Martinsson, and Tropp]{halko2011finding}
N.~Halko, P.-G. Martinsson, and J.~A. Tropp.
\newblock Finding structure with randomness: {P}robabilistic algorithms for
  constructing approximate matrix decompositions.
\newblock \emph{SIAM review}, 53\penalty0 (2):\penalty0 217--288, 2011.

\bibitem[Hu and Wang(2022)]{HW22}
J.~Hu and Y.~Wang.
\newblock An adaptive dynamical low rank method for the nonlinear boltzmann
  equation.
\newblock \emph{Journal of Scientific Computing}, 92\penalty0 (2):\penalty0 75,
  2022.
\newblock \doi{10.1007/s10915-022-01934-4}.
\newblock URL \url{https://doi.org/10.1007/s10915-022-01934-4}.

\bibitem[Karniadakis and Sherwin(2005)]{KS05}
G.~E. Karniadakis and S.~J. Sherwin.
\newblock \emph{Spectral/hp element methods for computational fluid dynamics}.
\newblock Oxford University Press, USA, 2005.

\bibitem[Kieri and Vandereycken(2019)]{kieri2019projection}
E.~Kieri and B.~Vandereycken.
\newblock Projection methods for dynamical low-rank approximation of
  high-dimensional problems.
\newblock \emph{Computational Methods in Applied Mathematics}, 19\penalty0
  (1):\penalty0 73--92, 2019.

\bibitem[Koch and Lubich(2007)]{KL07}
O.~Koch and C.~Lubich.
\newblock Dynamical low‐rank approximation.
\newblock \emph{SIAM Journal on Matrix Analysis and Applications}, 29\penalty0
  (2):\penalty0 434--454, 2017/04/02 2007.
\newblock \doi{10.1137/050639703}.
\newblock URL \url{http://dx.doi.org/10.1137/050639703}.

\bibitem[Kuo et~al.(2012)Kuo, Schwab, and Sloan]{kuo2012quasi}
F.~Y. Kuo, C.~Schwab, and I.~H. Sloan.
\newblock Quasi-{M}onte {C}arlo finite element methods for a class of elliptic
  partial differential equations with random coefficients.
\newblock \emph{SIAM Journal on Numerical Analysis}, 50\penalty0 (6):\penalty0
  3351--3374, 2012.

\bibitem[{Kusch, J.} and {Stammer, P.}(2023)]{KS23}
{Kusch, J.} and {Stammer, P.}
\newblock A robust collision source method for rank adaptive dynamical low-rank
  approximation in radiation therapy.
\newblock \emph{ESAIM: M2AN}, 57\penalty0 (2):\penalty0 865--891, 2023.
\newblock \doi{10.1051/m2an/2022090}.
\newblock URL \url{https://doi.org/10.1051/m2an/2022090}.

\bibitem[Lubich and Oseledets(2014)]{LO14}
C.~Lubich and I.~V. Oseledets.
\newblock A projector-splitting integrator for dynamical low-rank
  approximation.
\newblock \emph{BIT Numerical Mathematics}, 54\penalty0 (1):\penalty0 171--188,
  2014.
\newblock \doi{10.1007/s10543-013-0454-0}.
\newblock URL \url{http://dx.doi.org/10.1007/s10543-013-0454-0}.

\bibitem[Mahoney and Drineas(2009)]{mahoney2009cur}
M.~W. Mahoney and P.~Drineas.
\newblock {CUR} matrix decompositions for improved data analysis.
\newblock \emph{Proceedings of the National Academy of Sciences}, 106\penalty0
  (3):\penalty0 697--702, 2009.

\bibitem[Manohar et~al.(2018)Manohar, Brunton, Kutz, and Brunton]{MBKB18}
K.~Manohar, B.~W. Brunton, J.~N. Kutz, and S.~L. Brunton.
\newblock Data-driven sparse sensor placement for reconstruction: Demonstrating
  the benefits of exploiting known patterns.
\newblock \emph{IEEE Control Systems Magazine}, 38\penalty0 (3):\penalty0
  63--86, 2018.
\newblock \doi{10.1109/MCS.2018.2810460}.

\bibitem[Musharbash and Nobile(2018)]{MN18}
E.~Musharbash and F.~Nobile.
\newblock Dual dynamically orthogonal approximation of incompressible {N]avier
  [S}tokes equations with random boundary conditions.
\newblock \emph{Journal of Computational Physics}, 354:\penalty0 135--162,
  2018.
\newblock \doi{https://doi.org/10.1016/j.jcp.2017.09.061}.
\newblock URL
  \url{http://www.sciencedirect.com/science/article/pii/S0021999117307349}.

\bibitem[Naderi and Babaee(2023)]{NB23}
M.~H. Naderi and H.~Babaee.
\newblock Adaptive sparse interpolation for accelerating nonlinear stochastic
  reduced-order modeling with time-dependent bases.
\newblock \emph{Computer Methods in Applied Mechanics and Engineering},
  405:\penalty0 115813, 2023.
\newblock \doi{https://doi.org/10.1016/j.cma.2022.115813}.
\newblock URL
  \url{https://www.sciencedirect.com/science/article/pii/S0045782522007691}.

\bibitem[Nouri et~al.(2021)Nouri, Babaee, Givi, Chelliah, and Livescu]{NBGCL21}
A.~G. Nouri, H.~Babaee, P.~Givi, H.~K. Chelliah, and D.~Livescu.
\newblock Skeletal model reduction with forced optimally time dependent modes.
\newblock \emph{Combustion and Flame}, page 111684, 2021.
\newblock \doi{https://doi.org/10.1016/j.combustflame.2021.111684}.
\newblock URL
  \url{https://www.sciencedirect.com/science/article/pii/S0010218021004272}.

\bibitem[Patil and Babaee(2020)]{PB20}
P.~Patil and H.~Babaee.
\newblock Real-time reduced-order modeling of stochastic partial differential
  equations via time-dependent subspaces.
\newblock \emph{Journal of Computational Physics}, 415:\penalty0 109511, 2020.
\newblock \doi{https://doi.org/10.1016/j.jcp.2020.109511}.
\newblock URL
  \url{http://www.sciencedirect.com/science/article/pii/S0021999120302850}.

\bibitem[Patil and Babaee(2023)]{patil2023reduced}
P.~Patil and H.~Babaee.
\newblock Reduced-order modeling with time-dependent bases for pdes with
  stochastic boundary conditions.
\newblock \emph{SIAM/ASA Journal on Uncertainty Quantification}, 11\penalty0
  (3):\penalty0 727--756, 2023.

\bibitem[Peherstorfer(2020)]{doi:10.1137/19M1257275}
B.~Peherstorfer.
\newblock Model reduction for transport-dominated problems via online adaptive
  bases and adaptive sampling.
\newblock \emph{SIAM Journal on Scientific Computing}, 42\penalty0
  (5):\penalty0 A2803--A2836, 2020.
\newblock \doi{10.1137/19M1257275}.
\newblock URL \url{https://doi.org/10.1137/19M1257275}.

\bibitem[Peherstorfer and Willcox(2015)]{doi:10.1137/140989169}
B.~Peherstorfer and K.~Willcox.
\newblock Online adaptive model reduction for nonlinear systems via low-rank
  updates.
\newblock \emph{SIAM Journal on Scientific Computing}, 37\penalty0
  (4):\penalty0 A2123--A2150, 2015.
\newblock \doi{10.1137/140989169}.
\newblock URL \url{https://doi.org/10.1137/140989169}.

\bibitem[Peherstorfer et~al.(2020)Peherstorfer, Drmac, and
  Gugercin]{peherstorfer2020stability}
B.~Peherstorfer, Z.~Drmac, and S.~Gugercin.
\newblock Stability of discrete empirical interpolation and gappy proper
  orthogonal decomposition with randomized and deterministic sampling points.
\newblock \emph{SIAM Journal on Scientific Computing}, 42\penalty0
  (5):\penalty0 A2837--A2864, 2020.

\bibitem[Ramezanian et~al.(2021)Ramezanian, Nouri, and Babaee]{RNB21}
D.~Ramezanian, A.~G. Nouri, and H.~Babaee.
\newblock On-the-fly reduced order modeling of passive and reactive species via
  time-dependent manifolds.
\newblock \emph{Computer Methods in Applied Mechanics and Engineering},
  382:\penalty0 113882, 2021.
\newblock \doi{https://doi.org/10.1016/j.cma.2021.113882}.
\newblock URL
  \url{https://www.sciencedirect.com/science/article/pii/S004578252100219X}.

\bibitem[Rodgers et~al.(2022)Rodgers, Dektor, and Venturi]{rodgers2022adaptive}
A.~Rodgers, A.~Dektor, and D.~Venturi.
\newblock Adaptive integration of nonlinear evolution equations on tensor
  manifolds.
\newblock \emph{Journal of Scientific Computing}, 92\penalty0 (2):\penalty0 39,
  2022.

\bibitem[Ryckelynck(2005)]{R05}
D.~Ryckelynck.
\newblock A priori hyperreduction method: an adaptive approach.
\newblock \emph{Journal of Computational Physics}, 202\penalty0 (1):\penalty0
  346--366, 2005.
\newblock \doi{https://doi.org/10.1016/j.jcp.2004.07.015}.
\newblock URL
  \url{https://www.sciencedirect.com/science/article/pii/S002199910400289X}.

\bibitem[Sapsis and Lermusiaux(2009)]{SL09}
T.~Sapsis and P.~Lermusiaux.
\newblock Dynamically orthogonal field equations for continuous stochastic
  dynamical systems.
\newblock \emph{Physica D: Nonlinear Phenomena}, 238\penalty0 (23-24):\penalty0
  2347--2360, 2009.

\bibitem[Schotth{\"o}fer et~al.(2022)Schotth{\"o}fer, Zangrando, Kusch, Ceruti,
  and Tudisco]{SZK23}
S.~Schotth{\"o}fer, E.~Zangrando, J.~Kusch, G.~Ceruti, and F.~Tudisco.
\newblock Low-rank lottery tickets: finding efficient low-rank neural networks
  via matrix differential equations, 2022.

\bibitem[Sorensen and Embree(2016)]{sorensen2016deim}
D.~C. Sorensen and M.~Embree.
\newblock A {DEIM} induced {CUR} factorization.
\newblock \emph{SIAM Journal on Scientific Computing}, 38\penalty0
  (3):\penalty0 A1454--A1482, 2016.

\bibitem[Szyld(2006)]{S06}
D.~B. Szyld.
\newblock The many proofs of an identity on the norm of oblique projections.
\newblock \emph{Numerical Algorithms}, 42\penalty0 (3):\penalty0 309--323,
  2006.
\newblock \doi{10.1007/s11075-006-9046-2}.
\newblock URL \url{https://doi.org/10.1007/s11075-006-9046-2}.

\bibitem[Vidal et~al.(2005)Vidal, Ma, and Sastry]{vidal2005generalized}
R.~Vidal, Y.~Ma, and S.~Sastry.
\newblock Generalized principal component analysis ({GPCA}).
\newblock \emph{IEEE transactions on pattern analysis and machine
  intelligence}, 27\penalty0 (12):\penalty0 1945--1959, 2005.

\bibitem[Wright(1992)]{Wright:1992aa}
K.~Wright.
\newblock Differential equations for the analytic singular value decomposition
  of a matrix.
\newblock \emph{Numerische Mathematik}, 63\penalty0 (1):\penalty0 283--295,
  1992.
\newblock \doi{10.1007/BF01385862}.
\newblock URL \url{http://dx.doi.org/10.1007/BF01385862}.

\bibitem[Xiu and Hesthaven(2006)]{xiu2006high}
D.~Xiu and J.~Hesthaven.
\newblock High-order collocation methods for differential equations with random
  inputs.
\newblock \emph{SIAM Journal on Scientific Computing}, 27\penalty0
  (3):\penalty0 1118, 2006.

\bibitem[Yang and White(2020)]{yang2020time}
M.~Yang and S.~R. White.
\newblock Time-dependent variational principle with ancillary krylov subspace.
\newblock \emph{Physical Review B}, 102\penalty0 (9):\penalty0 094315, 2020.

\bibitem[Zimmermann and Willcox(2016)]{ZW16}
R.~Zimmermann and K.~Willcox.
\newblock An accelerated greedy missing point estimation procedure.
\newblock \emph{SIAM Journal on Scientific Computing}, 38\penalty0
  (5):\penalty0 A2827--A2850, 2023/08/15 2016.
\newblock \doi{10.1137/15M1042899}.
\newblock URL \url{https://doi.org/10.1137/15M1042899}.

\end{thebibliography}

\end{document}